\documentclass[preprint,3p]{elsarticle}

\usepackage{amssymb}
\usepackage{amsmath}
\usepackage[numbers]{natbib}
\usepackage{caption}
\usepackage{graphicx}
\usepackage{float} 
\usepackage{subfigure}
\usepackage{subcaption}
\usepackage{booktabs}

\newtheorem{lemma}{Lemma}
\newdefinition{defn}{Definition}
\newtheorem{theorem}{Theorem}
\newproof{pf}{Proof}
\newdefinition{rmk}{Remark}
\newtheorem{cor}{Corollary}
\usepackage[colorlinks=true,
linkcolor=blue, %%%修改此处为你想要的颜色
anchorcolor=black, %%修改此处为你想要的颜色
citecolor=blue, %%设置文中引用参考文献的颜色
urlcolor=blue%%设置DOI的颜色为蓝色
]{hyperref}
\allowdisplaybreaks

\journal{arXiv}

\begin{document}

\begin{frontmatter}

%% Title, authors and addresses

%% use the tnoteref command within \title for footnotes;
%% use the tnotetext command for theassociated footnote;
%% use the fnref command within \author or \affiliation for footnotes;
%% use the fntext command for theassociated footnote;
%% use the corref command within \author for corresponding author footnotes;
%% use the cortext command for theassociated footnote;
%% use the ead command for the email address,
%% and the form \ead[url] for the home page:
%% \title{Title\tnoteref{label1}}
%% \tnotetext[label1]{}
%% \author{Name\corref{cor1}\fnref{label2}}
%% \ead{email address}
%% \ead[url]{home page}
%% \fntext[label2]{}
%% \cortext[cor1]{}
%% \affiliation{organization={},
%%             addressline={},
%%             city={},
%%             postcode={},
%%             state={},
%%             country={}}
%% \fntext[label3]{}

\title{On a class of forward-backward reaction-diffusion systems with local and nonlocal coupling for image restoration} %% Article title

%% use optional labels to link authors explicitly to addresses:
%% \author[label1,label2]{}
%% \affiliation[label1]{organization={},
%%             addressline={},
%%             city={},
%%             postcode={},
%%             state={},
%%             country={}}
%%
%% \affiliation[label2]{organization={},
%%             addressline={},
%%             city={},
%%             postcode={},
%%             state={},
%%             country={}}

\author[1]{Yihui Tong} %% Author name
\ead{23B912019@stu.hit.edu.cn}
%% Author affiliation

\author[1]{Wenjie Liu} %% Author name
\ead{liuwenjie@hit.edu.cn}

\author[1]{Zhichang Guo}
\ead{mathgzc@hit.edu.cn}

\author[2]{Jingfeng Shao}
\ead{sjfmath@gxu.edu.cn}

\author[1]{Wenjuan Yao\corref{cor1}}
\ead{mathywj@hit.edu.cn}
\cortext[cor1]{Corresponding author}

%% Author affiliation
\affiliation[1]{organization={School of Mathematics},%Department and Organization
	addressline={Harbin Institute of Technology}, 
	city={Harbin},
	postcode={150001}, 
	%state={},
	country={China}}
	
\affiliation[2]{organization={School of Mathematics and Information Science},
	addressline={Guangxi University}, 
	city={Nanning},
	%          citysep={}, % Uncomment if no comma needed between city and postcode
	postcode={530004}, 
	%state={},
	country={China}}

%% Abstract
\begin{abstract}
%% Text of abstract
This paper investigates a class of novel nonlinear reaction-diffusion systems that couple forward-backward with fractional diffusion for image restoration, offering the advantage of preserving both contour features and textures. The existence of Young measure solutions to the proposed model is established using the regularization technique, Rothe's method, relaxation theorem, and Moser's iteration. Uniqueness follows from the independence property satisfied by the solution. Numerical experiments illustrate the effectiveness of our model in image denoising and deblurring, in comparison with existing methods.
\end{abstract}

%%Graphical abstract
%\begin{graphicalabstract}
%\includegraphics{grabs}
%\end{graphicalabstract}

%%Research highlights
%\begin{highlights}
%\item Research highlight 1
%\item Research highlight 2
%\end{highlights}

%% Keywords
\begin{keyword}
%% keywords here, in the form: keyword \sep keyword
Reaction-diffusion system \sep Forward-backward diffusion \sep Fractional $p$-Laplacian \sep Young measure solutions \sep Image restoration

%% PACS codes here, in the form: \PACS code \sep code

%% MSC codes here, in the form: \MSC code \sep code
%% or \MSC[2008] code \sep code (2000 is the default)

\end{keyword}

\end{frontmatter}

%% Add \usepackage{lineno} before \begin{document} and uncomment 
%% following line to enable line numbers
%% \linenumbers

%% main text
%%

\section{Introduction}\label{}
\subsection{Related works}
Image deblurring is a fundamental problem in image processing, widely applied in medical imaging, computer vision, and remote sensing. Noise and blur, which frequently arise during image acquisition and transmission, can significantly degrade image quality \cite{WELK2005}. The degradation process is generally modeled as
\begin{equation}
	f=Ku+n,\notag
\end{equation}
where $K$ is the blur operator, and $n$ denotes noise. In this paper, we focus on the non-blind image restoration problem, aiming to recover the clean image $u$ from a noisy, blurred observation $f$, assuming $K$ is known exactly. This typical ill-posed inverse problem can be addressed by minimizing the energy functional
\begin{equation}\label{Vese-formula}
	E_{\lambda}(u)=\int_{\Omega}G\left(\left|\nabla u\right|\right) dx + \frac{\lambda}{2}\int_{\Omega}\left(Ku-f\right)^2 dx,
\end{equation}
where the first term on the right-hand side is the regularization term that enforces the smoothness of the deconvolved image, the second term is the fidelity term, and $\lambda\geq 0$ is a weight parameter. Tikhonov et al. \cite{Tikhonov1977}
introduced the regulariser $G(s)=s^2$, which results in over-smoothing of the deblurred image. To preserve image edges, total variation (TV) regularization was proposed in \cite{Gilles1997,Osher2003}, where $G\left(\left|\nabla u\right|\right)dx=\left|Du\right|$. However, TV-based variational models often yield piecewise constant solutions, leading to the well-known staircase effect \cite{XU2019}. To mitigate this, the adaptive regularizer $G(s)=s^{p(x)}$ was developed in \cite{CHEN2006,SONG2003}. Additionally, various improved models incorporating different regularization terms have been proposed, including sparse regularization models \cite{Carlavan2012,LIU2015}, higher-order TV-based methods \cite{ADAM2020,LV2013}, and total generalized variation (TGV) approaches \cite{Bredies2010,GAO2023}. 

There are two primary approaches to implementing the variational model \eqref{Vese-formula}. One relies on optimization algorithms \cite{Huang2008,Osher2005Iterative,REN2020}, and the other transforms the minimization problem into solving a partial differential equation (PDE) numerically, specifically by considering the gradient descent flow associated with \eqref{Vese-formula}
\begin{equation}\label{GDFlow-Vese}
	\frac{\partial u}{\partial t}=\operatorname{div}\left(g\left(\left|\nabla u\right|\right)\nabla u\right)-\lambda K^{\prime}\left(Ku-f\right),
\end{equation}
where $G^{\prime}(s)=g(s)s$, and $K^{\prime}$ denotes the adjoint operator of $K$. When $\lambda=0$, \eqref{GDFlow-Vese} reduces to a pure diffusion equation, which has been widely used in image denoising. A popular choice of the diffusivity is $g(s)=1/s$, resulting in the TV flow \cite{RUDIN1992}. With the choice $g(s) = \frac{1}{1 + s^2}$ for the diffusivity, \eqref{GDFlow-Vese} is known as the Perona-Malik (PM) equation \cite{Perona&Malik1990}. The forward-backward diffusion characteristic of the PM equation enhances image edges effectively, but the presence of backward diffusion renders the equation ill-posed, leading to staircase instability in numerical discretization \cite{Kichenassamy1996}. To reduce the noise sensitivity of PM diffusion and improve its mathematical well-posedness, various regularization methods \cite{CATTE1992,GUIDOTTI2009nonlocal,GUIDOTTI2012backward,GUIDOTTI2013restoration,SHAO2022,Weickert1998} have been proposed. Fourth-order PDEs \cite{Hajiaboli2011,YK2000} can be utilized to suppress staircase artifacts along ramp edges but may promote piecewise planar solutions. For further advancements, we refer the reader to \cite{GUIDOTTI2009two,JIDESH2012,LLT2003,WEN2022}. In practical applications, the diffusion coefficient in \eqref{GDFlow-Vese} can be tailored to specific image processing requirements and is not necessarily derived from evolution equation associated with a variational model. Zhou et al. \cite{ZHOU2015} introduced the following general reaction-diffusion framework for image restorationstoration, initially applied to multiplicative noise removal
\begin{equation}
	\frac{\partial u}{\partial t}=\operatorname{div}\left(a\left(\left|\nabla u\right|,u\right)\nabla u\right)-\lambda h\left(f,u\right).\notag
\end{equation}
Their framework reduces image recovery tasks to designing suitable forms of the diffusion coefficient $a(|\nabla u|, u)$ and the source term $h(f, u)$. Considering the sharpening effect induced by backward diffusion, Guidotti et al. \cite{GUIDOTTI2013restoration} proposed a mild forward-backward regularization of the PM equation (FBRPM) for image restoration, which is analytically tractable and improves robustness to noise and blur, expressed as
\begin{equation}\label{Guidotti-eqn}
	\frac{\partial u}{\partial t}=\operatorname{div}\left[\left(\frac{1}{1+\left|\nabla u\right|^2}+\delta\left|\nabla u\right|^{\gamma-2}\right)\nabla u\right]-\lambda K^{\prime}\left(Ku-f\right),
\end{equation}
where $\delta>0$, $1<\gamma\leq 2$.
For $\delta<\frac{1}{8}$, the equation remains of the forward-backward type. Unlike the PM equation, the backward regime of \eqref{Guidotti-eqn} is confined to a bounded domain, whose size expands as $\delta$ vanishes. From a dynamical perspective, the regularization replaces the staircasing phenomenon observed in the PM equation with a milder, microstructured ramping characterized by alternating gradients of finite size \cite{Guidotti2015131}. The tunable $\delta$ and $\gamma$ allow \eqref{Guidotti-eqn} to control the locations of forward and backward regimes, facilitating the desired sharpening and smoothing effects. However, previous theoretical research has primarily focused on the case $\lambda = 0$ with sufficiently small $\delta$, under which \eqref{Guidotti-eqn} has no classical weak solutions. For $\gamma = 2$, Demoulini \cite{Demoulini1996} established the existence and uniqueness of the Young measure solution to problem \eqref{Guidotti-eqn} under homogeneous Dirichlet boundary conditions by employing Rothe's method, the relaxation theorem, and the Young measure representation introduced by Kinderlehrer et al. \cite{Kichenassamy1996,Kinderlehrer1992}. The analysis in \cite{Demoulini1996} were further extended to Neumann boundary conditions in \cite{GUIDOTTI2012backward}. The qualitative properties of the Young measure solution to \eqref{Guidotti-eqn} for $1\leq \gamma<2$ were subsequently studied by Yin and Wang \cite{YIN2003Youngmeasure}. For further theoretical studies on regularization schemes for the PM equation, we refer the reader to \cite{Amann2007015,Bellettini2006783,Giacomelli2019374,THANH20141403}.

%This effectively reduces the staircase effect, replacing it with a rather milder micro-structured ramping featuring alternate gradients of finite size 

%Subsequently, Yin et al. \cite{YIN2003Youngmeasure} studied the qualitative properties of the Young measure solution to (3) for $1\leq \gamma<2$.

%Although \eqref{Guidotti-eqn} admits no classical solution, its well-posedness in the sense of Young measure solutions for $\gamma=2$ and $\lambda=0$ has been established in \cite{GUIDOTTI2012backward}, where the boundary conditions are extended from the Dirichlet setting \cite{Demoulini1996} to the Neumann setting.

%In their framework, the image recovery task is reduced to designing an appropriate diffusion coefficient $a\left(\left|\nabla u\right|,u\right)$ and source term $h\left(f,u\right)$.

As an extension of the traditional diffusion equation, the reaction-diffusion system offers greater flexibility in modulating the diffusion process, thus better accommodating the requirements of image restoration. This class of methods originates from the following variational model
\begin{equation}\label{System variation}
	R_\lambda(u)=\varPhi(u)+\frac{\lambda}{2}\left\|Ku-f\right\|_{H^{-1}(\Omega)}^2.
\end{equation}
By selecting the regularization term $\varPhi(u)=\text{TV}(u)$, the authors of \cite{OSHER2003decomposition} proposed a model for image deblurring and cartoon-texture decomposition that combines TV minimization from \cite{Osher2003} with the $H^{-1}$ norm to capture oscillatory components. Guo et al. \cite{GUO2011reaction} introduced the auxiliary variable $w=\Delta^{-1}(f-u)$ into the gradient descent flow of \eqref{System variation}, and derived the following reaction-diffusion system for removing additive noise in images (i.e., the case $K=I$)
\begin{equation}\label{Guosystem2011}
	\left\{\begin{array}{ll}  
		\displaystyle\frac{\partial u}{\partial t}=\operatorname{div}\left(\frac{D u}{\left|D u\right|}\right)-\lambda w, & x\in \Omega,~0<t<T,\\
		\displaystyle\frac{\partial w}{\partial t} = \Delta w - \left(f-u\right), & x\in \Omega,~0<t<T,\\
		\displaystyle u(0,x)=f,~w(0,x)=0, & x\in\Omega,\\
		\displaystyle \frac{\partial u}{\partial \vec{n}}(t,x)=0,~\frac{\partial w}{\partial \vec{n}}(t,x)=0, & x\in \partial\Omega,~0<t<T,
	\end{array} \right.
\end{equation}
where $\partial\Omega$ is the boundary of $\Omega$ and $\vec{n}$ is the unit outward normal vector to $\partial\Omega$. The existence and uniqueness of entropy solutions were established in \cite{GUO2011reaction} using regularization methods. In contrast to the classical TV-based model, system \eqref{Guosystem2011} incorporates a time-evolving fidelity term $w$, which enables effective separation of the cartoon component $u$ from the noise or texture component $f-u$, including oscillatory fractal-like boundaries. With an appropriate scale $\lambda$, the model achieves enhanced performance in noise removal and texture preservation. For $\varPhi(u)=G(\left|\nabla u\right|)$, \eqref{System variation} can be viewed as a generalization of the classical variational framework \eqref{Vese-formula}. To suppress the staircase effect induced by TV flow, the modification of \eqref{Guosystem2011} involves not only optimizing the regularization term $\varPhi$ \cite{BENDAIDA2023,ZHANG2019elliptic}, but also redesigning the primary diffusion term. For further advancements in diffusion-based methods, we refer readers to \cite{GAO2021system,Nokrane2023fourth}.

%Additionally, modifications to \eqref{Guosystem2011} have been explored from the perspective of diffusion equations. For further advancements in diffusion-based methods, we refer readers to \cite{GAO2021system,Nokrane2023fourth}. 

% Building on this, 

Compared to classical PDE methods based on local information, nonlocal methods are better suited for handling textures and repetitive structures in images. The pioneering work in this area is the nonlocal means (NLM) filter proposed by Buades et al. \cite{BUADES2005Ieee,BUADES2005Review}, which makes full use of redundant information in the image, thereby better preserves detailed features while removing noise. Gilboa et al. \cite{GILBOA2007Segmentation,GILBOA2009} generalized the NLM into a variational framework by introducing a nonlocal gradient operator, leading to the nonlocal TV (NLTV) regularization model. In \cite{Li2024variable}, the authors described the nonlocal diffusion process as
\begin{equation}\label{Nonlocal-diffusion}
	\frac{\partial u}{\partial t}=\int_{\Omega}c(x,y)\left|u(t,y)-u(t,x)\right|^{p-2}\left(u(t,y)-u(t,x)\right)dy,
\end{equation}
where $p\geq 1$, and $c(x,y)\geq 0$ is a nonlocal weight function. In the framework of \eqref{Nonlocal-diffusion}, various variants can be derived by selecting appropriate forms for $c(x,y)$. For instance, in the NLTV model, $c(x,y)$ is designed to measure the similarity between two image patches, making it well suited for texture image restoration \cite{Jung2011,Kindermann2005,SHI2021nonlocal}. For $c(x,y)=\left|x-y\right|^{-(N+sp)}$, \eqref{Nonlocal-diffusion} becomes a fractional $p$-Laplacian equation
\begin{equation}\label{fractional-diffusion}
	\frac{\partial u}{\partial t}=\int_{\Omega}\frac{1}{\left|x-y\right|^{N+sp}}\left|u(t,y)-u(t,x)\right|^{p-2}\left(u(t,y)-u(t,x)\right)dy,
\end{equation}
with $0<s<1$. Theoretical results including well-posedness of \eqref{fractional-diffusion} and convergence of the solution to the classical $p$-Laplacian equation as $s\rightarrow 1^{-}$ have been studied in \cite{ANDREU2008Neumann}, where the fractional $p$-Laplacian operator $(-\Delta)_p^s$ is defined as
\begin{equation}\notag
	(-\Delta)_p^s u(t,x)=\text{P.V.}\int_{\Omega}\frac{1}{\left|x-y\right|^{N+sp}}\left|u(t,x)-u(t,y)\right|^{p-2}\left(u(t,x)-u(t,y)\right)dy,
\end{equation}
which reduces to the fractional Laplacian $(-\Delta)^s$ for $p=2$. Applications of the fractional Laplacian include image processing, modeling of Brownian motion, L\'{e}vy flight processes, and anomalous diffusion phenomena \cite{GATTO2015,Devillanova2016,ZOUGAR2024}. In particular, in the case of $p=1$, equation \eqref{fractional-diffusion} becomes singular when $u(t,x)=u(t,y)$, and the main difficulty lies in properly defining the term $\frac{u(t,x)-u(t,y)}{\left|u(t,x)-u(t,y)\right|}$. To overcome this, Maz\'{o}n et al. \cite{MAZON2016} adopted the approach proposed in \cite{ANDREU2008Neumann} by replacing the original expression with an antisymmetric $L^{\infty}$-function $\varsigma$, such that $\left\|\varsigma\right\|_{L^{\infty}((0,T)\times\Omega\times\Omega)}\leq 1$, and $\varsigma(t,x,y)\in\operatorname{sign}\left(u(t,y)-u(t,x)\right)$, where $\operatorname{sign}(r)$ denotes the multivalued function, defined as $\operatorname{sign}(r) = \frac{r}{|r|}$ for $r \neq 0$, and $\operatorname{sign}(r) = [-1, 1]$ for $r = 0$. Based on maximum a posteriori (MAP) estimation, the weighted fractional $1$-Laplacian diffusion equation has been employed to remove multiplicative noise from images \cite{GAO2022Fractional}. By selecting spatially dependent $s(x,y)$, the authors of \cite{Li2024variable} investigated a class of variable-order fractional $1$-Laplacian diffusion equations, which enhanced the model's adaptability and achieved satisfactory denoising performance across the entire image.

%For further details, see \cite{Karami2019new,Kindermann2005}. 
%\begin{equation}\notag
%	(-\Delta)_1^su(x)=\text{P.V.}\int_{\Omega}\frac{1}{\left|x-y\right|^{N+s}}\frac{u(x)-u(y)}{\left|u(x)-u(y)\right|}dy,
%\end{equation}

%Applications of the fractional Laplacian include image processing \cite{GATTO2015,LIAN2024fractional}, and it also arises naturally in the modeling of Brownian motion, jump L\'{e}vy processes, and anomalous diffusion \cite{Caffarelli2012,Devillanova2016,ZOUGAR2024}.

%a common approach is to integrate local and nonlocal methods, as demonstrated in

Although nonlocal methods perform well in denoising homogeneous regions and preserving texture features, they tend to over-smooth low-contrast areas and leave residual noise near image edges \cite{Sutour2014}. To address these limitations and improve performance across the entire image, the idea of coupling local and nonlocal methods has gained increasing attention \cite{BOUKDIR2024system,SHI2021NARWA}. Recently, Liu et al. \cite{LIU2019system} replaced the primary TV diffusion term in \eqref{Guosystem2011} with fractional $1$-Laplacian diffusion and studied the following reaction-diffusion system
\begin{equation}\label{Guo system}
	\left\{\begin{array}{ll}  
		\displaystyle\frac{\partial u}{\partial t}=\Delta_1^s u-\lambda w, & x\in \Omega,~0<t<T,\\
		\displaystyle\frac{\partial w}{\partial t} = \Delta w - \left(f-u\right), & x\in \Omega,~0<t<T,\\
		\displaystyle u(0,x)=f,~w(0,x)=0, & x\in\Omega,\\
		\displaystyle \frac{\partial w}{\partial \vec{n}}(t,x)=0, & x\in \partial\Omega,~0<t<T.
	\end{array} \right.
\end{equation}
The well-posedness was established by employing the regularization methods adopted in \cite{MAZON2016}. To better capture oscillatory patterns while preserving image edges, Vese et al. \cite{VESE2008negative} modified the $H^{-1}$ norm in the variational model \eqref{System variation} to the $H^{-s}$ norm, where $s\geq 0$, and applied this model to image denoising, deblurring, and decomposition. Within their framework, the noise or texture component of the original image is modeled by tempered distributions in $H^{-s}(\Omega)$. G\'{a}rriz et al. \cite{Garriz2020coupling} proposed several hybrid energy functionals that benefit from local regularization in homogeneous regions and nonlocal regularization in textured regions. The existence and uniqueness of solutions to the gradient flows, which couple local and nonlocal diffusion operators acting in distinct subdomains, were established via the fixed point method. A higher-order version of the model introduced in \cite{Garriz2020coupling} was subsequently investigated in \cite{SHI2023coupling}. Numerical results demonstrate that the coupled method restores textures in images without causing new features and artifacts.

%Numerical experiments demonstrate that model \eqref{Guo system} achieves notable results in image decomposition by effectively preserving high-frequency contour structures while enhancing low-frequency textures in smooth areas.
%To capture oscillation patterns, Vese et al. \cite{VESE2008negative} transformed $H^{-1}(\Omega)$ to $H^{-s}(\Omega)$, $s\geq 0$, and applied \eqref{System variation} to tasks such as image denoising, deblurring, and decomposition. In their framework, the noise (or texture) component of the original image is modeled  by tempered distributions within the space $H^{-s}(\Omega)$. 

\subsection{The proposed model}

As far as we know, there is few literature about PDE system devoted to image deblurring. Inspired by the above ideas, we aim to restore a clear image with sharp edges from a degraded texture image corrupted by noise and blur. To this end, we combine forward-backward diffusion with fractional diffusion and propose the following nonlinear reaction-diffusion system:
\begin{equation}\label{New system}
	\left\{\begin{array}{ll}  
		\displaystyle\frac{\partial u}{\partial t}=\mathrm{div}\left[\left(\frac{1}{1+\left|\nabla u\right|^2}+\delta\left|\nabla u\right|^{\gamma-2}\right)\nabla u\right]-\lambda_1 u + \lambda_2 w, & x\in \Omega,~0<t<T,\\
		\displaystyle\frac{\partial w}{\partial t} = \Delta_p^s w - \lambda_3 K^{\prime}\left(Ku-f\right), & x\in \Omega,~0<t<T,\\
		\displaystyle u(0,x)=f,~w(0,x)=0, & x\in\Omega,\\
		\displaystyle \frac{\partial u}{\partial \vec{n}}(t,x)=0, & x\in \partial\Omega,~0<t<T,\\
	\end{array} \right.
\end{equation}
where $\Omega\subset\mathbb{R}^N$ ($N\geq 2$) is a bounded domain with appropriately smooth boundary $\partial\Omega$, $K\in\mathcal{L}\left(L^1(\Omega),L^2(\Omega)\right)$, and $\delta>0$, $1<\gamma\leq 2$, $0<s<1$, $1\leq p\leq 2$, $\lambda_1\geq 0$, $\lambda_2,\lambda_3>0$ are some given constants. Compared with the FBRPM model \cite{GUIDOTTI2013restoration}, the proposed model offers two main advantages:  long range information of the fidelity term between $u$ and
$f$ is taken into account, and additional scales are introduced for more flexible adjustment of the fidelity term.

%Numerical experiments demonstrate that our model exhibits favorable restoration performance across the entire image.

Different from previous models, the new system achieves image restoration through the interaction of two reaction-diffusion equations. Specifically, the fractional $p$-Laplacian simulates oscillation patterns and distinctive textures, enabling $w$ to leverage image structure and dynamically adjust the fidelity between $u$ and $f$. The inclusion of the regularized PM diffusion term provides adaptive control over forward and backward diffusion modes based on the local gradient $\nabla u$, facilitating smoothing in homogeneous regions and sharpening at edges. The enhanced image $u$, in turn, impacts the evolution and updating of $w$. In other words, the first equation in \eqref{New system} incorporates a different fidelity term $w$, derived from an evolution equation, which captures high-frequency redundant information of the image and is iteratively refined over time. As a result, we can select appropriate scale parameters to achieve an effective balance between denoising, deblurring, and the preservation of fine structures. Numerical experiments demonstrate that our model exhibits favorable restoration performance across the entire image.
%In comparison with the GLRPM model, our model retains texture features more effectively.

This paper focuses on the theoretical analysis of the proposed model. Although problem \eqref{New system} admits no classical weak solutions, we employ the Young measure representation developed by Demoulini et al. \cite{Demoulini1996} and regularization method to define Young measure solutions. Most research on equations of the forward-backward diffusion type has concentrated on evolution equations or systems without source terms \cite{Rieger2003}. In a way, our work also serves as a theoretical complement to equations of the type in \cite{GUIDOTTI2012backward,Wang2014,YIN2003Youngmeasure}. Unlike previous work, establishing the well-posedness of problem \eqref{New system} poses two major challenges.
\begin{itemize}
	\item[(i)] The presence of the term $w$ prevents us from applying Rothe's method to derive $L^{\infty}$ a-priori estimates for approximation solutions based on the variational structure of the steady-state version of problem \eqref{New system}. The existence of Young measure solutions cannot be proved following the method in \cite{YIN2003Youngmeasure}. To address this, we employ Moser's iteration or Stampacchia's truncation method to analyze the relaxed problem associated with the steady-state system, which is obtained via the relaxation theorem.
	
	%The presence of the function $w$ makes it impossible to derive uniform $L^{\infty}$ a-priori estimates for approximation solutions from the variational structure, after the steady-state system is obtained via time semi-discretization using Rothe's method. The existence of Young measure solutions cannot be proved following the method in \cite{YIN2003Youngmeasure}. To address this, we apply the relaxation theorem developed by Dacorogna et al. \cite{Dacorogna2007} to analyze the relaxed problem associated with the steady-state system. By employing Moser's iteration or Stampacchia's truncation method, we derive estimates for the approximation solutions to problem \eqref{New system}.
	\item[(ii)] For $p=1$, the singular nature of problem \eqref{New system} hinders direct analysis of its well-posedness. To overcome this, we adopt the regularization method in \cite{MAZON2016} and study a family of approximation problems involving fractional $p$-Laplacian ($1<p<2$) diffusion terms. By establishing uniform estimates for the approximate solutions and applying the standard limiting process $p\rightarrow 1^{+}$, we obtain the existence of Young measure solutions.
	
	%For $p=1$, problem \eqref{New system} may be singular, which prevents suitable definition of solutions and direct analysis of well-posedness. To overcome this, we adopt the regularization method in \cite{MAZON2016}. We consider a sequence of approximation problems involving fractional $p$-Laplacian diffusion terms. Through uniform estimates of the approximation solutions and the standard limiting process as $p\rightarrow 1^{+}$, we obtain the existence of Young measure solutions.
\end{itemize}
%study its well-posedness.

%Moreover, while \eqref{New system} lacks a classical solution, it is well-posed in the weaker sense of Young measure solutions, which can be constructed via time discretization and successive minimization.

The structure of the paper is as follows. Section \ref{sec2} provides some preliminaries on fractional Sobolev spaces and Young measures, and states the main result. In Section \ref{sec3.1}, we establish the existence, uniqueness, and continuous dependence of the Young measure solution to problem \eqref{New system} for $1<p\leq 2$, while the case $p=1$ is addressed in Section \ref{sec3.2}. In Section \ref{sec4}, we present numerical experiments and demonstrate the effectiveness of the proposed model through comparisons with other existing approaches.

\section{Preliminaries and main result}\label{sec2}
In this section, we first present some necessary preliminaries on fractional Sobolev spaces for later use. Let $\Omega\subset\mathbb{R}^N$, $N\geq 2$ be a bounded smooth enough open domain. 
For brevity, we denote $Q_T=(0,T)\times\Omega$ and $Q_T^{1,2}=(0,T)\times\Omega\times\Omega$. For any $1 \leq p<\infty$ and $0 < s < 1$, the Gagliardo seminorm of a measurable function $u$ in $\Omega$ is given by
\begin{equation}
	\left[u\right]_{W^{s,p}(\Omega)}=\left(\int_{\Omega}\int_{\Omega}\frac{\left|u(x)-u(y)\right|^p}{\left|x-y\right|^{N+sp}} dxdy\right)^{\frac{1}{p}}.\notag
\end{equation}
The fractional Sobolev space is defined as
\begin{equation}
	W^{s,p}(\Omega)=\left\{u\in L^p(\Omega):\left[u\right]_{W^{s,p}(\Omega)}<+\infty \right\},\notag
\end{equation}
which forms an intermediary Banach space between $L^p(\Omega)$ and $W^{1,p}(\Omega)$, equipped with the norm
\begin{equation}
	\left\|u\right\|_{W^{s,p}(\Omega)}=\left\|u\right\|_{L^p(\Omega)}+\left[u\right]_{W^{s,p}(\Omega)}.\notag
\end{equation}
For $1<p<\infty$, the fractional Sobolev space $W^{s,p}(\Omega)$ is reflexive and separable. Furthermore, as $s\rightarrow 1^{-}$, the space $W^{s,p}(\Omega)$ becomes $W^{1,p}(\Omega)$. Consider the functional $\mathcal{B}_p^s(u):L^2(\Omega)\rightarrow [0,+\infty]$ defined by
\begin{equation}
	\mathcal{B}_p^s(u)=
	\begin{cases}
		\left[u\right]^p_{W^{s,p}(\Omega)}, & \text{if } u \in W^{s,p}(\Omega), \\
		+\infty, & \text{elsewhere}.
	\end{cases}\notag
\end{equation}
Then $\mathcal{B}_p^s$ is convex and weak lower-semi continuous in $L^2(\Omega)$. Moreover, the subdifferential of $\mathcal{B}_p^s$ is maximal monotone in $L^2(\Omega)$. To properly express the fractional operator $(-\Delta)_1^s$, we define the multivalued function $\operatorname{sign}(r)$ as follows
\begin{equation}
	\operatorname{sign}(r)=
	\begin{cases}
		\frac{r}{\left|r\right|}, & \text{if } r \neq 0, \\
		[-1,1], & \text{if } r = 0.
	\end{cases}\notag
\end{equation}
It is straightforward to verify that $\operatorname{sign}(r)$ is monotone, that is, for any $a,b\in\mathbb{R}$, $A\in\operatorname{sign}(a)$, $B\in\operatorname{sign}(b)$, it holds that $(A-B)(a-b)\geq 0$. For more details, we refer to \cite{BRASCO2014,DINEZZA2012}. Now, we present the compact embedding results for fractional Sobolev spaces.
\begin{lemma}\cite{DINEZZA2012}\label{embedding}
	Let $0<s<1$ and $p\geq 1$ such that $sp<N$. Let $\Omega\subset\mathbb{R}^N$ be a smooth bounded domain. Assume that
	\begin{equation}
		1\leq q<p_s^{\ast}:=\frac{Np}{N-sp}. \notag
	\end{equation}
	Then there exists a positive constant $C=C(N,p,q,s,\Omega)$ such that for any $f\in W^{s,p}(\Omega)$, it holds that
	\begin{equation}
		\left\|f\right\|_{L^q(\Omega)} \leq C\left\|f\right\|_{W^{s, p}(\Omega)}.\notag
	\end{equation}
	That is, the space $W^{s,p}(\Omega)$ is continuously embedded in $L^q(\Omega)$. Moreover, this embedding is compact.
\end{lemma}

For readers' convenience, we now introduce some essential definitions and results on Young measures, beginning with common function spaces. Let $p\geq 1$, we define
\begin{align}
	&\mathcal{E}_0^p(\mathbb{R}^N)=\left\{\phi\in C(\mathbb{R}^N):\lim_{\left|\theta\right|\rightarrow \infty}\frac{\left|\phi(\theta)\right|}{1+\left|\theta\right|^p}\text{ exists}\right\},\notag\\
	&
	\mathcal{E}^p(\mathbb{R}^N)=\left\{\phi\in C(\mathbb{R}^N):\sup_{\theta\in\mathbb{R}^N}\frac{\left|\phi(\theta)\right|}{1+\left|\theta\right|^p}<+\infty\right\},\notag
\end{align}
endowed with the norm
\begin{equation}
	\left\| \phi\right\|_{\mathcal{E}^p(\mathbb{R}^N)}=\sup_{\theta\in\mathbb{R}^N} \frac{\left|\phi(\theta)\right|}{1+\left|\theta\right|^p},\quad\phi\in \mathcal{E}^p(\mathbb{R}^N).\notag
\end{equation}
As mentioned in \cite{Kinderlehrer1994gradient}, $\mathcal{E}_0^p(\mathbb{R}^N)$ is a separable Banach space, while $\mathcal{E}^p(\mathbb{R}^N)$ is an inseparable space. Let $C_0(\mathbb{R}^N)$ be the closure of continuous functions on $\mathbb{R}^N$ with compact support. The dual of $C_0(\mathbb{R}^N)$ can be identified with $\mathcal{M}(\mathbb{R}^N)$, the space of signed Radon measures with finite mass, through the pairing
\begin{equation}
	\langle\mu,f\rangle=\int_{\mathbb{R}^N}f(\theta)d\mu(\theta),~f\in C_0(\mathbb{R}^N),~\mu\in\mathcal{M}(\mathbb{R}^N).\notag
\end{equation}
A map $\nu:\Omega\rightarrow\mathcal{M}(\mathbb{R}^N)$ is called weakly $\ast$ measurable if the functions $x\mapsto\langle\nu_x,f\rangle$ are measurable for all $f\in C_0(\mathbb{R}^N)$, where $\nu_x=\nu(x)$.

\begin{lemma}[Fundamental theorem on Young measures]\label{fundamental theorem}\cite{BALL2005}
	Let $u_j:\Omega\rightarrow\mathbb{R}^N$ \emph{(}$j\geq 1$\emph{)} be a sequence of measurable functions. Then there exists a subsequence $\{u_{j_k}\}_{k=1}^{\infty}$ and a weakly $\ast$ measurable map $\nu:\Omega\rightarrow\mathcal{M}(\mathbb{R}^N)$, such that
	\begin{itemize}
		\item [\emph{(i)}] $\nu(x) \geq 0,\|\nu(x)\|_{\mathcal{M}\left(\mathbb{R}^N\right)}=\int_{\mathbb{R}^N} d \nu_x \leq 1, \text { a.e. } x \in \Omega $;
		\item [\emph{(ii)}] For all $\varphi\in C_0(\mathbb{R})^N$, $\varphi(u_{j_k})\rightharpoonup\langle\nu,\varphi\rangle$ weakly $\ast$ in $L^{\infty}(\Omega)$;
		\item [\emph{(iii)}] Furthermore, one has $\|\nu(x)\|_{\mathcal{M}\left(\mathbb{R}^N\right)}=1$, a.e. $x\in\Omega$ if and only if
		\begin{align}
			\lim _{L \rightarrow \infty} \sup _{k \geq 1} \left|\{x \in\Omega:|u_{j_k}(x)| \geq L\}\right|=0.\notag
		\end{align}
	\end{itemize}
	The map $\nu:\Omega\rightarrow\mathcal{M}(\mathbb{R}^N)$ is called the Young measure in $\mathbb{R}^N$ generated by the sequence $\{u_{j_k}\}_{k=1}^{\infty}$.
\end{lemma}

\begin{defn}\cite{Demoulini1996}\label{DW1q}
	A Young measure $\nu=(\nu_x)_{x\in \Omega}$ on $\mathbb{R}^N$ is called a $W^{1,p}(\Omega)$-gradient Young measure ($ p\geq 1$) if
	\begin{itemize}
		\item [(i)] the functions $x\mapsto \langle\nu_x,f\rangle$ are measurable for all $f$ bounded continuous on $\mathbb{R}^N$;
		\item [(ii)] there exists a sequence of functions $\{u^k\}_{k=1}^{\infty} \subseteq W^{1,p}(\Omega)$ such that the following representation formula holds:
		\begin{equation}
			\lim _{k \rightarrow \infty} \int_E g(\nabla u^k(x)) d x=\int_E\langle \nu_x, g\rangle d x,\label{rfYM}
		\end{equation}
		for all measurable subsets $E \subseteq \Omega$ and all $g \in \mathcal{E}_0^p (\mathbb{R}^N)$. We also refer to $\nu$ as the $W^{1,p}(\Omega)$-gradient Young measure generated by $\{\nabla u^k\}_{k=1}^{\infty}$ and to $\{\nabla u^k\}_{k=1}^{\infty}$ as the $W^{1,p}(\Omega)$-gradient generating sequence of $\nu$. Moreover, the representation formula (\ref{rfYM}) remains valid for all $g\in\mathcal{E}^p(\mathbb{R}^N)$.
	\end{itemize}
\end{defn}

The $W^{1,p}(\Omega)$-gradient Young measures are characterized by the following properties:
\begin{lemma}\cite{Kinderlehrer1994gradient}\label{criteria W1p Young}
	Let $\nu=(\nu_x)_{x\in\Omega}$ be a Young measure on $\mathbb{R}^N$. Then $\nu$ is a $W^{1,p}(\Omega)$-gradient Young measure if and only if
	\begin{itemize}
		\item[\emph{(i)}] there exists $u\in W^{1,p}(\Omega)$ such that $\nabla u(x)=\langle\nu_x,\mathrm{id}\rangle$, \emph{a.e.} $x\in \Omega$, where $\mathrm{id}$ is the unit mapping in $\mathbb{R}^N$;
		\item[\emph{(ii)}] Jensen's inequality $\phi(\nabla u(x))\leq \langle\nu_x,\phi\rangle$, \emph{a.e.} $x\in \Omega$, holds for all $\phi\in\mathcal{E}^p(\mathbb{R}^N)$ continuous, quasiconvex, and bounded below; 
		\item[\emph{(iii)}] the function $x\mapsto\int_{\mathbb{R}^N}\left|\theta\right|^p d\nu_x(\theta)$ is in $L^1(\Omega)$.
	\end{itemize}
\end{lemma}
\begin{lemma}\cite{Kinderlehrer1992,Kinderlehrer1994gradient}\label{kp1}
	Suppose $f \in \mathcal{E}^p(\mathbb{R}^N)$ with $p \geq 1$, is quasiconvex and bounded below, and let $\{u^k\}_{k=1}^{\infty}$ converge weakly to $u$ in $W^{1,p}(\Omega)$. Then
	\begin{itemize}
		\item [\emph{(i)}] for all measurable $E \subseteq \Omega$,
		\begin{equation}
			\int_E f(\nabla u) dx \leq\liminf_{k \rightarrow \infty} \int_E f(\nabla u^k) dx;\notag
		\end{equation}
		\item [\emph{(ii)}] if, in addition,
		\begin{equation}
			\lim _{k \rightarrow \infty} \int_{\Omega} f(\nabla u^k) d x=\int_{\Omega} f(\nabla u) dx,\notag
		\end{equation}
		then $\{f(\nabla u^k)\}_{k=1}^{\infty}$ converges weakly to $f(\nabla u)$ in $L^1(\Omega)$;
		\item [\emph{(iii)}] Furthermore, let $\nu$ be a Young measure generated by the sequence $\{\nabla u^k\}_{k=1}^{\infty}$ and assume that
		\begin{equation}
			\max\{c\left|\theta\right|^p-1,0\} \leq f(\theta) \leq C\left|\theta\right|^p+1,\notag
		\end{equation}
		for some $0<c \leq C$. Then $\nu$ is a $W^{1,p}\left(\Omega\right)$-gradient Young measure.
	\end{itemize}
\end{lemma}

We now state the weak compactness of sequences of gradient-generated Young measures.
\begin{lemma}\cite{Demoulini1996,YIN2003Youngmeasure}\label{wcp}
	Let $1 \leq p\leq 2$. Suppose $\{\nu^\beta=(\nu_{t,x}^\beta)_{(t,x) \in Q_T}\}_{\beta>0}$ is a family of $W^{1,p}(Q_T)$-gradient Young measures, each generated by $\{\nabla u^{\beta, m}\}_{m=1}^{\infty}$, where $u^{\beta, m}\in W^{1,p}(Q_T)$ is uniformly bounded in $\beta$ and $m$. Then, there exist a subsequence $\{\nu^{\beta_j}\}_{j=1}^{\infty}$ of $\{\nu^\beta\}_{\beta>0}$ and a $W^{1,p}(Q_T)$-gradient Young measure $\nu$ such that
	\begin{itemize}
		\item [\emph{(i)}] $\{\langle\nu^{\beta_j},\varphi\rangle\}_{j=1}^{\infty}$ converges weakly $\ast$ to $\langle\nu,\varphi\rangle$ in $L^{\infty}(Q_T)$ for all $\varphi\in C_0(\mathbb{R}^N)$;
		\item [\emph{(ii)}] for $1 \leq q<p$, $\{\langle\nu^{\beta_j},\varphi\rangle\}_{j=1}^{\infty}$ converges weakly to $\langle\nu,\varphi\rangle$ in $L^1(Q_T)$ for all $\varphi\in \mathcal{E}_0^q(\mathbb{R}^N)$;
		\item [\emph{(iii)}]  $\{\langle\nu^{\beta_j},\varphi\rangle\}_{j=1}^{\infty}$ converges to $\langle\nu,\varphi\rangle$ in the biting sense for all $\varphi\in \mathcal{E}_0^p(\mathbb{R}^N)$, namely, there exists a decreasing sequence of subsets $E_{j+1}\subset E_j$ of $Q_T$ with $\lim_{j \rightarrow \infty}|E_j|=0$, such that $\{\langle\nu^{\beta_j},\varphi\rangle\}_{j=1}^{\infty}$ converges weakly to $\langle\nu,\varphi\rangle$ in $L^1(Q_T\backslash E_j)$ for all $j$.
	\end{itemize}
\end{lemma}
\begin{rmk}\label{rmk1}
	Suppose $(\nu^{\beta})_{\beta>0}$, with $\nu^{\beta} = (\nu^{\beta}_{t,x})_{(t,x)\in Q_T}$, is a sequence of Young measures bounded
	in $L^1(Q_T;(\mathcal{E}_0^p(\mathbb{R}^N))^{\prime})$. For each $\beta$, let $\{\nabla v^{\beta,k}\}_{k=1}^{\infty}$ be the generating gradients, where $v^{\beta,k}\in L^p(0,T;W^{1,p}(\Omega))$. Then $(\nu^{\beta})_{\beta>0}$ is bounded in $L^{p}(Q_T;(\mathcal{E}_0^1(\mathbb{R}^N))^{\prime})\cap L^{\infty}(Q_T;\mathcal{M}(\mathbb{R}^N))$ and a subsequence of the $\{v^{\beta,k}\}_{\beta,k}$ is bounded in $L^p(0,T;W^{1,p}(\Omega))$ uniformly in $\beta$ and $k$. Thus Lemma \ref{wcp} applies.
\end{rmk}

Recall that the heat flux $\vec{q}_{\gamma}$ is expressed as:
\begin{equation}
	\vec{q}_{\gamma}(\theta)=\frac{\theta}{1+\left|\theta\right|^2}+\delta\left|\theta\right|^{\gamma-2}\theta,\quad\theta\in\mathbb{R}^N,\notag
\end{equation}
where $\vec{q}_{\gamma}=\nabla\varphi_{\gamma}$ with some non-convex potential $\varphi_{\gamma}\in C^1(\mathbb{R}^N)$ satisfying the structure conditions:
\begin{equation}\label{structure condition 1}
	\left|\vec{q}_{\gamma}(\theta)\right|\leq\tilde{C}_{\gamma}\left|\theta\right|^{\gamma-1},\quad\max\left\{\tilde{c}_{\gamma}\left|\theta\right|^{\gamma}-1,0\right\} \leq\varphi_{\gamma}(\theta) \leq\tilde{C}_{\gamma}\left|\theta\right|^{\gamma}+1,\quad\theta\in\mathbb{R}^N,
\end{equation}
for some $0<\tilde{c}_{\gamma}\leq\tilde{C}_{\gamma}$.
Let $\varphi_{\gamma}^{**}$ be the convexification of $\varphi_{\gamma}$, defined as $\varphi_{\gamma}^{**}(\theta)=\{g(\theta):g\leq\varphi,~g \text{ convex}\}$. Since $\varphi_{\gamma}\in C^1(\mathbb{R}^N)$, it follows that $\varphi_{\gamma}^{**}\in C^1(\mathbb{R}^N)$ and is convex. We set $\vec{\varrho}_{\gamma}=\nabla\varphi_{\gamma}^{**}$. Note that $\vec{q}_{\gamma}=\vec{\varrho}_{\gamma}$ on the set $\{\theta\in\mathbb{R}^N:\varphi_{\gamma}(\theta)=\varphi_{\gamma}^{**}(\theta)\}$. Additionally, $\varphi_{\gamma}^{**}$ and $\vec{\varrho}_{\gamma}$ satisfy the same growth conditions \eqref{structure condition 1} as $\varphi_{\gamma}$ and $\vec{q}_{\gamma}$, respectively, i.e.
\begin{equation}\label{structure condition 2}
	\left|\vec{\varrho}_{\gamma}(\theta)\right|\leq\tilde{C}_{\gamma}\left|\theta\right|^{\gamma-1}+1,\quad\max\left\{\tilde{c}_{\gamma}\left|\theta\right|^{\gamma}-1,0\right\} \leq\varphi_{\gamma}^{**}(\theta) \leq\tilde{C}_{\gamma}\left|\theta\right|^{\gamma}+1,\quad\theta\in\mathbb{R}^N.
\end{equation}

Next, we define the Young measure solution to system \eqref{New system} and present the main result as follows:
\begin{defn}\label{defn solution system}
	Given $f\in L^{\infty}(\Omega)\cap W^{1,\gamma}(\Omega)$, a couple of functions $(u,w)$ defined on $Q_T$ is called a Young measure solution of problem \eqref{New system}, if it meets the following conditions:
	\begin{itemize}
		\item [(i)] $u\in L^{\infty}(0,T;W^{1,\gamma}(\Omega))\cap L^{\infty}(Q_T)$, $\frac{\partial u}{\partial t}\in L^2(Q_T)$, with $u(0,x)=f(x)$, a.e. $x\in\Omega$;
		
		$w\in L^{\infty}(0,T; W^{s,p}(\Omega))\cap L^{\infty}(Q_T)$, $\frac{\partial w}{\partial t}\in L^2(Q_T)$, satisfying $w(0,x)=0$, a.e. $x\in\Omega$.
		
		\item [(ii)] There exists a $W^{1,\gamma}(Q_T)$-gradient Young measure $\nu=(\nu_{t,x})_{(t,x)\in Q_T}$ such that for every $\zeta\in L^2(Q_T)\cap L^{1}(0,T;W^{1,\gamma}(\Omega))$, the following integral equality holds:
		\begin{equation}\label{equation-for-u-nu}
			\iint_{Q_T}\left(\frac{\partial u}{\partial t}\zeta+\langle\nu,\vec{q}_{\gamma}\rangle\cdot\nabla\zeta+\left(\lambda_1 u-\lambda_2 w\right)\zeta\right) dxdt=0,
		\end{equation}
		and
		\begin{align}
			&\nabla u(t,x)=\langle\nu_{t,x},\mathrm{id}\rangle,~\text{a.e. }(t,x)\in Q_T,\label{gradientforu}\\
			&\langle \nu_{t,x},\vec{q}_{\gamma}\cdot\mathrm{id}\rangle=\langle \nu_{t,x},\vec{q}_{\gamma}\rangle\cdot\langle \nu_{t,x},\mathrm{id}\rangle,~\text{a.e. }(t,x)\in Q_T,\label{independentfornu}\\
			&\mathrm{supp} \nu_{t,x}\subseteq \{\theta\in\mathbb{R}^N:\varphi_{\gamma}(\theta)=\varphi_{\gamma}^{**}(\theta)\},~\text{a.e. }(t,x)\in Q_T,\label{suppfornu}
		\end{align}
		where $\mathrm{id}$ is the unit mapping in $\mathbb{R}^N$.
		\item [(iii)] ($1<p\leq 2$) For every $\psi\in L^2(Q_T)\cap L^1(0,T;W^{s,p}(\Omega))$, the following integral equality holds:
		\begin{align}
			\iint_{Q_T}&\frac{\partial w}{\partial t}\psi dxdt+\frac{1}{2}\iiint_{Q_T^{1,2}}\frac{\left|w(t,y)-w(t,x)\right|^{p-2}\left(w(t,y)-w(t,x)\right)}{\left|x-y\right|^{N+s p}}\left(\psi(t,y)-\psi(t,x)\right) dxdydt\notag\\
			&=\lambda_3 \iint_{Q_T}\left(f-Ku\right)K\psi dxdt,\label{equation-w}
		\end{align}
		\item [(iv)] ($p=1$) There exist a function $\varsigma(t,x,y)\in L^{\infty}(Q_T^{1,2})$ with $\left\|\varsigma(t,x,y)\right\|_{L^{\infty}(Q_T^{1,2})}\leq 1$ that satisfies $\varsigma(t,x,y)=-\varsigma(t,y,x)$ and $\varsigma(t,x,y)\in\operatorname{sign}\left(w(t,y)-w(t,x)\right)$ for a.e. $(t,x,y)\in Q_T^{1,2}$, such that for every $\psi\in L^2(Q_T)\cap L^1(0,T;W^{s,1}(\Omega))$, the following identity holds:
		\begin{align}
			\iint_{Q_T}\frac{\partial w}{\partial t}\psi dxdt+\frac{1}{2}\iiint_{Q_T^{1,2}}\frac{\varsigma(t,x,y)}{\left|x-y\right|^{N+s}}\left(\psi(t,y)-\psi(t,x)\right)dxdydt=\lambda_3 \iint_{Q_T}\left(f-Ku\right)K\psi dxdt,\label{equation-w-1}
		\end{align}
	\end{itemize}
\end{defn}
\begin{theorem}\label{theorem 1}
	If $f\in L^{\infty}(\Omega)\cap W^{1,\gamma}(\Omega)$, the system \eqref{New system} admits a unique Young measure solution.
\end{theorem}

\section{Proof of the main result}\label{sec3}
In this section, we divide the proof of Theorem \ref{theorem 1} into two parts. We first analyze the existence, uniqueness, and continuous dependence of Young measure solutions to problem \eqref{New system} for $1< p \leq 2$. The case 
$p=1$ will be treated by the regularization method based on the well-posedness of the problem in the regular case.
\subsection{The regular case $1< p \leq 2$}\label{sec3.1}

%We employ Rothe's method \cite{WU2006} and the variational principle, in conjunction with the relaxation theorem \cite{Dacorogna2007} and Moser's iteration, to establish the existence of Young measure solutions for system \eqref{New system}. The continuous dependence of the solutions will be discussed by the independence property (15).

%Subsequently, we demonstrate the continuous dependence of the Young measure solution using the independence property \eqref{independentfornu}. 

We employ Rothe's method \cite{WU2006} and variational principle, together with relaxation theorem \cite{Dacorogna2007} and Moser's iteration, to prove the existence of Young measure solutions for system \eqref{New system}. Continuous dependence will be addressed using independence property \eqref{independentfornu}. The proof of existence is primarily inspired by references \cite{Demoulini1996,Wang2014,YIN2003Youngmeasure}. To apply Rothe's method, we first consider the semi-discretization in time for \eqref{New system}. Prior to this, we define the discretization of functions on $Q_T$ with respect to $t$. For any positive integer $n$ and function $\omega:Q_T\rightarrow\mathbb{R}$, let
\begin{equation}
	\omega_n^i(x)=\omega\left(x,i\frac{T}{n}\right),\label{func-discretization}
\end{equation}
where $i=1,2,\cdots,n$. By discretizing system \eqref{New system}, we obtain the corresponding steady-state system. However, due to the non-convexity of $\varphi_{\gamma}$, we study the following (discretized) relaxed system derived from relaxation theorem:
\begin{align}
	& \frac{m}{T}(z-w_m^{j-1})+ (-\Delta)_{p}^s z= \lambda_3 K^{\prime}\left(f-Ku^{j-1}_m\right),~w_m^0=0,\label{Discrete-Delta_p^s}\\
	& \frac{m}{T}(v-u_m^{j-1})=\operatorname{div}\left(\vec{\varrho}_{\gamma}(\nabla v)\right)-\lambda_1 u_m^{j-1}+\lambda_2  w_m^j,~u_m^0=f,\label{Discrete-fb}
\end{align}
where $m$ is a positive integer, and $u_m^{j-1}\in L^2(\Omega)\cap W^{1,\gamma}(\Omega)$, $w_m^{j-1}\in L^2(\Omega)\cap W^{s,p}(\Omega)$ are defined as in \eqref{func-discretization}, with $j=1,2,\cdots,m$. Next, we analyze the above steady-state system starting from \eqref{Discrete-Delta_p^s}, where the weak solution is defined as follows:
\begin{defn}\label{weak solution to Discrete Delta-p-s}
	Given $f\in L^2(\Omega)\cap W^{1,\gamma}(\Omega)$, a function $z\in L^2(\Omega)\cap W^{s,p}(\Omega)$ is called a weak solution to problem \eqref{Discrete-Delta_p^s} if, for any $\phi\in L^2(\Omega)\cap W^{s,p}(\Omega)$, it satisfies
	\begin{align}
		\int_{\Omega}&\left(z-w_m^{j-1}\right)\phi dx +  \frac{T}{2m}\int_{\Omega}\int_{\Omega}\frac{1}{\left|x-y\right|^{N+sp}}\left|z(x)-z(y)\right|^{p-2}\left(z(x)-z(y)\right)\left(\phi(x)-\phi(y)\right)dxdy\notag\\
		& = \frac{T}{m}\lambda_3\int_{\Omega}\left(f-K u_m^{j-1}\right) K \phi dx,\label{Defn-p^s}
	\end{align}
	and we say $v\in L^2(\Omega)\cap W^{1,\gamma}(\Omega)$ is a weak solution to problem \eqref{Discrete-fb} if, for any $\eta\in L^2(\Omega)\cap W^{1,\gamma}(\Omega)$, the following integral equality holds:
	\begin{equation}\label{Defn-fb}
		\int_{\Omega}\left( \left(v-u_m^{j-1}\right)\eta + \frac{T}{m}\vec{\varrho}_{\gamma}(\nabla v)\cdot\nabla\eta +\frac{T}{m}\lambda_1 u_m^{j-1}\eta-\frac{T}{m}\lambda_2 w_m^j\eta \right) dx = 0.
	\end{equation}
\end{defn}
\begin{lemma}\label{discrete system infty}
	Assume $f\in L^2(\Omega)\cap W^{1,\gamma}(\Omega)$, and that $u_m^{j-1}\in L^2(\Omega)\cap W^{1,\gamma}(\Omega)$ and $w_m^{j-1}\in L^2(\Omega)\cap W^{s,p}(\Omega)$, then the discretized problems \eqref{Discrete-Delta_p^s} and \eqref{Discrete-fb} admit unique weak solutions $w_m^{j}\in L^2(\Omega)\cap W^{s,p}(\Omega)$ and $u_m^{j}\in L^2(\Omega)\cap W^{1,\gamma}(\Omega)$, respectively. Furthermore, if $f\in L^{\infty}(\Omega)$, we have $w_m^j\in L^{\infty}(\Omega)$ and $u_m^{j}\in L^{\infty}(\Omega)$ for $j=1,2,\cdots,m$, with the following estimates:
	\begin{align}
		& \left\|w_m^j\right\|_{L^{\infty}(\Omega)}+\left\|w_m^j\right\|_{W^{s,p}(\Omega)}+\frac{m}{T}\sum_{i=1}^{j}\left\|w_m^i-w_m^{i-1} \right\|_{L^2(\Omega)}^2\leq M_1,\label{estimate w_mj}\\
		& \left\|u_m^j\right\|_{L^{\infty}(\Omega)}+\left\|u_m^j\right\|_{W^{1,\gamma}(\Omega)}+\frac{m}{T}\sum_{i=1}^{j}\left\|u_m^i-u_m^{i-1} \right\|_{L^2(\Omega)}^2\leq M_1,\label{estimate u_mj}
	\end{align}
	where $M_1$ is a positive constant depending only on $T$, $\Omega$, $K$, $\tilde{c}_{\gamma}$, $\lambda_1$, $\lambda_2$, $\lambda_3$, and $\left\|f\right\|_{L^{\infty}(\Omega)}$, $j=1,2,\cdots,m$.
\end{lemma}
\begin{pf}
	We begin by studying the existence of weak solutions to problem \eqref{Discrete-Delta_p^s}. For any $z\in L^2(\Omega)\cap W^{s,p}(\Omega)$, we define the following functional:
	\begin{equation}
		F_m(z;w_m^{j-1})=\frac{1}{2}\int_{\Omega}\left(z-w_m^{j-1}\right)^2 dx+\frac{T}{2mp}\left[z\right]_{W^{s,p}(\Omega)}^p
		-\frac{T}{m}\lambda_3 \int_{\Omega}K^{\prime}\left(f-K u_m^{j-1}\right) z dx,\notag
	\end{equation}
	which is weak lower semi-continuous in $L^2(\Omega)$ \cite{GAO2022Fractional,LIU2019system}. In view of $K\in \mathcal{L}(L^1(\Omega),L^2(\Omega))$, and  by applying Cauchy's inequality, we obtain
	\begin{align}
		F_m(z;w_m^{j-1})\geq& \frac{1}{2}\left\|z\right\|_{L^2(\Omega)}^2-\int_{\Omega}w_m^{j-1}z dx -\frac{T}{m}\lambda_3\int_{\Omega}\left(f-Ku_m^{j-1}\right)Kz dx \notag\\
		\geq& -\left\|w_m^{j-1}\right\|_{L^2(\Omega)}^2-\frac{2T^2}{m^2}\lambda_3^2\left\| K\right\|^2\left|\Omega\right|\left\| f\right\|_{L^2(\Omega)}^2-\frac{2T^2}{m^2}\lambda_3^2\left\| K\right\|^4\left|\Omega\right|^2\left\| u^{j-1}_m\right\|_{L^2(\Omega)}^2.\notag
	\end{align}
	This shows that $F_m(z;w_m^{j-1})$ is bounded from below, implying that its infimum in $L^2(\Omega)\cap W^{s,p}(\Omega)$ is finite. Assume that $\{w_m^{j,k}\}_{k=1}^{\infty}\subseteq L^2(\Omega)\cap W^{s,p}(\Omega)$ is a minimizing sequence for $F_m(z;w_m^{j-1})$, satisfying
	\begin{equation}
		F_m(w_m^{j,k};w_m^{j-1})<\inf\left\{F_m(z;w_m^{j-1}):z\in L^2(\Omega)\cap W^{s,p}(\Omega)\right\}+\frac{T}{m^2 k},~k=1,2,\cdots,m.\notag
	\end{equation}
	Using Cauchy's inequality again, we have
	\begin{align}
		& \frac{1}{2}\int_{\Omega}\left(w_m^{j,k}-w_m^{j-1}\right)^2 dx+\frac{T}{2mp}\left[w_m^{j,k}\right]_{W^{s,p}(\Omega)}^p - \frac{T}{m^2 k}\notag\\
		\leq &  \frac{T}{2mp}\left[w_m^{j-1}\right]_{W^{s,p}(\Omega)}^p + \frac{T}{m}\lambda_3\int_{\Omega}\left(f-K u_m^{j-1}\right)\left(Kw_m^{j,k}-K w_m^{j-1}\right) dx \notag \\ 
		\leq & \frac{T}{2mp}\left[w_m^{j-1}\right]_{W^{s,p}(\Omega)}^p + \frac{2T^2}{m^2}\lambda_3^2\left\| K\right\|^2\left|\Omega\right|\int_{\Omega}\left|f\right|^2 dx + \frac{2T^2}{m^2}\lambda_3^2\left\| K\right\|^4\left|\Omega\right|^2\int_{\Omega}\left|u_m^{j-1}\right|^2 dx\notag\\
		&+ \frac{1}{4}\int_{\Omega}\left(w_m^{j,k}-w_m^{j-1}\right)^2 dx. \notag
	\end{align}
	Thus, we obtain
	\begin{align}
		&\left[w_m^{j,k}\right]_{W^{s,p}(\Omega)}^p+\frac{mp}{2T}\int_{\Omega}\left(w_m^{j,k}-w_m^{j-1}\right)^2 dx\notag\\
		\leq & \left[w_m^{j-1}\right]_{W^{s,p}(\Omega)}^p + \frac{4Tp}{m}\lambda_3^2\left\| K\right\|^2\left|\Omega\right|\left\|f\right\|_{L^2(\Omega))}^2 + \frac{4Tp}{m}\lambda_3^2\left\| K\right\|^4\left|\Omega\right|^2\left\|u_m^{j-1}\right\|_{L^2(\Omega))}^2 +\frac{2p}{m k}.\label{Wsp w_mjk}
	\end{align}
	For any $z\in L^2(\Omega)\cap W^{s,p}(\Omega)$, applying H\"{o}lder's inequality, we derive
	\begin{align}
		F_m(z;w_m^{j-1})\geq&\frac{1}{2}\left\|z\right\|_{L^2(\Omega)}^2 - \left\|w_m^{j-1} \right\|_{L^2(\Omega)}\left\|z\right\|_{L^2(\Omega)} - \frac{T}{m}\lambda_3 \left\|f - K u_m^{j-1} \right\|_{L^{2}(\Omega)}\left\|Kz\right\|_{L^2(\Omega)}\notag\\
		\geq & \frac{1}{2}\left\|z\right\|_{L^2(\Omega)}^2 - \left\|w_m^{j-1} \right\|_{L^2(\Omega)}\left\|z\right\|_{L^2(\Omega)} - \frac{T}{m}\lambda_3 \left\|K\right\|\sqrt{\left|\Omega\right|} \left\|f\right\|_{L^2(\Omega)}\left\|z\right\|_{L^2(\Omega)}\notag\\
		& - \frac{T}{m}\lambda_3 \left\| K\right\|^2\left|\Omega\right|\left\| u_m^{j-1}\right\|_{L^2(\Omega)} \left\|z\right\|_{L^2(\Omega)},\notag
	\end{align}
	which shows that $F_m(z;w_m^{j-1})$ is coercive on $L^2(\Omega)$. From \eqref{Wsp w_mjk}, the sequence $\{w_m^{j,k}\}_{k=1}^{\infty}$ is uniformly bounded in $L^2(\Omega)\cap W^{s,p}(\Omega)$ with respect to $k$. The embedding of $W^{s,p}(\Omega)$ into $L^p(\Omega)$ is compact, as shown by Lemma \ref{embedding}. Therefore, there exists $w_m^j\in L^2(\Omega)\cap W^{s,p}(\Omega)$ and a subsequence of $\{w_m^{j,k}\}_{k=1}^{\infty}$, still denoted by itself, such that
	\begin{equation}\label{convergence w_mjk}
		w_m^{j,k}\rightarrow w_m^j\text{ strongly in } L^{p}(\Omega),~w_m^{j,k}\rightharpoonup w_m^j\text{ weakly in } L^{2}(\Omega),~w_m^{j,k}\rightharpoonup w_m^j\text{ weakly in } W^{s,p}(\Omega).
	\end{equation}
	Then, we have
	\begin{align}
		&\left[w_m^j\right]_{W^{s,p}(\Omega)}^p+\frac{mp}{2T}\int_{\Omega}\left(w_m^{j}-w_m^{j-1}\right)^2 dx\notag\\
		\leq & \liminf_{k\rightarrow\infty}\left(\left[w_m^{j,k}\right]_{W^{s,p}(\Omega)}^p+\frac{mp}{2T}\int_{\Omega}\left(w_m^{j,k}-w_m^{j-1}\right)^2 dx\right)\notag\\
		\leq & \left[w_m^{j-1}\right]_{W^{s,p}(\Omega)}^p + \frac{4Tp}{m}\lambda_3^2\left\| K\right\|^2\left|\Omega\right|\left\|f\right\|_{L^2(\Omega))}^2 + \frac{4Tp}{m}\lambda_3^2\left\| K\right\|^4\left|\Omega\right|^2\left\|u_m^{j-1}\right\|_{L^2(\Omega))}^2 +\frac{2p}{m k}.\label{iteration w_mj}
	\end{align}
	Due to the weak lower semi-continuity of $F_m(z;w_m^{j-1})$ in $L^2(\Omega)$ and its strict convexity in $L^2(\Omega)\cap W^{s,p}(\Omega)$, $w_m^j$ is the unique minimizer of $F_m(z;w_m^{j-1})$. Note that for any $\phi\in L^2(\Omega)\cap W^{s,p}(\Omega)$, since $\epsilon\mapsto F_m(w_m^j+\epsilon\phi;w_m^{j-1})$ attains its minimum at $\epsilon=0$, we compute
	\begin{align}
		0=&\frac{d}{d\epsilon}F_m(w_m^j+\epsilon\phi;w_m^{j-1})\Big|_{\epsilon=0}\notag\\
		=&\int_{\Omega}\left(w_m^j-w_m^{j-1}\right)\phi dx-\frac{T}{m}\lambda_3\int_{\Omega}\left(f-K u_m^{j-1}\right)K \phi dx\notag\\
		&+\frac{T}{2m}\int_{\Omega}\int_{\Omega}\frac{1}{\left|x-y\right|^{N+sp}}\left|w_m^j(x)-w_m^j(y)\right|^{p-2}\left(w_m^j(x)-w_m^j(y)\right)\left(\phi(x)-\phi(y)\right) dxdy.\notag
	\end{align}
	By Definition \ref{weak solution to Discrete Delta-p-s}, it follows that $w_m^j$ is a weak solution to problem \eqref{Discrete-Delta_p^s}. Now, we prove the uniqueness of the weak solution $w_m^j$. Assume $w_m^j$, $\tilde{w}_m^j\in L^2(\Omega)\cap W^{s,p}(\Omega)$ are two weak solutions to \eqref{Discrete-Delta_p^s}. From \eqref{Defn-p^s}, we have
	\begin{align}
		\frac{T}{2m} & \int_{\Omega}\int_{\Omega}\frac{1}{\left|x-y\right|^{N+sp}} \left|w_m^j(x)-w_m^j(y)\right|^{p-2}\left(w_m^j(x)-w_m^j(y)\right) \left(w_m^j(x)-w_m^j(y)-\left(\tilde{w}_m^j(x)-\tilde{w}_m^j(y)\right)\right)dxdy\notag\\
		& = -\int_{\Omega}w_m^j  \left(w_m^j-\tilde{w}_m^j\right) dx + \frac{T}{m}\lambda_3\int_{\Omega}\left(f-Ku_m^{j-1}\right)\left(K w_m^j- K \tilde{w}_m^j\right) dx + \int_{\Omega}w_m^{j-1}\left(w_m^j-\tilde{w}_m^j\right) dx,\notag
	\end{align}
	and
	\begin{align}
		\frac{T}{2m}&\int_{\Omega}\int_{\Omega}\frac{1}{\left|x-y\right|^{N+sp}}\left|\tilde{w}_m^j(x)-\tilde{w}_m^j(y)\right|^{p-2}\left(\tilde{w}_m^j(x)-\tilde{w}_m^j(y)\right)\left(\tilde{w}_m^j(x)-\tilde{w}_m^j(y)-\left(w_m^j(x)-w_m^j(y)\right)\right)dxdy\notag\\
		& =-\int_{\Omega}\tilde{w}_m^j \left(\tilde{w}_m^j-w_m^j\right) dx+ \frac{T}{m}\lambda_3\int_{\Omega}\left(f-K u_m^{j-1}\right)\left(K \tilde{w}_m^j-K w_m^j\right) dx + \int_{\Omega}w_m^{j-1}\left(\tilde{w}_m^j-w_m^j\right) dx.\notag
	\end{align}
	By summing the two equalities above, we get
	\begin{align}
		\frac{T}{2m} & \int_{\Omega}\int_{\Omega}\frac{1}{\left|x-y\right|^{N+sp}}\left[\left|w_m^j(x)-w_m^j(y)\right|^{p-2}\left(w_m^j(x)-w_m^j(y)\right) - \left|\tilde{w}_m^j(x)-\tilde{w}_m^j(y)\right|^{p-2}\left(\tilde{w}_m^j(x)-\tilde{w}_m^j(y)\right)\right]\times\notag\\ &\left(w_m^j(x)-w_m^j(y)-\left(\tilde{w}_m^j(x)-\tilde{w}_m^j(y)\right)\right) dxdy+\int_{\Omega}\left(w_m^j-\tilde{w}_m^j\right)^2 dx = 0,\notag
	\end{align}
	which, together with the monotonicity of the function $t\mapsto t|t|^{p-2}$ for $p>1$, implies that $w_m^j=\tilde{w}_m^j$, a.e. $x\in\Omega$, thus proving the uniqueness of $w_m^j$.
	
	Next, we use the weak solution $w_m^j\in L^2(\Omega)\cap W^{s,p}(\Omega)$ obtained in the previous step to analyze problem \eqref{Discrete-fb}. For any $v\in L^2(\Omega)\cap W^{1,\gamma}(\Omega)$, we define
	\begin{equation}
		E_m^{**}(v;u_m^{j-1})=\int_{\Omega}\left(\frac{T}{m}\varphi_{\gamma}^{**}(\nabla v)+\frac{1}{2}\left(v-u_m^{j-1}\right)^2+\frac{T}{m}\left(\lambda_1 u_m^{j-1}-\lambda_2 w_m^j\right)v\right) dx. \notag
	\end{equation}
	According to the relaxation theorem \cite{Dacorogna2007}, $E_m^{**}(v;u_m^{j-1})$ attains the same infimum as the following functional in $L^2(\Omega)\cap W^{1,\gamma}(\Omega)$:
	\begin{equation}
		E_m(v;u_m^{j-1})=\int_{\Omega}\left(\frac{T}{m}\varphi_{\gamma}(\nabla v)+\frac{1}{2}\left(v-u_m^{j-1}\right)^2+\frac{T}{m}\left(\lambda_1 u_m^{j-1}-\lambda_2 w_m^j\right)v\right) dx .\notag
	\end{equation}
	From the structural conditions \eqref{structure condition 1} and \eqref{structure condition 2}, it follows that the two functionals above are bounded from below. Let $\{u_m^{j,k}\}_{k=1}^{\infty}\subset L^2(\Omega)\cap W^{1,\gamma}(\Omega)$	denote the minimizing sequence of $E_m(v;u_m^{j-1})$, satisfying
	\begin{equation}
		E_m(u_m^{j,k};u_m^{j-1})<\inf\left\{E_m(v;u_m^{j-1}):v\in L^2(\Omega)\cap W^{1,\gamma}(\Omega)\right\}+\frac{T}{m^2 k},~k=1,2,\cdots,m.\notag
	\end{equation}
	Since $\varphi_{\gamma}^{**}\leq \varphi_{\gamma}$, we have
	\begin{equation}
		E_m^{**}(u_m^{j,k};u_m^{j-1})<\inf\left\{E_m^{**}(v;u_m^{j-1}):v\in L^2(\Omega)\cap W^{1,\gamma}(\Omega)\right\}+\frac{T}{m^2 k},~k=1,2,\cdots,m,\notag
	\end{equation}
	which, together with Cauchy's inequality,  implies that
	\begin{align}
		&\frac{T}{m}\int_{\Omega}\varphi_{\gamma}^{**}(\nabla u_m^{j,k}) dx+\frac{1}{2}\int_{\Omega}\left(u_m^{j,k}-u_m^{j-1}\right)^2dx-\frac{T}{m^2 k}\notag\\ \leq&\frac{T}{m}\int_{\Omega}\varphi_{\gamma}^{**}(\nabla u_m^{j-1}) dx +\frac{T}{m}\int_{\Omega}\left(\lambda_1 u_m^{j-1}-\lambda_2 w_m^j\right)\left(u_m^{j-1}-u_m^{j,k}\right) dx\notag\\
		\leq & \frac{T}{m}\int_{\Omega}\varphi_{\gamma}^{**}(\nabla u_m^{j-1}) dx +\frac{T^2}{m^2}\int_{\Omega}\left(\lambda_1 u_m^{j-1}-\lambda_2 w_m^j\right)^2 dx+\frac{1}{4}\int_{\Omega}\left(u_m^{j-1}-u_m^{j,k}\right)^2 dx.\notag
	\end{align}
	Hence, we obtain
	\begin{align}
		&\int_{\Omega}\varphi_{\gamma}^{**}(\nabla u_m^{j,k}) dx+\frac{m}{4T}\int_{\Omega}\left(u_m^{j,k}-u_m^{j-1}\right)^2dx\notag\\
		\leq & \int_{\Omega}\varphi_{\gamma}^{**}(\nabla u_m^{j-1}) dx +\frac{T}{m}\int_{\Omega}\left(\lambda_1 u_m^{j-1}-\lambda_2 w_m^j\right)^2 dx+\frac{1}{m k}\notag\\
		\leq & \int_{\Omega}\varphi_{\gamma}^{**}(\nabla u_m^{j-1}) dx +\frac{2T}{m}\lambda_1^2\left\|u_m^{j-1} \right\|_{L^{2}(\Omega)}^2 +\frac{2T}{m}\lambda_2^2 \left\|w_m^j\right\|_{L^2(\Omega)}^2+\frac{1}{m k}. \label{iteration u_mjk}
	\end{align}
	In view of the inequality
	\begin{align}
		E_m^{**}(v;u_m^{j-1})\geq\frac{1}{2}\left\|v \right\|_{L^2(\Omega)}^2-\left(1+\frac{T}{m}\lambda_1\right) \left\|u_m^{j-1} \right\|_{L^2(\Omega)} \left\|v \right\|_{L^2(\Omega)}  -\frac{T}{m}\lambda_2 \left\|w_m^{j} \right\|_{L^2(\Omega)} \left\|v \right\|_{L^2(\Omega)}, \notag
	\end{align}
	we know that $E_m^{**}(v;u_m^{j-1})$ is coercive on $L^2(\Omega)$. Combining \eqref{structure condition 2} and \eqref{iteration u_mjk}, we conclude that the minimizing sequence $\{u_m^{j,k}\}_{k=1}^{\infty}$ is uniformly bounded in $L^2(\Omega)\cap W^{1,\gamma}(\Omega)$ with respect to $k$. By the Rellich-Kondrachov theorem, there exists a convergent subsequence of $\{u_m^{j,k}\}_{k=1}^{\infty}$, which we denote by the sequence itself, and a function $u_m^j\in L^2(\Omega)\cap W^{1,\gamma}(\Omega)$, such that
	\begin{equation}\label{convergence u_mjk}
		u_m^{j,k}\rightarrow u_m^j \text{ strongly in } L^{\gamma}(\Omega),~u_m^{j,k}\rightharpoonup u_m^j \text{ weakly in } L^{2}(\Omega),~u_m^{j,k}\rightharpoonup u_m^j \text{ weakly in } W^{1,\gamma}(\Omega).
	\end{equation}
	From Lemma \ref{kp1}, we know that $E_m^{**}$ is weakly lower semi-continuous in $L^2(\Omega)\cap W^{1,\gamma}(\Omega)$, so $u_m^j$ is a minimizer of $E_m^{**}$. Based on \eqref{iteration u_mjk}, \eqref{convergence u_mjk}, and \eqref{integrablelimit-for-umjk-to-umj}, we deduce
	\begin{align}
		&\int_{\Omega}\varphi_{\gamma}^{**}(\nabla u_m^{j}) dx+\frac{m}{4T}\int_{\Omega}\left(u_m^{j}-u_m^{j-1}\right)^2dx\notag\\
		\leq & \liminf_{k\rightarrow\infty}\left(\int_{\Omega}\varphi_{\gamma}^{**}(\nabla u_m^{j,k}) dx+\frac{m}{4T}\int_{\Omega}\left(u_m^{j,k}-u_m^{j-1}\right)^2dx\right) \notag\\
		\leq & \int_{\Omega}\varphi_{\gamma}^{**}(\nabla u_m^{j-1}) dx +\frac{2T}{m}\lambda_1^2\left\|u_m^{j-1} \right\|_{L^{2}(\Omega)}^2 +\frac{2T}{m}\lambda_2^2 \left\|w_m^j\right\|_{L^2(\Omega)}^2+\frac{1}{m k}.\label{new iteration u_mj}
	\end{align}
	Now, we derive the Euler-Lagrange equation satisfied by $u_m^j$. For any test function $\eta\in L^2(\Omega)\cap W^{1,\gamma}(\Omega)$, since the map $\varepsilon\mapsto E_m^{**}(u_m^j+\varepsilon\eta;u_m^{j-1})$ has a minimum at $\varepsilon=0$, we have
	\begin{align}
		0=&\frac{d}{d\epsilon}E_m^{**}(u_m^j+\varepsilon\eta;u_m^{j-1})\Big|_{\varepsilon=0}\notag\\
		=& \frac{T}{m}\int_{\Omega} \vec{\varrho}_{\gamma}(\nabla u_m^j)\cdot\nabla\eta dx + \int_{\Omega}\left(u_m^j-u_m^{j-1}\right)\eta dx + \frac{T}{m}\int_{\Omega}\left(\lambda_1 u_m^{j-1}-\lambda_2 w_m^j\right)\eta dx,\label{equation u_mj}
	\end{align}
	which shows that $u_m^j$ is a weak solution of problem \eqref{Discrete-fb}. Note that $\varphi_{\gamma}^{**}$ is convex, implying that $\vec{\varrho}_{\gamma}$ is monotone, and thus the uniqueness of the weak solution $u_m^j$ follows.
	
	%The above process outlines the iterative procedure for generating the weak solutions $w_m^j \in L^2(\Omega)\cap W^{s,p}(\Omega)$ of problem \eqref{Discrete-Delta_p^s} and $u_m^j \in L^2(\Omega)\cap W^{1,\gamma}(\Omega)$ of problem \eqref{Discrete-fb}, assuming the prior solutions $u_m^{j-1} \in L^2(\Omega)\cap W^{1,\gamma}(\Omega)$ and $w_m^{j-1} \in L^2(\Omega)\cap W^{s,p}(\Omega)$.
	%The above process describes the iterative construction of weak solutions $w_m^j\in L^2(\Omega)\cap W^{s,p}(\Omega)$ for problem (19) and $u_m^j\in L^2(\Omega)\cap W^{1,r}(\Omega)$ for problem (20), based on the previously obtained solutions $u_m^{j-1}$ and $w_m^{j-1}$.
	
	The above process describes the iterative construction of weak solutions $w_m^j \in L^2(\Omega) \cap W^{s,p}(\Omega)$ and $u_m^j \in L^2(\Omega) \cap W^{1,r}(\Omega)$ to problems \eqref{Discrete-Delta_p^s} and \eqref{Discrete-fb}, respectively, based on the previously obtained solutions $w_m^{j-1}$ and $u_m^{j-1}$. Specifically, we use $u_m^{j-1}$ and $w_m^{j-1}$ to compute the weak solution $w_m^j$ of \eqref{Discrete-Delta_p^s}, and then combine $w_m^j$ with $u_m^{j-1}$ to obtain the weak solution $u_m^j$ of \eqref{Discrete-fb}. This iteration reflects the image restoration  mechanism of the proposed model.
	
	%Specifically, we first use $u_m^{j-1}$ and $w_m^{j-1}$ to compute the unique weak solution $w_m^j$ of \eqref{Discrete-Delta_p^s}, then use it with $u_m^{j-1}$ to obtain the unique weak solution $u_m^j$ of \eqref{Discrete-fb}. This iteration reflects the image restoration  mechanism of the proposed model.

	Next, we prove that $w_m^j \in L^{\infty}(\Omega)$ and $u_m^j \in L^{\infty}(\Omega)$, and derive uniform boundedness estimates under the condition $f \in L^{\infty}(\Omega) \cap W^{1,\gamma}(\Omega)$. When $j=1$, we set $u_m^0 = f$ and $w_m^0 = 0$, and then compute $w_m^1$ and $u_m^1$. Repeating the above process in turn, we obtain $w_m^j$ satisfying \eqref{iteration w_mj} and $u_m^j $ satisfying \eqref{new iteration u_mj} for $j = 1, 2, \dots, m$. Now, we apply the Moser's iteration method to establish $L^{\infty}(\Omega)$ estimates for $w_m^j$ and $u_m^j$. Let $r\geq 1$, and choose $\phi=|w_m^j|^r w_m^j$ as a test function in \eqref{Defn-p^s}. Then, we obtain
	\begin{align}
		\int_{\Omega}\left|w_m^j\right|^{r+2}dx&+\frac{T}{2m}\int_{\Omega}\int_{\Omega}\frac{1}{\left|x-y\right|^{N+sp}}\left|w_m^j(x)-w_m^j(y)\right|^{p-2}\left(w_m^j(x)-w_m^j(y)\right) \left(\left|w_m^j(x)\right|^{r}w_m^j(x)\right.\notag\\
		& \left.-\left|w_m^j(y)\right|^{r}w_m^j(y)\right)dxdy
		=\frac{T}{m}\lambda_3\int_{\Omega}K^{\prime}\left(f-K u_m^{j-1}\right)\left|w_m^j\right|^{r}w_m^j dx+ \int_{\Omega}w_m^{j-1}\left|w_m^j\right|^{r}w_m^j dx.\notag
	\end{align}
	Using H\"{o}lder's inequality yields
	\begin{align}
		\left\| w_m^j \right\|_{L^{r+2}(\Omega)}^{r+2}\leq & \left\|w_m^{j-1}\right\|_{L^{r+2}(\Omega)}\left\|w_m^j\right\|_{L^{r+2}(\Omega)}^{r+1}+\frac{T}{m}\lambda_3 \left\| K^{\prime}\left(f-K u_m^{j-1}\right)\right\|_{L^{\infty}(\Omega)}\left\| w_m^j \right\|_{L^{r+1}(\Omega)}^{r+1} \notag\\
		\leq & \left\|w_m^{j-1}\right\|_{L^{r+2}(\Omega)}\left\|w_m^j\right\|_{L^{r+2}(\Omega)}^{r+1}+\frac{T}{m}\lambda_3 \left\|K^{\prime}\right\|\left|\Omega\right|^{\frac{1}{r+2}}\left\| f-K u_m^{j-1}\right\|_{L^{2}(\Omega)}\left\| w_m^j \right\|_{L^{r+2}(\Omega)}^{r+1} \notag\\
		\leq & \left\|w_m^{j-1}\right\|_{L^{r+2}(\Omega)}\left\|w_m^j\right\|_{L^{r+2}(\Omega)}^{r+1} + \frac{T}{m}\lambda_3 \left\|K^{\prime}\right\|\left|\Omega\right|^{\frac{1}{r+2}}\left\|f\right\|_{L^{2}(\Omega)}\left\| w_m^j \right\|_{L^{r+2}(\Omega)}^{r+1} \notag\\
		&+\frac{T}{m}\lambda_3 \left\|K^{\prime}\right\|\left\|K\right\|\left|\Omega\right|^{\frac{1}{r+2}}\left\| u_m^{j-1}\right\|_{L^{1}(\Omega)}\left\| w_m^j \right\|_{L^{r+2}(\Omega)}^{r+1}  \notag\\
		\leq & \left\|w_m^{j-1}\right\|_{L^{r+2}(\Omega)}\left\|w_m^j\right\|_{L^{r+2}(\Omega)}^{r+1} + \frac{T}{m}\lambda_3 \left\|K^{\prime}\right\|\left|\Omega\right|^{\frac{1}{r+2}}\left\|f\right\|_{L^{2}(\Omega)}\left\| w_m^j \right\|_{L^{r+2}(\Omega)}^{r+1} \notag\\
		& + \frac{T}{m}\lambda_3 \left\|K^{\prime}\right\|\left\|K\right\|\left|\Omega\right|\left\| u_m^{j-1}\right\|_{L^{r+2}(\Omega)}\left\| w_m^j \right\|_{L^{r+2}(\Omega)}^{r+1} \notag
	\end{align}
	which implies that
	\begin{align}
		\left\| w_m^j \right\|_{L^{r+2}(\Omega)}\leq \left\| w_m^{j-1} \right\|_{L^{r+2}(\Omega)}+\frac{T}{m}\lambda_3 \left\|K^{\prime}\right\|\left(\left|\Omega\right|+1\right)\left\|f\right\|_{L^{2}(\Omega)}+ \frac{T}{m}\lambda_3 \left\|K^{\prime}\right\|\left\|K\right\|\left|\Omega\right|\left\| u_m^{j-1}\right\|_{L^{r+2}(\Omega)}.\label{r+2 bound w_mj}
	\end{align}
	Similarly, since $u_m^j$ is the weak solution of \eqref{Discrete-fb}, we substitute the test function $\eta=|u_m^j|^{r}u_m^j$, where $r\geq 1$, into \eqref{Defn-fb}, and obtain
	\begin{align}
		&\int_{\Omega}\left|u_m^j\right|^{r+2} dx+\frac{T}{m}\int_{\Omega}\vec{\varrho}_{\gamma}(\nabla u_m^j)\cdot\left(r\left|u_m^j\right|^{r}\nabla u_m^j +\left|u_m^j\right|^{r}\nabla u_m^j\right) dx\notag\\
		=&\int_{\Omega}u_m^{j-1}\left|u_m^j\right|^r u_m^j dx-\frac{T}{m}\lambda_1\int_{\Omega}u_m^{j-1}\left|u_m^j\right|^r u_m^j dx + \frac{T}{m}\lambda_2\int_{\Omega}w_m^{j}\left|u_m^j\right|^r u_m^j dx . \notag
	\end{align}
	Given the convexity of $\varphi_{\gamma}(\theta)$ near $\theta = 0$, we have $\vec{\varrho}_{\gamma}(0) = 0$ and $\vec{\varrho}_{\gamma}(\nabla u_m^j)\cdot\nabla u_m^j\geq 0$. From \eqref{r+2 bound w_mj}, we derive
	\begin{align}
		\left\| u_m^j \right\|_{L^{r+2}(\Omega)}^{r+2}\leq &\left\| u_m^{j-1}\right\|_{L^{r+2}(\Omega)}\left\|u_m^j\right\|_{L^{r+2}(\Omega)}^{r+1}+\frac{T}{m}\lambda_1\left\|u_m^{j-1}\right\|_{L^{r+2}(\Omega)}\left\|u_m^j\right\|_{L^{r+2}(\Omega)}^{r+1}\notag\\
		&+\frac{T}{m}\lambda_2\left\|w_m^{j}\right\|_{L^{r+2}(\Omega)}\left\|u_m^j\right\|_{L^{r+2}(\Omega)}^{r+1} \notag\\
		\leq &\left\| u_m^{j-1}\right\|_{L^{r+2}(\Omega)}\left\|u_m^j\right\|_{L^{r+2}(\Omega)}^{r+1}+\frac{T}{m}\lambda_1\left\|u_m^{j-1}\right\|_{L^{r+2}(\Omega)}\left\|u_m^j\right\|_{L^{r+2}(\Omega)}^{r+1}\notag\\
		&+\frac{T}{m}\lambda_2\left\|u_m^j\right\|_{L^{r+2}(\Omega)}^{r+1}\left(\left\| w_m^{j-1} \right\|_{L^{r+2}(\Omega)}+\frac{T}{m}\lambda_3 \left\|K^{\prime}\right\|\left(\left|\Omega\right|+1\right)\left\|f\right\|_{L^{2}(\Omega)}\right.\notag\\
		&\left.+ \frac{T}{m}\lambda_3 \left\|K^{\prime}\right\|\left\|K\right\|\left|\Omega\right|\left\| u_m^{j-1}\right\|_{L^{r+2}(\Omega)}\right) \notag\\
		= & \left\| u_m^{j-1}\right\|_{L^{r+2}(\Omega)}\left\|u_m^j\right\|_{L^{r+2}(\Omega)}^{r+1}+\frac{T}{m}\lambda_1\left\|u_m^{j-1}\right\|_{L^{r+2}(\Omega)}\left\|u_m^j\right\|_{L^{r+2}(\Omega)}^{r+1} \notag\\
		& + \frac{T}{m}\lambda_2 \left\|w_m^{j-1}\right\|_{L^{r+2}(\Omega)}\left\|u_m^j\right\|_{L^{r+2}(\Omega)}^{r+1} \notag\\
		& + \frac{T^2}{m^2}\lambda_2 \lambda_3  \left\|K^{\prime}\right\|\left(\left|\Omega\right|+1\right)\left\|f\right\|_{L^{2}(\Omega)}\left\|u_m^j\right\|_{L^{r+2}(\Omega)}^{r+1} \notag\\
		& + \frac{T^2}{m^2}\lambda_2 \lambda_3  \left\|K^{\prime}\right\|\left\|K\right\|\left|\Omega\right|\left\| u_m^{j-1}\right\|_{L^{r+2}(\Omega)} \left\|u_m^j\right\|_{L^{r+2}(\Omega)}^{r+1}.
		\notag
	\end{align}
	Then
	\begin{align}
		\left\| u_m^j \right\|_{L^{r+2}(\Omega)}\leq & \left\| u_m^{j-1}\right\|_{L^{r+2}(\Omega)}+\frac{T}{m}\lambda_1\left\|u_m^{j-1}\right\|_{L^{r+2}(\Omega)} + \frac{T}{m}\lambda_2 \left\|w_m^{j-1}\right\|_{L^{r+2}(\Omega)}\notag\\
		& + \frac{T^2}{m}\lambda_2 \lambda_3  \left\|K^{\prime}\right\|\left(\left|\Omega\right|+1\right)\left\|f\right\|_{L^{2}(\Omega)} + \frac{T^2}{m}\lambda_2 \lambda_3  \left\|K^{\prime}\right\|\left\|K\right\|\left|\Omega\right|\left\| u_m^{j-1}\right\|_{L^{r+2}(\Omega)}. \notag
	\end{align}
	Thus, we obtain
	\begin{align}
		\left\| u_m^j \right\|_{L^{r+2}(\Omega)}+\left\| w_m^j \right\|_{L^{r+2}(\Omega)}\leq & \left(1+\frac{T}{m}\lambda_1 + \frac{T^2}{m}\lambda_2 \lambda_3  \left\|K^{\prime}\right\|\left\|K\right\|\left|\Omega\right|+ \frac{T}{m}\lambda_3 \left\|K^{\prime}\right\|\left\|K\right\|\left|\Omega\right|\right)\left\| u_m^{j-1} \right\|_{L^{r+2}(\Omega)}\notag\\
		& + \left(1+\frac{T}{m}\lambda_2\right)\left\| w_m^{j-1} \right\|_{L^{r+2}(\Omega)}+ \frac{T^2}{m}\lambda_2 \lambda_3  \left\|K^{\prime}\right\|\left(\left|\Omega\right|+1\right)\left\|f\right\|_{L^{2}(\Omega)} \notag\\
		& + \frac{T}{m}\lambda_3 \left\|K^{\prime}\right\|\left(\left|\Omega\right|+1\right)\left\|f\right\|_{L^{2}(\Omega)},\notag
	\end{align}
	which yields
	\begin{equation}\label{r+2 bound u_mj + w_mj}
		\left\| u_m^j \right\|_{L^{r+2}(\Omega)}+\left\| w_m^j \right\|_{L^{r+2}(\Omega)} \leq \left(1+\frac{C_0}{m}\right) \left(\left\| u_m^{j-1} \right\|_{L^{r+2}(\Omega)} + \left\| w_m^{j-1} \right\|_{L^{r+2}(\Omega)}\right)  + \frac{C_0}{m},
	\end{equation}
	where
	\begin{align}
		C_0=&\max\left\{T\lambda_1+T^2\lambda_2\lambda_3\left\| K^{\prime} \right\|\left\| K \right\|\left|\Omega\right|+T\lambda_3\left\| K^{\prime} \right\|\left\| K \right\|\left|\Omega\right|,T\lambda_2,T^2 \lambda_2\lambda_3 \left\|K^{\prime}\right\|\left(\left|\Omega\right|+1\right)\left\|f\right\|_{L^{2}(\Omega)}
		\right.\notag\\
		&\left.\quad \qquad+T\lambda_3 \left\|K^{\prime}\right\|\left(\left|\Omega\right|+1\right)\left\|f\right\|_{L^{2}(\Omega)}
		\right\}.\notag
	\end{align}
	Direct calculation shows that
	\begin{equation}
		\left(1+\frac{C_0}{m}\right)^m\leq e^{C_0},\quad \sum_{j=0}^{m-1} \left(1+\frac{C_0}{m}\right)^j\leq\frac{m}{C_0}e^{C_0}.\notag
	\end{equation}
	Using inequality \eqref{r+2 bound u_mj + w_mj}, by induction, we get
	\begin{equation}
		\left\| u_m^j \right\|_{L^{r+2}(\Omega)}+\left\| w_m^j \right\|_{L^{r+2}(\Omega)} \leq \left(1+\frac{C_0}{m}\right)^j \left\| f \right\|_{L^{r+2}(\Omega)} + \frac{C_0}{m} \sum_{i=0}^{j-1}\left(1+\frac{C_0}{m}\right)^i \leq e^{C_0}\left(\left|\Omega\right|+1\right)\left\| f \right\|_{L^{\infty}(\Omega)} + e^{C_0} .\notag
	\end{equation}
	Taking the limit $r\rightarrow\infty$ in the above inequality, we obtain
	\begin{equation}\label{infty u_mj + w_mj}
		\left\| u_m^j \right\|_{L^{\infty}(\Omega)}+\left\| w_m^j \right\|_{L^{\infty}(\Omega)} \leq e^{C_0}\left(\left|\Omega\right|+1\right) \left\| f \right\|_{L^{\infty}(\Omega)} + e^{C_0},~j=1,2,\cdots,m,
	\end{equation}
	which completes the $L^{\infty}(\Omega)$ estimates for $w_m^j$ and $u_m^j$. 
	
	Finally, we prove inequalities \eqref{estimate w_mj} and \eqref{estimate u_mj}. As $1<p\leq 2$, we deduce from \eqref{iteration w_mj} and \eqref{infty u_mj + w_mj} that there exists a positive constant $C_1$, depending only on $T$, $\Omega$, $K$, $\lambda_1$, $\lambda_2$, $\lambda_3$, and $\left\|f \right\|_{L^{\infty}(\Omega)}$, such that for $i=1,2,\cdots,j$
	\begin{equation}
		\left[w_m^i\right]_{W^{s,p}(\Omega)}^p+\frac{mp}{2T}\int_{\Omega}\left(w_m^{i}-w_m^{i-1}\right)^2 dx\leq \left[w_m^{i-1}\right]_{W^{s,p}(\Omega)}^p + \frac{C_1}{m}.\notag
	\end{equation}
	Summing the above inequality from $i=1$ to $j$ leads to
	\begin{equation}
		\left[w_m^j\right]_{W^{s,p}(\Omega)}^p+\frac{m}{2T}\sum_{i=1}^{j}\int_{\Omega}\left(w_m^i-w_m^{j-1}\right)^2 dx\leq j\frac{C_1}{m}\leq C_1.\notag
	\end{equation}
	Then, from Lemma \ref{embedding}, we obtain inequality \eqref{estimate w_mj}. Based on \eqref{new iteration u_mj} and \eqref{infty u_mj + w_mj}, it follows that
	\begin{equation}\label{iteration varphi u_mi}
		\int_{\Omega}\varphi_{\gamma}^{**}(\nabla u_m^{i}) dx+\frac{m}{4T}\int_{\Omega}\left(u_m^{i}-u_m^{i-1}\right)^2dx
		\leq   \int_{\Omega}\varphi_{\gamma}^{**}(\nabla u_m^{i-1}) dx +\frac{C_2}{m},
	\end{equation}
	where $C_2$ is a positive constant depending only on $T$, $\Omega$, $K$, $\lambda_1$, $\lambda_2$, $\lambda_3$, and $\left\| f \right\|_{L^{\infty}(\Omega)}$. Summing up \eqref{iteration varphi u_mi} from $i=1$ to $j$, we have
	\begin{equation}\label{bound varphi u_mj}
		\int_{\Omega}\varphi_{\gamma}^{**}(\nabla u_m^{j}) dx +\frac{m}{4T}\sum_{i=1}^{j}\int_{\Omega}\left(u_m^i-u_m^{i-1}\right)^2 dx \leq  \int_{\Omega}\varphi_{\gamma}^{**}(\nabla f) dx + j\frac{C_2}{m}\leq \int_{\Omega}\varphi_{\gamma}^{**}(\nabla f) dx + C_2.
	\end{equation}
	Combining \eqref{structure condition 2} and \eqref{bound varphi u_mj}, we derive \eqref{estimate u_mj}, thus completing the proof.
\end{pf}

To construct the approximate solution of the Young measure solution for system \eqref{New system} and the corresponding Young measure sequence, we use Lemma \ref{fundamental theorem} to introduce the Young measure generated by the (spatial) gradients of the minimization sequence $u_m^{j,k}$ in Lemma \ref{discrete system infty}. It follows from the assumption $f \in L^{\infty}(\Omega) \cap W^{1,\gamma}(\Omega)$ and estimate \eqref{infty u_mj + w_mj} that there exists a positive constant $l_0$, depending only on $T$, $\Omega$, $K$, $\lambda_1$, $\lambda_2$, $\lambda_3$, and $\left\| f \right\|_{L^{\infty}(\Omega)}$, such that
\begin{equation}
	E_m^{**}\left(T_{l_0}(u_m^{j,k});u_m^{j-1}\right)\leq E_m^{**}\left(u_m^{j,k};u_m^{j-1}\right),\notag
\end{equation}
where the truncation function $T_l$ is given by $T_l(t)=\max\left\{-l,\min\left\{l,t\right\}\right\}$. Without loss of generality, we may assume that $\|u_m^{j,k}\|_{L^{\infty}(\Omega)}\leq l_0$. Moreover, from \eqref{convergence u_mjk}, we derive that $u_m^{j,k}$ converges strongly to $u_m^j$ in $L^2(\Omega)$, which leads to
\begin{equation}\label{integrablelimit-for-umjk-to-umj}
	\int_{\Omega}\varphi_{\gamma}^{**}(\nabla u_m^j) dx=\lim_{k\rightarrow\infty}\int_{\Omega}\varphi_{\gamma}^{**}(\nabla u_m^{j,k}) dx=\lim_{k\rightarrow\infty}\int_{\Omega}\varphi_{\gamma}(\nabla u_m^{j,k}) dx.
\end{equation}

For $j=1,2,\cdots,m$, let $\nu^{m,j}=(\nu^{m,j}_x)_{x\in\Omega}$ denote the Young measure generated by $\{\nabla u_m^{j,k}\}_{k=1}^{\infty}$. By Lemma \ref{kp1}, and using \eqref{convergence u_mjk} and \eqref{integrablelimit-for-umjk-to-umj}, we know that $\nu^{m,j}$ is a $W^{1,\gamma}(\Omega)$-gradient Young measure. From the representation formula \eqref{rfYM}, we observe that $\nabla u_m^{j,k}$ converges weakly to $\langle\nu^{m,j},\mathrm{id}\rangle$ in $L^1(\Omega)$, where $\mathrm{id}$ is the unit mapping in $\mathbb{R}^N$. Then
\begin{equation}
	\nabla u_m^j(x)=\langle \nu_x^{m,j},\text{id}\rangle,\text{ a.e. }x\in\Omega,\label{gradient-for-umj}
\end{equation}
Since $\varphi_{\gamma}^{**}\in\mathcal{E}^{\gamma}(\mathbb{R}^N)$, $\varphi_{\gamma}\in\mathcal{E}_0^{\gamma}(\mathbb{R}^N)$, and (\ref{integrablelimit-for-umjk-to-umj}) holds, we deduce from \eqref{rfYM} that
\begin{equation}
	\int_{\Omega}\langle \nu^{m,j}, \varphi_{\gamma}^{**}\rangle dx=\lim_{k\rightarrow\infty}\int_{\Omega}\varphi_{\gamma}^{**}(\nabla u_m^{j,k}) dx=\lim_{k\rightarrow\infty}\int_{\Omega}\varphi_{\gamma}(\nabla u_m^{j,k}) dx=\int_{\Omega}\langle \nu^{m,j}, \varphi_{\gamma}\rangle dx,\notag
\end{equation}
which, together with $\varphi_{\gamma}^{**}\leq\varphi_{\gamma}$, implies that $\varphi_{\gamma}^{**}(\nabla u_m^{j}(x))=\langle \nu_x^{m,j}, \varphi_{\gamma}^{**}\rangle=\langle \nu_x^{m,j}, \varphi_{\gamma}\rangle$, a.e. $x\in\Omega$. Thus, we get
\begin{equation}\label{supp-for-umj}
	\operatorname{supp}\nu_x^{m,j}\subseteq\{\theta\in\mathbb{R}^N:\varphi_{\gamma}^{**}(\theta)=\varphi_{\gamma}(\theta)\},\text{ a.e. } x\in \Omega,~j=1,2,\cdots,m.
\end{equation}
Note that $\{\theta\in\mathbb{R}^N:\varphi_{\gamma}^{**}(\theta)=\varphi_{\gamma}(\theta)\}\subseteq\{\theta\in\mathbb{R}^N:\vec{\varrho}_{\gamma}(\theta)=\vec{q}_{\gamma}(\theta)\}$, so we obtain
\begin{equation}\label{nux-alpha=nux-alpha**}
	\langle \nu_x^{m,j},\vec{\varrho}_{\gamma}\rangle=\langle \nu_x^{m,j},\vec{q}_{\gamma}\rangle,\text{ a.e. }x\in\Omega,~ j=1,2,\cdots,m.
\end{equation}
Given that $\vec{\varrho}_{\gamma}\in\mathcal{E}_0^{\gamma}(\mathbb{R}^N)$ and $\vec{\varrho}_{\gamma}\cdot\mathrm{id}\in\mathcal{E}^{\gamma}(\mathbb{R}^N)$, it follows from \eqref{rfYM} and \eqref{convergence u_mjk} that
\begin{align}
	&\nabla u_m^{j,k}\rightharpoonup\langle\nu^{m,j},\mathrm{id}\rangle\text{ weakly in } L^{\gamma}(\Omega),~\vec{\varrho}_{\gamma}(\nabla u_m^{j,k})\rightharpoonup\langle\nu^{m,j},\vec{\varrho}_{\gamma}\rangle\text{ weakly in } L^{\frac{\gamma}{\gamma-1}}(\Omega),\notag\\
	&\vec{\varrho}_{\gamma}(\nabla u_m^{j,k})\cdot\nabla u_m^{j,k}\rightharpoonup\langle\nu^{m,j},\vec{\varrho}_{\gamma}\cdot\mathrm{id}\rangle \text{ weakly in } L^1(\Omega).\label{varrho_gamma_u_mjk}
\end{align}
Then, for any $\hat{\eta}\in  W^{1,\gamma}(\Omega)$, we have
\begin{align}
	&\lim_{k \rightarrow \infty}\int_{\Omega}\left(\vec{\varrho}_{\gamma}(\nabla u_m^{j,k})\cdot\nabla\hat{\eta}+\frac{m}{T}\left(u_m^{j,k}-u_m^{j-1}\right)\hat{\eta}+\left(\lambda_1 u_m^{j-1}-\lambda_2 w_m^j\right)\hat{\eta}\right) dx \notag\\
	=& \int_{\Omega}\left(\langle\nu^{m,j},\vec{\varrho}_{\gamma}\rangle\cdot\nabla\hat{\eta}+\frac{m}{T}\left(u_m^j-u_m^{j-1}\right)\hat{\eta}+\left(\lambda_1 u_m^{j-1}-\lambda_2 w_m^j\right)\hat{\eta}\right) dx.\notag
\end{align}
Taking $\hat{\eta}=(u_m^{j,k}-u_m^j)\tau$ with $\tau\in C_0^{\infty}(\Omega)$, and applying H\"{o}lder's inequality, we derive from \eqref{iteration u_mjk} that
\begin{align}
	& \limsup_{k\rightarrow\infty}\left|\int_{\Omega}\left(\vec{\varrho}_{\gamma}(\nabla u_m^{j,k})-\langle\nu^{m,j},\vec{\varrho}_{\gamma}\rangle \right)\cdot\left(\nabla u_m^{j,k}-\nabla u_m^{j}\right)\tau dx\right|\notag\\
	\leq & \limsup_{k\rightarrow\infty}\left|\int_{\Omega}\left(\vec{\varrho}_{\gamma}(\nabla u_m^{j,k})-\langle\nu^{m,j},\vec{\varrho}_{\gamma}\rangle \right)\cdot\nabla\left(u_m^{j,k}\tau-u_m^{j}\tau\right) dx\right|\notag\\
	& + \lim_{k\rightarrow\infty}\left|\int_{\Omega}\left(u_m^{j,k}-u_m^{j}\right)\left(\vec{\varrho}_{\gamma}(\nabla u_m^{j,k})-\langle\nu^{m,j},\vec{\varrho}_{\gamma}\rangle \right)\cdot\nabla\tau dx\right| \notag\\
	\leq & \lim_{k \rightarrow \infty} \left|\frac{m}{T}\int_{D}\left(u_m^j-u_m^{j,k}\right)^2 \tau dx\right|+ \lim_{k\rightarrow\infty}\int_{\Omega}\left|\left(u_m^{j,k}-u_m^{j}\right)\left(\vec{\varrho}_{\gamma}(\nabla u_m^{j,k})-\langle\nu^{m,j},\vec{\varrho}_{\gamma}\rangle \right)\cdot\nabla\tau\right| dx\notag\\
	\leq & \lim_{k \rightarrow \infty} \frac{m}{T}\left\|u_m^{j,k}-u_m^j\right\|^2_{L^{2}(\Omega)}\left\|\tau\right\|_{L^{\infty}(\Omega)}+\lim_{k \rightarrow \infty}C_3 \left\|u_m^{j,k}-u_m^j\right\|_{L^{\gamma}(\Omega)}\left\|\nabla\tau\right\|_{L^{\infty}(\Omega)},\notag
\end{align}
where $C_3$
is a positive constant depending only on $\tilde{C}_{\gamma}$, $M_1$, and $\Omega$. Using \eqref{convergence u_mjk} and \eqref{varrho_gamma_u_mjk}, we obtain
\begin{equation}
	\int_{\Omega}\left(\langle\nu^{m,j},\vec{\varrho}_{\gamma}\cdot\mathrm{id}\rangle-\langle\nu^{m,j},\vec{\varrho}_{\gamma}\rangle\cdot\langle\nu^{m,j},\mathrm{id}\rangle \right)\tau dx = 0,\notag
\end{equation}
which yields
\begin{equation}
	\langle\nu_x^{m,j},\vec{\varrho}_{\gamma}\cdot\mathrm{id}\rangle=\langle\nu_x^{m,j},\vec{\varrho}_{\gamma}\rangle\cdot\langle\nu_x^{m,j},\mathrm{id}\rangle,\text{ a.e. }x\in\Omega.\notag
\end{equation}
Combining \eqref{supp-for-umj} and \eqref{nux-alpha=nux-alpha**}, we deduce that
\begin{equation}\label{indenpendence nu_mj}
	\langle\nu_x^{m,j},\vec{q}_{\gamma}\cdot\mathrm{id}\rangle=\langle\nu_x^{m,j},\vec{q}_{\gamma}\rangle\cdot\langle\nu_x^{m,j},\mathrm{id}\rangle,\text{ a.e. }x\in\Omega.
\end{equation}

We now define the approximation solutions for the system \eqref{New system}. For $j=1,2,\cdots,m$, let $\chi_m^j$ denote the characteristic function of the interval $[(j-1)\frac{T}{m},j\frac{T}{m})$, and define $\lambda_m^j(t)$ as
\begin{equation}
	\lambda_m^j(t)=\left(\frac{m}{T}t-(j-1)\right)\chi_m^j(t),\quad 0\leq t\leq T.\notag
\end{equation}
For $(t,x)\in Q_T$, we define the following (Rothe function):
\begin{align}
	&u_m(t,x)=\sum_{j=1}^{m}\chi_m^j(t)\left\{u_m^{j-1}(x)+\lambda_m^j(t)\left(u_m^j(x)-u_m^{j-1}(x)\right)\right\},\notag\\
	&\bar{u}_m(t,x)=\sum_{j=1}^{m}\chi_m^j(t)u_m^j(x),~\tilde{u}_m(x,t)=\sum_{j=1}^{m}\chi_m^j(t)u_m^{j-1}(x),\notag\\
	&v_m^k(t,x)=\sum_{j=1}^{m}\chi_m^j(t)u_m^{j,k}(x),~k=1,2,\cdots,\notag\\
	&w_m(t,x)=\sum_{j=1}^{m}\chi_m^j(t)\left\{w_m^{j-1}(x)+\lambda_m^j(t)\left(w_m^j(x)-w_m^{j-1}(x)\right)\right\},\notag\\
	&\bar{w}_m(t,x)=\sum_{j=1}^{m}\chi_m^j(t)w_m^j(x),\notag
\end{align}
and the Young measure:
\begin{equation}
	\nu^m=(\nu_{t,x}^m)_{(t,x)\in Q_T},~\nu_{t,x}^m=\sum_{j=1}^{m}\chi_m^j(t)\nu_x^{m,j}.\notag
\end{equation}
From Lemma \ref{discrete system infty}, we conclude that $u_m,~\bar{u}_m,~\tilde{u}_m\in L^\infty(0,T; W^{1,\gamma}(\Omega))\cap L^{\infty}(Q_T)$ with $\frac{\partial u_m}{\partial t} \in L^2(Q_T)$, and $w_m,~\bar{w}_m \in L^\infty(0,T; W^{s,p}(\Omega))\cap L^{\infty}(Q_T)$ with $\frac{\partial w_m}{\partial t} \in L^2(Q_T)$. Furthermore, $\nu^{m}\in L^1(Q_T;(\mathcal{E}_0^{\gamma}(\mathbb{R}^N))^{\prime})$ is the $W^{1,\gamma}(Q_T)$-gradient Young measure generated by $\{\nabla v_m^k\}_{k=1}^{\infty}$. It follows from (\ref{gradient-for-umj}) and (\ref{supp-for-umj}) that
\begin{align}
	&\nabla \bar{u}_m(t,x)=\langle \nu_{t,x}^m,\mathrm{id}\rangle,\text{ a.e. } (t,x)\in Q_T,\label{gradientforum}\\
	&\operatorname{supp}\nu_{t,x}^m\subseteq\{\theta\in\mathbb{R}^N:\varphi^{**}_{\gamma}(\theta)=\varphi_{\gamma}(\theta)\},\text{ a.e. } (t,x)\in Q_T,\label{suppform}
\end{align}

% thereby constructing the Young measure solution of the system.
Next, we derive uniformly bounded estimates for both the approximate solutions of system \eqref{New system} and the corresponding Young measures, serving as preparation for constructing the Young measure solution to the system.
\begin{lemma}\label{estimates-for-w-m}
	Under the assumption of Theorem \ref{theorem 1}, there exists a positive constant $M_2$, independent of $m$, $j$, $k$, $s$, and $p$ such that the following estimates hold:
	\begin{align}
		&\left\|w_m\right\|_{L^{\infty}(0,T; W^{s,p}(\Omega))}+\left\|w_m\right\|_{L^{\infty}(Q_T)}+\left\|\frac{\partial w_m}{\partial t} \right\|_{L^{2}(Q_T)}\leq M_2,\label{bound-wm-Linftynorm}\\
		&\left\|\bar{w}_m\right\|_{L^{\infty}(0,T ;W^{s,p}(\Omega))}+\left\|\bar{w}_m\right\|_{L^{\infty}(Q_T)}\leq M_2,\label{bound-barwm-Linftynorm}\\
		&\left\|w_m-\bar{w}_m\right\|_{L^2(Q_T)}^2\leq\frac{1}{m^2}M_2,\label{bound-w_m-barw_m-L2norm}\\
		&	\left\|u_m\right\|_{L^{\infty}(0,T;W^{1,\gamma}(\Omega))}+\left\|u_m\right\|_{L^{\infty}(Q_T)}+\left\|\frac{\partial u_m}{\partial t}\right\|_{L^2(Q_T)}\leq M_2,\label{u_m-Linfty0TW_01p-Linfty-du_mdt-L2-bound}\\
		&\left\|\bar{u}_m\right\|_{L^{\infty}(0,T;W^{1,\gamma}(\Omega))}+\left\|\bar{u}_m\right\|_{L^{\infty}(Q_T)}\leq M_2,~\left\|\tilde{u}_m\right\|_{L^{\infty}(0,T;W^{1,\gamma}(\Omega))}+\left\|\tilde{u}_m\right\|_{L^{\infty}(Q_T)}\leq M_2,\label{baru_m-Linfty0TW_01p-Linfty-bound}\\
		&\left\|u_m-\bar{u}_m\right\|_{L^2(Q_T)}^2\leq\frac{1}{m^2}M_2,~ \left\|u_m-\tilde{u}_m\right\|_{L^2(Q_T)}^2\leq\frac{1}{m^2}M_2,\label{bound-u_m-baru_m-tildeum-L2norm}\\
		&
		\left\|\nu^m\right\|_{L^{1}(Q_T;(\mathcal{E}_0^{\gamma}(\mathbb{R}^N))^{\prime})}+\left\|\langle\nu^m,\vec{q}_{\gamma}\rangle\right\|_{L^{\infty}(0,T;L^{\frac{\gamma}{\gamma-1}}(\Omega))}\leq M_2.\label{nu_m-L1Q_TE_0p-nu_m,alpha-Lpprime-bound}
	\end{align}
\end{lemma}
\begin{pf}
	Since $0\leq\lambda_m^j\leq 1$, we have
	\begin{equation}
		\sup_{0\leq t\leq T}\left\| w_m(t,\cdot)\right\|_{W^{s,p}(\Omega)} = \sup_{1\leq j\leq m}\left\|\lambda_m^j(t)w_m^j+\left(1-\lambda_m^j(t)\right)w_m^{j-1}\right\|_{W^{s,p}(\Omega)}\leq 2 M_1.\notag
	\end{equation}
	Considering that $\frac{\partial w_m}{\partial t}=\sum_{j=1}^{m}\frac{m}{T}\chi_m^j(w_m^j-w_m^{j-1})$, we derive
	\begin{equation}
		\left\| \frac{\partial w_m}{\partial t}\right\|_{L^2(Q_T)}^2=\sum_{j=1}^{m}\int_{(j-1)\frac{T}{m}}^{j\frac{T}{m}}\int_{\Omega}\frac{m^2}{T^2}\left(w_m^j-w_m^{j-1}\right)^2 dx dt=\frac{m}{T}\sum_{j=1}^{m}\left\|w_m^j-w_m^{j-1}\right\|_{L^2(\Omega)}^2, \notag
	\end{equation}
	which, when combined with \eqref{estimate w_mj} and \eqref{infty u_mj + w_mj}, yields \eqref{bound-wm-Linftynorm}. Using the result from Lemma \ref{discrete system infty}, we obtain
	\begin{align}
		\left\|w_m-\bar{w}_m\right\|_{L^2(Q_T)}^2&=\sum_{j=1}^{m}\int_{(j-1)\frac{T}{m}}^{j\frac{T}{m}}\int_{\Omega}\left(\frac{m}{T}t-j\right)^2 \left(w_m^j-w_m^{j-1}\right)^2 dxdt\notag\\
		&\leq\sum_{j=1}^{m}\int_{(j-1)\frac{T}{m}}^{j\frac{T}{m}}\int_{\Omega}\left(w_m^j-w_m^{j-1}\right)^2 dxdt\notag\\
		&=\frac{T}{m}\sum_{j=1}^{m}\left\|w_m^j-w_m^{j-1}\right\|_{L^2(\Omega)}^2\leq \frac{T^2}{m^2}M_1.\notag
	\end{align}
	Thus, inequality \eqref{bound-w_m-barw_m-L2norm} is proven. Similarly, by applying the same method, inequalities \eqref{bound-barwm-Linftynorm}, \eqref{u_m-Linfty0TW_01p-Linfty-du_mdt-L2-bound}, \eqref{baru_m-Linfty0TW_01p-Linfty-bound}, and \eqref{bound-u_m-baru_m-tildeum-L2norm} can be proved.
	
	We now proceed to prove \eqref{nu_m-L1Q_TE_0p-nu_m,alpha-Lpprime-bound}. First, we establish that $\langle\nu_x^{m,j},\vec{q}_{\gamma} \rangle=\vec{\varrho}_{\gamma}(\nabla u_m^j(x))$, a.e. $x\in\Omega$. Let $\tilde{\eta}\in C_0^{\infty}(\Omega)$ and $-1<\epsilon <1$. Then, there exists a positive constant $C_4$ such that for $\xi\in\mathbb{R}^N$, $\varphi_{\gamma}^{**}(\xi+\epsilon\nabla\tilde{\eta})\leq C_4(1+|\xi|^2)$, i.e. $\varphi_{\gamma}^{**}(\cdot+\epsilon\nabla\tilde{\eta})\in \mathcal{E}^{\gamma}(\mathbb{R}^N)$. For $j=1,2,\cdots,m$, using the representation formula \eqref{rfYM} along with part (i) of Lemma \ref{kp1}, we derive
	\begin{align}
		&\int_{\Omega}\left(\frac{T}{m}\langle\nu^{m,j},\varphi_{\gamma}^{**}\rangle+\frac{1}{2}\left(u_m^j-u_m^{j-1}\right)^2+\frac{T}{m}\left(\lambda_1 u_m^{j-1}-\lambda_2 w_m^j\right)u_m^j\right)dx\notag\\
		=&\int_{\Omega}\left(\frac{T}{m}\varphi_{\gamma}^{**}(\nabla u_m^j)+\frac{1}{2}\left(u_m^j-u_m^{j-1}\right)^2+\frac{T}{m}\left(\lambda_1 u_m^{j-1}-\lambda_2 w_m^j\right)u_m^j\right)dx\notag\\
		=&E_m^{**}\left(u_m^j;u_m^{j-1}\right) \leq E_m^{**}\left(u_m^j+\epsilon\nabla\tilde{\eta};u_m^{j-1}\right)\notag\\
		=&\int_{\Omega}\left(\frac{T}{m}\varphi_{\gamma}^{**}(\nabla u_m^j+\epsilon\nabla\tilde{\eta})+\frac{1}{2}\left(u_m^j+\epsilon\tilde{\eta}-u_m^{j-1}\right)^2+\frac{T}{m}\left(\lambda_1 u_m^{j-1}-\lambda_2 w_m^j\right)\left(u_m^j+\epsilon \tilde{\eta}\right)\right)dx \notag\\
		\leq& \liminf_{k \rightarrow \infty}\frac{T}{m}\int_{\Omega}\varphi_{\gamma}^{**}(\nabla u_m^{j,k}+\epsilon\nabla\tilde{\eta})dx+\frac{1}{2}\int_{\Omega}\left(u_m^j+\epsilon\tilde{\eta}-u_m^{j-1}\right)^2 dx\notag\\
		&+\frac{T}{m}\int_{\Omega}\left(\lambda_1 u_m^{j-1}-\lambda_2 w_m^j\right)\left(u_m^j+\epsilon \tilde{\eta}\right)dx\notag\\
		=& \int_{\Omega}\left(\frac{T}{m}\langle\nu^{m,j},\varphi_{\gamma}^{**}(\cdot+\epsilon\nabla\tilde{\eta}) \rangle +\frac{1}{2}\left(u_m^j+\epsilon\tilde{\eta}-u_m^{j-1}\right)^2+\frac{T}{m}\left(\lambda_1 u_m^{j-1}-\lambda_2 w_m^j\right)\left(u_m^j+\epsilon \tilde{\eta}\right)\right) dx.\notag
	\end{align}
	Differentiating at $\epsilon=0$ yields
	\begin{equation}
		\frac{T}{m}\int_{\Omega}\langle\nu_x^{m,j},\vec{\varrho}_{\gamma}\rangle\cdot\nabla\tilde{\eta} dx+\int_{\Omega}\left(u_m^j-u_m^{j-1}\right)\tilde{\eta} dx+\frac{T}{m}\int_{\Omega}\left(\lambda_1 u_m^{j-1}-\lambda_2 w_m^j\right)\tilde{\eta} dx=0.\notag
	\end{equation}
	From \eqref{equation u_mj} and the equilibrium equation above, we obtain
	\begin{equation}
		\vec{\varrho}_{\gamma}(\nabla u_m^j(x))=\langle\nu^{m,j}_x,\vec{\varrho}_{\gamma}\rangle=\langle \nu_x^{m,j},\vec{q}_{\gamma}\rangle,\text{ a.e. }x\in\Omega.\label{add1}
	\end{equation}
	Furthermore, we have
	\begin{equation}\label{equilibrumforPsi}
		\int_{\Omega} \left(\langle\nu_x^{m,j},\vec{q}_{\gamma}\rangle\cdot\nabla\eta+\frac{m}{T}\left(u_m^j-u_m^{j-1}\right)\eta+\left(\lambda_1 u_m^{j-1}-\lambda_2 w_m^{j}\right)\eta\right) dx=0,~\forall\eta\in L^2(\Omega)\cap W^{1,\gamma}(\Omega).
	\end{equation}
	For any $t\in [0,T]$, there exists a positive integer $j=1,2,\cdots,m$ such that $(j-1)\frac{T}{m}\leq t<j\frac{T}{m}$. By direct calculation, we derive from \eqref{structure condition 2} and \eqref{estimate u_mj} that
	\begin{align}
		\sup_{0\leq t\leq T} \left\|\langle\nu^m_{t,\cdot},\vec{q}_{\gamma}\rangle\right\|_{L^{\frac{\gamma}{\gamma-1}}(\Omega)}^{\frac{\gamma}{\gamma-1}} & = \sup_{1\leq j\leq m} \int_{\Omega}\left|\langle\nu_x^{m,j},\vec{q}_{\gamma}\rangle\right|^{\frac{\gamma}{\gamma-1}} dx\notag\\
		& = \sup_{1\leq j\leq m}\int_{\Omega}\left|\vec{\varrho}_{\gamma}(\nabla u_m^j)\right|^{\frac{\gamma}{\gamma-1}} dx\notag\\
		& \leq \sup_{1\leq j\leq m}\int_{\Omega}\left(2^{\frac{1}{\gamma-1}}\tilde{C}_{\gamma}^{\frac{\gamma}{\gamma-1}}\left|\nabla u_m^j\right|^{\gamma}+2^{\frac{1}{\gamma-1}}\right) dx\notag\\
		& \leq 2^{\frac{1}{\gamma-1}}\tilde{C}_{\gamma}^{\frac{\gamma}{\gamma-1}} M_1^{\gamma}+2^{\frac{1}{\gamma-1}}\left|\Omega\right|.\label{nu_m,q2}
	\end{align}
	Recalling that $\{\nabla v_m^k\}_{k=1}^{\infty}$ is the $W^{1,\gamma}(Q_T)$-gradient generating sequence of $\nu^m$, and utilizing \eqref{rfYM}, \eqref{iteration u_mjk}, and \eqref{bound varphi u_mj}, we have
	\begin{align}
		\left\| \nu^m\right\|_{L^1(Q_T;(\mathcal{E}_0^{\gamma}(\mathbb{R}^N))^{\prime})}&=\iint_{Q_T}\sup\left\{\left|\langle \nu_{t,x}^m,\tilde{\varphi}\rangle\right|:\tilde{\varphi}\in\mathcal{E}_0^{\gamma}(\mathbb{R}^N),\left\|\tilde{\varphi}\right\|_{\mathcal{E}^{\gamma}(\mathbb{R}^N)}\leq 1\right\}dxdt\notag\\
		&\leq\iint_{Q_T}\int_{\mathbb{R}^N}\left(1+\left|\theta\right|^{\gamma}\right)d\nu_{t,x}^m(\theta)dxdt\notag\\
		&=\iint_{Q_T}\langle \nu_{t,x}^m,\left|\mathrm{id}\right|^{\gamma}\rangle dxdt+T\left|\Omega\right|\notag\\
		&=\lim_{k\rightarrow\infty}\iint_{Q_T}|\nabla v_m^k|^{\gamma} dxdt+T\left|\Omega\right|\notag\\
		&\leq \left(M_1^{\gamma}+\left|\Omega\right|\right)T.\notag
	\end{align}
	Combining this with \eqref{nu_m,q2}, we obtain \eqref{nu_m-L1Q_TE_0p-nu_m,alpha-Lpprime-bound}, completing the proof.
\end{pf}

According to \eqref{Defn-p^s} and \eqref{equilibrumforPsi}, for any $\psi\in L^2(Q_T)\cap L^1(0,T;W^{s,p}(\Omega))$, we derive
\begin{align}
	\iint_{Q_T} & \frac{\partial w_m}{\partial t}\psi dx dt+\frac{1}{2}\iiint_{Q_T^{1,2}}\frac{\left|\bar{w}_m(t,x)-\bar{w}_m(t,y)\right|^{p-2}\left(\bar{w}_m(t,x)-\bar{w}_m(t,y)\right)}{\left|x-y\right|^{N+sp}}\left(\psi(t,x)-\psi(t,y)\right)dx dy dt\notag\\
	& = \lambda_3 \iint_{Q_T}\left(f-K\tilde{u}_m\right)K\psi dx dt,\label{equation w_m}
\end{align}
and for any $\zeta\in L^2(Q_T)\cap  L^{1}(0,T;W^{1,\gamma}(\Omega))$, we have
\begin{equation}\label{equation for u_m}
	\iint_{Q_T} \left(\frac{\partial u_m}{\partial t}\zeta+\langle\nu^m,\vec{q}_{\gamma}\rangle\cdot\nabla\zeta+\left(\lambda_1 \tilde{u}_m-\lambda_2 \bar{w}_m\right)\zeta\right) dx dt = 0.
\end{equation}

Using the results from \eqref{bound-wm-Linftynorm} and \eqref{bound-barwm-Linftynorm} in Lemma \ref{estimates-for-w-m}, together with the Aubin-Lions Lemma, there exist functions $w \in L^\infty(0,T; L^2(\Omega)\cap W^{s,p}(\Omega))\cap L^{\infty}(Q_T)$ with $\frac{\partial w}{\partial t} \in L^2(Q_T)$, and $\bar{w} \in L^{\infty}(0,T; L^2(\Omega)\cap W^{s,p}(\Omega))\cap L^{\infty}(Q_T)$, and subsequences (not relabeled) of $\{w_m\}_{m=1}^{\infty}$ and $\{\bar{w}_m\}_{m=1}^{\infty}$, respectively, such that
\begin{align}
	& w_m\rightarrow w\text{ strongly in } L^p(Q_T),~\frac{\partial w_m}{\partial t}\rightharpoonup \frac{\partial w}{\partial t}\text{ weakly in } L^2(Q_T),\label{convergen w_m}\\
	&w_m\rightharpoonup w \text{ weakly }{\ast} \text{ in } L^{\infty}(0,T; W^{s,p}(\Omega)),~\bar{w}_m\rightharpoonup \bar{w} \text{ weakly }{\ast} \text{ in } L^{\infty}(0,T; W^{s,p}(\Omega)).\label{convergen barw_m}
\end{align}
From \eqref{bound-wm-Linftynorm}, \eqref{bound-w_m-barw_m-L2norm}, and \eqref{convergen w_m}, we deduce that $\bar{w}_m$ converges strongly to $w$ in $L^2(Q_T)$, which implies that $w=\bar{w}$, a.e. $(t,x)\in Q_T$. Due to \eqref{u_m-Linfty0TW_01p-Linfty-du_mdt-L2-bound},  \eqref{baru_m-Linfty0TW_01p-Linfty-bound}, and the classical compactness results in Sobolev spaces, there exist functions $u\in L^{\infty}(0,T;W^{1,\gamma}(\Omega))\cap W^{1,\gamma}(Q_T)\cap L^{\infty}(Q_T)$ with $\frac{\partial u_m}{\partial t}\in L^2(Q_T)$, and $\bar{u},~ \tilde{u}\in L^{\infty}(0,T;W^{1,\gamma}(\Omega))\cap L^{\infty}(Q_T)$, and subsequences can be extracted from $\{u_m\}_{m=1}^{\infty}$, $\{\bar{u}_m\}_{m=1}^{\infty}$, and $\{\tilde{u}_m\}_{m=1}^{\infty}$, respectively, still denoted by themselves, such that
\begin{align}
	& u_m\rightarrow u \text{ strongly in } L^{\gamma}(Q_T),~u_m\rightharpoonup u \text{ weakly in }W^{1,\gamma}(Q_T),~ \frac{\partial u_m}{\partial t}\rightharpoonup\frac{\partial u}{\partial t}\text{ weakly in } L^2(Q_T),\label{convergence-of-u_m}\\
	& \bar{u}_m\rightharpoonup\bar{u} \text{ weakly }{\ast}\text{ in } L^{\infty}(0,T,W^{1,\gamma}(\Omega)),~\tilde{u}_m\rightharpoonup\tilde{u} \text{ weakly }{\ast}\text{ in } L^{\infty}(0,T,W^{1,\gamma}(\Omega)).\label{convergence-of-baru_m-tildeu_m}
\end{align}
Combining \eqref{bound-u_m-baru_m-tildeum-L2norm}, \eqref{convergence-of-u_m}, and \eqref{convergence-of-baru_m-tildeu_m}, it follows that $\bar{u}_m$ and $\tilde{u}_m$ converge strongly to $u$ in $L^2(Q_T)$, leading to $u=\bar{u}=\tilde{u}$, a.e. $(t,x)\in Q_T$. Furthermore, applying the Aubin-Lions Lemma once more, we deduce that $u,w\in C([0,T];L^2(\Omega))$. Since
\begin{equation}
	\left\|K\tilde{u}_m-Ku\right\|^2_{L^2(Q_T)}\leq \left\|K\right\|^2\left|\Omega\right|\left\|\tilde{u}_m-u\right\|^2_{L^2(Q_T)}\leq 2\left\|K\right\|^2\left|\Omega\right|\left\|\tilde{u}_m-u_m\right\|^2_{L^2(Q_T)}+ 2\left\|K\right\|^2\left|\Omega\right|\left\|u_m-u\right\|^2_{L^2(Q_T)},\notag
\end{equation}
we conclude that $K\tilde{u}_m$ converges strongly to $ Ku$ in $L^2(Q_T)$. Likewise, the strong convergence of $K\bar{w}_m$ in $L^2(Q_T)$ can also be established.

On the other hand, by Lemma \ref{wcp}, Remark \ref{rmk1}, and \eqref{nu_m-L1Q_TE_0p-nu_m,alpha-Lpprime-bound}, there exists a $W^{1,\gamma}(Q_T)$-gradient Young measure $\nu$ and a subsequence (not relabeled) of $\{\nu^m\}_{m=1}^{\infty}$ such that
\begin{align}
	&\langle\nu^m,\mathrm{id}\rangle\rightharpoonup\langle\nu,\mathrm{id}\rangle\text{ weakly in }L^{\gamma}(Q_T),~\langle\nu^m,\vec{q_{\gamma}}\rangle\rightharpoonup\langle\nu,\vec{q_{\gamma}}\rangle\text{ weakly }{\ast}\text{ in }L^{\infty}(0,T;L^{\frac{\gamma}{\gamma-1}}(\Omega)),\notag\\
	& \langle\nu^m,\vec{q}_{\gamma}\cdot\mathrm{id} \rangle\rightharpoonup\langle\nu,\vec{q}_{\gamma}\cdot\mathrm{id}\rangle\text{ in the biting sense}.\label{Lq/q-1convergence-for-nu_m,alpha}
\end{align}

Based on the convergence results established above, we now present the following theorem.

\begin{theorem}\label{thm2}
	If $f\in L^{\infty}(\Omega)\cap W^{1,\gamma}(\Omega)$, the system \eqref{New system} admits a unique Young measure solution for $1<p\leq 2$.
\end{theorem}
\begin{pf}
	By \eqref{bound-barwm-Linftynorm}, we have
	\begin{equation}
		\sup_{0\leq t\leq T}\left\|\frac{\left|\bar{w}_m(t,x)-\bar{w}_m(t,y)\right|^{p-2}}{\left|x-y\right|^{(N+sp)\frac{p-1}{p}}}(\bar{w}_m(t,x)-\bar{w}_m(t,y))\right\|_{L^{\frac{p}{p-1}}(\Omega\times\Omega)}\leq \left\| \bar{w}_m\right\|^{p-1}_{L^{\infty}(0,T;W_0^{s,p}(\Omega))}\leq M_2^{p-1}.\notag
	\end{equation}
	Hence, there exists a function $\xi(t,x,y)\in L^{\infty}(0,T;L^{\frac{p}{p-1}}(\Omega\times\Omega))$ such that
	\begin{equation}
		\frac{\left|\bar{w}_m(t,x)-\bar{w}_m(t,y)\right|^{p-2}}{\left|x-y\right|^{(N+sp)\frac{p-1}{p}}}(\bar{w}_m(t,x)-\bar{w}_m(t,y))\rightharpoonup\xi\text{ weakly }{\ast}\text{ in }L^{\infty}(0,T;L^{\frac{p}{p-1}}(\Omega\times\Omega)).\notag
	\end{equation}
	Let $\mathcal{A}:L^{\infty}(0,T;W^{s,p}(\Omega))\rightarrow L^{\infty}(0,T;L^{\frac{p}{p-1}}(\Omega\times\Omega))$ be an operator defined by
	\begin{equation}
		\omega\mapsto\frac{\left|\omega(t,x)-\omega(t,y)\right|^{p-2}}{\left|x-y\right|^{(N+sp)\frac{p-1}{p}}}\left(\omega(t,x)-\omega(t,y)\right).\notag
	\end{equation}
	For any $\psi\in L^2(Q_T)\cap L^1(0,T;W^{s,p}(\Omega))$, taking the limit as $m\rightarrow\infty$ in \eqref{equation w_m}, we have
	\begin{equation}
		\iint_{Q_T}\frac{\partial w}{\partial t}\psi dx dt+\frac{1}{2}\iiint_{Q_T^{1,2}}\xi\frac{\psi(t,x)-\psi(t,y)}{\left|x-y\right|^{(N+sp)\frac{1}{p}}} dx dy dt=\lambda_3\iint_{Q_T}\left(f-Ku\right)K\psi dx dt.\notag
	\end{equation}
	Choosing $\psi=w$ yields
	\begin{equation}\label{equation psi = w}
		\iint_{Q_T}\frac{\partial w}{\partial t} w dx dt+\frac{1}{2}\iiint_{Q_T^{1,2}}\xi\frac{w(t,x)-w(t,y)}{\left|x-y\right|^{(N+sp)\frac{1}{p}}} dx dy dt=\lambda_3\iint_{Q_T}\left(f-Ku\right)Kw dx dt.
	\end{equation}
	Taking $\psi=\bar{w}_m$ in equilibrium equation \eqref{equation w_m}, we get
	\begin{equation}\label{equation psi = bar wm}
		\iint_{Q_T}\frac{\partial w_m}{\partial t}\bar{w}_m dx dt + \frac{1}{2}\iiint_{Q_T^{1,2}}\frac{\left|\bar{w}_m(t,x)-\bar{w}_m(t,y)\right|^p}{\left|x-y\right|^{N+sp}} dx dy dt =\lambda_3 \iint_{Q_T}\left(f-K\tilde{u}_m\right)K\bar{w}_m dx dt.
	\end{equation}
	Note that
	\begin{align}
		&\left| \iint_{Q_T} \left( \left(f - K\tilde{u}_m\right)K\bar{w}_m-\left(f-Ku\right)Kw \right)dxdt \right|\notag\\
		\leq & \left| \iint_{Q_T}\left(  K\bar{w}_m - K w\right) f dxdt \right|+\left|\iint_{Q_T}\left(Ku Kw - K\tilde{u}_m K\bar{w}_m\right)dxdt\right|  \notag\\
		\leq & \left| \iint_{Q_T}\left(  K\bar{w}_m - K w\right) f dxdt \right| + \left| \iint_{Q_T}\left(Ku - K\tilde{u}_m\right)Kw dxdt \right| + \left| \iint_{Q_T}\left(Kw-K\bar{w}_m \right)K\tilde{u}_m dxdt \right|. \notag
	\end{align}
	As $K\bar{w}_m\rightarrow Kw$ and $K\tilde{u}_m\rightarrow Ku$ strongly in $L^2(Q_T)$, it follows that
	\begin{equation}
		\lim_{m\rightarrow\infty}\iint_{Q_T}\left(f-K\tilde{u}_m\right)K\bar{w}_m dx dt = \iint_{Q_T}\left(f-Ku\right)Kw dx dt.\notag
	\end{equation}
	Combining \eqref{convergen w_m}, \eqref{convergence-of-baru_m-tildeu_m}, \eqref{equation psi = w}, and \eqref{equation psi = bar wm}, we derive
	\begin{align}
		&\limsup\limits_{m\rightarrow\infty}\iiint_{Q_T^{1,2}}\mathcal{A}\bar{w}_m\frac{\bar{w}_m(t,x)-\bar{w}_m(t,y)}{\left|x-y\right|^{(N+sp)\frac{1}{p}}} dx dy dt \notag\\
		=& \limsup\limits_{m\rightarrow\infty}\iiint_{Q_T^{1,2}} \frac{\left|\bar{w}_m(t,x)-\bar{w}_m(t,y)\right|^p}{\left|x-y\right|^{N+sp}} dx dy dt\notag\\
		\leq & -2\iint_{Q_T}\frac{\partial w}{\partial t} w dx dt + 2\lambda_3\iint_{Q_T}\left(f-Ku\right) Kw dx dt \notag\\
		= & \iiint_{Q_T^{1,2}}\xi\frac{w(t,x)-w(t,y)}{\left|x-y\right|^{(N+sp)\frac{1}{p}}} dx dy dt.\label{inequality mathcal A w_m}
	\end{align}
	Using the monotonicity of the function $t\mapsto t|t|^{p-2}$ for $p>1$, we have
	\begin{align}
		&\iiint_{Q_T^{1,2}}\mathcal{A}\bar{w}_m
		\left(\frac{\bar{w}_m(t,x)-\bar{w}_m(t,y)}{\left|x-y\right|^{(N+sp)\frac{1}{p}}}-\frac{\psi(t,x)-\psi(t,y)}{\left|x-y\right|^{(N+sp)\frac{1}{p}}}\right)
		dx dy dt\notag\\
		\geq&\iiint_{Q_T^{1,2}}\mathcal{A}\psi \left(\frac{\bar{w}_m(t,x)-\bar{w}_m(t,y)}{\left|x-y\right|^{(N+sp)\frac{1}{p}}}-\frac{\psi(t,x)-\psi(t,y)}{\left|x-y\right|^{(N+sp)\frac{1}{p}}}\right)dx dy dt.\notag
	\end{align}
	By \eqref{inequality mathcal A w_m}, we obtain
	\begin{align}
		\iiint_{Q_T^{1,2}}\left(\xi-\mathcal{A}\psi\right)\left(\frac{w(t,x)-w(t,y)}{\left|x-y\right|^{(N+sp)\frac{1}{p}}}-\frac{\psi(t,x)-\psi(t,y)}{\left|x-y\right|^{(N+sp)\frac{1}{p}}}\right) dx dy dt \geq 0.\notag
	\end{align}
	Then, a standard monotone method leads to
	\begin{equation}\label{xi}
		\xi=\mathcal{A}w=\frac{\left|w(t,x)-w(t,y)\right|^{p-2}}{\left|x-y\right|^{(N+sp)\frac{p-1}{p}}}(w(t,x)-w(t,y
		)).
	\end{equation}
	
	Now, we verify that $(u,w)$ is the Young measure solution to the system \eqref{New system}. Letting $m\rightarrow\infty$ in \eqref{equation w_m}, and using \eqref{convergen w_m}, \eqref{convergence-of-baru_m-tildeu_m}, and \eqref{xi}, we can deduce
	\begin{align}\notag
		\iint_{Q_T} & \frac{\partial w}{\partial t}\psi dx dt+\frac{1}{2}\iiint_{Q_T^{1,2}}\frac{\left|w(t,x)-w(t,y)\right|^{p-2}\left(w(t,x)-w(t,y)\right)}{\left|x-y\right|^{N+sp}}\left(\psi(t,x)-\psi(t,y)\right)dx dy dt\notag\\
		& =\lambda_3 \iint_{Q_T}\left(f-Ku\right)K\psi dx dt,~\forall\psi\in L^2(Q_T)\cap L^1(0,T;W^{s,p}(\Omega)).\notag
	\end{align}
	From \eqref{convergen barw_m}, \eqref{convergence-of-u_m}, \eqref{Lq/q-1convergence-for-nu_m,alpha}, and taking the limit as $m\rightarrow\infty$ in \eqref{equation for u_m}, we have
	\begin{equation}
		\iint_{Q_T} \left(\frac{\partial u}{\partial t}\zeta+\langle\nu,\vec{q}_{\gamma}\rangle\cdot\nabla\zeta+\left(\lambda_1 u-\lambda_2 w\right)\zeta\right) dx dt = 0,~\forall\zeta\in L^2(Q_T) L^{1}(0,T;W^{1,\gamma}(\Omega)).\notag
	\end{equation}
	Since $\nabla \bar{u}_m = \langle\nu^m,\mathrm{id}\rangle$, we derive from Lemma \ref{wcp} and \eqref{convergence-of-baru_m-tildeu_m} that $\langle\nu^m,\mathrm{id}\rangle$ converges weakly to $\langle\nu,\mathrm{id}\rangle$ in $L^2(Q_T)$, and $\nabla \bar{u}_m$ converges weakly to $\nabla \bar{u}$ in $L^2(Q_T)$, which together with $u=\bar{u}$, a.e. $(t,x)\in Q_T$, implies that $\nabla u(t,x)=\langle\nu_{t,x},\mathrm{id}\rangle$, a.e. $(t,x)\in Q_T$. Given that $u_m(0,x)=f(x)$ and $w_m(0,x)=0$, a.e. $x\in\Omega$, and by employing \eqref{convergence-of-u_m}, we obtain that the initial conditions $u(0,x)=f(x)$ and $w(0,x)=0$, a.e. $x\in\Omega$ hold. To prove \eqref{suppfornu}, we define the set $S\subseteq\mathbb{R}^N$ as follows:
	\begin{equation}
		S=\bigcup_{m=1}^{\infty}\bigcup_{(t,x)\in Q_T} \operatorname{supp}\nu_{t,x}^m. \notag
	\end{equation}
	Given \eqref{suppform}, we deduce that $S\subseteq\{\theta\in\mathbb{R}^N:\varphi_{\gamma}^{**}(\theta)=\varphi_{\gamma}(\theta)\}$. For any non-negative function $g\in C_0(\mathbb{R}^N)$ such that $g(\theta)=0$ for $\theta\in S$, and from (i) in Lemma \ref{wcp} and the fact that $\langle\nu^m_{t,x},g\rangle=0$, a.e. $(t,x)\in Q_T$, we have
	\begin{equation}
		\iint_{Q_T}\int_{\mathbb{R}^N}g(\theta)d\nu_{t,x}(\theta)dxdt=\lim_{m\rightarrow\infty}\iint_{Q_T}\int_{\mathbb{R}^N}g(\theta)d\nu^m_{t,x}(\theta)dxdt=0.\notag
	\end{equation}
	Therefore, $\langle\nu_{t,x},g\rangle=0$, a.e. $(t,x)\in Q_T$, implying that
	\begin{equation}
		\operatorname{supp}\nu_{t,x}\subseteq\{\theta\in\mathbb{R}^N:\varphi_{\gamma}^{**}(\theta)=\varphi_{\gamma}(\theta)\},\text{ a.e. } (t,x)\in Q_T.\notag
	\end{equation}
	
	Next, we prove the independence property \eqref{independentfornu}. From \eqref{indenpendence nu_mj} and the definition of the Young measure $\nu^m$, we obtain 
	\begin{equation}
		\langle\nu_{t,x}^{m},\vec{q}_{\gamma}\cdot\mathrm{id}\rangle=\langle\nu_{t,x}^{m},\vec{q}_{\gamma}\rangle\cdot\langle\nu_{t,x}^{m},\mathrm{id}\rangle,\text{ a.e. }(t,x)\in Q_T.\notag
	\end{equation}
	We proceed to analyze the convergence of $\langle\nu^{m},\vec{q}_{\gamma}\rangle\cdot\langle\nu^{m},\mathrm{id}\rangle$. Choosing $\zeta = \bar{u}_m\tilde{\zeta}$ for $\tilde{\zeta}\in C_0^{\infty}(Q_T)$ as a test function in \eqref{equation-for-u-nu} and \eqref{equation for u_m}, we derive
	\begin{align}
		&\left|\iint_{Q_T}\left(\langle\nu^m,\vec{q}_{\gamma}\rangle\cdot\langle\nu^m,\mathrm{id}\rangle-\langle\nu,\vec{q}_{\gamma}\rangle\cdot\langle\nu,\mathrm{id}\rangle\right)\tilde{\zeta} dxdt\right|\notag\\
		\leq&\left| \iint_{Q_T}\left( \langle \nu^m,\vec{q}_{\gamma} \rangle - \langle \nu,\vec{q}_{\gamma} \rangle \right)\cdot\langle \nu^m,\mathrm{id} \rangle \tilde{\zeta} dxdt \right|+\left| \iint_{Q_T}\left( \langle \nu^m,\mathrm{id} \rangle - \langle \nu,\mathrm{id} \rangle \right)\cdot\langle \nu,\vec{q}_{\gamma} \rangle \tilde{\zeta} dxdt \right|\notag\\
		\leq & \left| \iint_{Q_T}\left( \langle \nu^m,\vec{q}_{\gamma} \rangle - \langle \nu,\vec{q}_{\gamma} \rangle \right)\cdot\nabla\left(\bar{u}_m \tilde{\zeta}\right) dxdt \right|+\left| \iint_{Q_T}\bar{u}_m\left( \langle \nu^m,\vec{q}_{\gamma} \rangle - \langle \nu,\vec{q}_{\gamma} \rangle \right)\cdot\nabla\tilde{\zeta} dxdt \right| \notag\\
		&+\left| \iint_{Q_T}\left( \langle \nu^m,\mathrm{id} \rangle - \langle \nu,\mathrm{id} \rangle \right)\cdot\langle \nu,\vec{q}_{\gamma} \rangle \tilde{\zeta} dxdt \right|\notag\\
		\leq & \left|\iint_{Q_T}\bar{u}_m\left(\frac{\partial u}{\partial t}-\frac{\partial u_m}{\partial t}\right)\tilde{\zeta}dxdt\right|+\lambda_1\left|\iint_{Q_T}\bar{u}_m\left(u-\tilde{u}_m\right)\tilde{\zeta} dxdt\right| +\lambda_2 \left|\iint_{Q_T}\bar{u}_m\left(\bar{w}_m-w\right)\tilde{\zeta} dxdt\right| \notag\\
		& + \left| \iint_{Q_T}\bar{u}_m\left( \langle \nu^m,\vec{q}_{\gamma} \rangle - \langle \nu,\vec{q}_{\gamma} \rangle \right)\cdot\nabla\tilde{\zeta} dxdt \right|+\left| \iint_{Q_T}\left( \langle \nu^m,\mathrm{id} \rangle - \langle \nu,\mathrm{id} \rangle \right)\cdot\langle \nu,\vec{q}_{\gamma} \rangle \tilde{\zeta} dxdt \right| \notag\\
		\leq & \left|\iint_{Q_T}u\left(\frac{\partial u}{\partial t}-\frac{\partial u_m}{\partial t}\right)\tilde{\zeta}dxdt\right| + \left|\iint_{Q_T}\left(\bar{u}_m-u\right)\left(\frac{\partial u}{\partial t}-\frac{\partial u_m}{\partial t}\right)\tilde{\zeta}dxdt\right| +\lambda_1\left|\iint_{Q_T}\bar{u}_m\left(u-\tilde{u}_m\right)\tilde{\zeta} dxdt\right| \notag\\
		& + \lambda_2 \left|\iint_{Q_T}\bar{u}_m\left(\bar{w}_m-w\right)\tilde{\zeta} dxdt\right|+\left| \iint_{Q_T}u\left( \langle \nu^m,\vec{q}_{\gamma} \rangle - \langle \nu,\vec{q}_{\gamma} \rangle \right)\cdot\nabla\tilde{\zeta} dxdt \right| \notag\\
		& + \left| \iint_{Q_T}\left(\bar{u}_m-u\right)\left( \langle \nu^m,\vec{q}_{\gamma} \rangle - \langle \nu,\vec{q}_{\gamma} \rangle \right)\cdot\nabla\tilde{\zeta} dxdt \right|+ \left| \iint_{Q_T}\left( \langle \nu^m,\mathrm{id} \rangle - \langle \nu,\mathrm{id} \rangle \right)\cdot\langle \nu,\vec{q}_{\gamma} \rangle \tilde{\zeta} dxdt \right|.\notag
	\end{align}
	From \eqref{bound-w_m-barw_m-L2norm}, \eqref{bound-u_m-baru_m-tildeum-L2norm}, \eqref{convergen w_m}, \eqref{convergence-of-u_m}, and \eqref{Lq/q-1convergence-for-nu_m,alpha}, and applying Hölder's inequality, it follows that
	\begin{equation}
		\lim_{m\rightarrow\infty}\left|\iint_{Q_T}\left(\langle\nu^m,\vec{q}_{\gamma}\rangle\cdot\langle\nu^m,\mathrm{id}\rangle-\langle\nu,\vec{q}_{\gamma}\rangle\cdot\langle\nu,\mathrm{id}\rangle\right)\tilde{\zeta} dxdt\right|=0,~\forall\tilde{\zeta}\in C_0^{\infty}(Q_T).\notag
	\end{equation}
	By the weak $\ast$ density of $C_0^{\infty}(Q_T)$ in $L^{\infty}(Q_T)$, we have
	\begin{equation}
		\langle\nu^m,\vec{q}_{\gamma}\rangle\cdot\langle\nu^m,\mathrm{id}\rangle \rightharpoonup \langle\nu,\vec{q}_{\gamma}\rangle\cdot\langle\nu,\mathrm{id}\rangle\text{ weakly in }L^1(Q_T).\notag
	\end{equation}
	Since $\langle\nu^m,\vec{q}_{\gamma}\cdot\mathrm{id} \rangle$ converges to $\langle\nu,\vec{q}_{\gamma}\cdot\mathrm{id}\rangle$ in the biting sense (see Lemma \ref{wcp}), there exists a decreasing sequence of subsets $F_{j+1}\subset F_j$ of $Q_T$ with $\lim_{j \rightarrow \infty}|F_j|=0$ such that
	\begin{equation}
		\langle\nu_{t,x},\vec{q}_{\gamma}\cdot\mathrm{id}\rangle=\langle\nu_{t,x},\vec{q}_{\gamma}\rangle\cdot\langle\nu_{t,x},\mathrm{id}\rangle,\text{ a.e. }(t,x)\in Q_T\backslash F_j,~j=1,2,\cdots.\notag
	\end{equation}
	From (iii) in Lemma \ref{criteria W1p Young}, we know that $\langle\nu,\vec{q}_{\gamma}\cdot\mathrm{id}\rangle\in L^1(Q_T)$. Then, we have
	\begin{equation}
		\langle\nu_{t,x},\vec{q}_{\gamma}\cdot\mathrm{id}\rangle=\langle\nu_{t,x},\vec{q}_{\gamma}\rangle\cdot\langle\nu_{t,x},\mathrm{id}\rangle, \text{ a.e. }(t,x)\in Q_T.\notag
	\end{equation}
	Therefore, the couple of functions $(u,w)$ is the desired Young measure solution to the reaction-diffusion system \eqref{New system}.
	In the following, we focus on the uniqueness and continuous dependence of the Young measure solutions. It is worth noting that uniqueness pertains only to the functions $u$ and $w$, while the Young measure $\nu=(\nu_{t,x})_{(t,x)\in Q_T}$ is generally not unique. 
	
	Suppose $(u_1,w_1)$ and $(u_2,w_2)$ are two Young measure solutions to the problem \eqref{New system} with the
	initial values $f_{10}$ and $f_{20}$, and let $\nu^1$ and $\nu^2$ be the $W^{1,\gamma}(Q_T)$-gradient Young measures with respect to the functions $u_1$ and $u_2$, respectively, i.e., $\nabla u_i(t,x)=\langle\nu_{t,x}^{i},\mathrm{id}\rangle$, a.e. $(t,x)\in Q_T$, for $i=1,2$. 
	For
	any $s\in [0,T]$, substituting the test functions $\zeta(t,x)=(u_1(t,x)-u_2(t,x))\chi_{[0,s]}(t)$ and $\psi(t,x)=(w_1(t,x)-w_2(t,x))\chi_{[0,s]}(t)$ into equations \eqref{equation-for-u-nu} and \eqref{equation-w}, respectively, we obtain
	\begin{align}
		\iint_{Q_s}&\left( \frac{\partial u_1}{\partial t}-\frac{\partial u_2}{\partial t} \right)\left(u_1-u_2\right) dxdt + \iint_{Q_s}\left( \langle\nu^1,\vec{q}_{\gamma}\rangle - \langle\nu^2,\vec{q}_{\gamma}\rangle\right)\cdot\left(\nabla u_1 - \nabla u_2\right) dxdt\notag\\
		& +\lambda_1 \iint_{Q_s} \left(u_1 - u_2 \right)^2 dxdt =\lambda_2 \iint_{Q_s} \left(u_1 - u_2\right)\left( w_1 - w_2 \right) dxdt, \notag
	\end{align}
	and
	\begin{align}
		\iint_{Q_s}&\left( \frac{\partial w_1}{\partial t} - \frac{\partial w_2}{\partial t}\right)\left( w_1 - w_2 \right) dxdt+\frac{1}{2}\iiint_{Q_s^{1,2}}
		\frac{1}{\left|x-y
			\right|^{N+sp}}\left[ \left|w_1(t,x)-w_1(t,y)\right|^{p-2}\left(w_1(t,x)-w_1(t,y)\right) \right. \notag\\
		& \left. -\left|w_2(t,x)-w_2(t,y)\right|^{p-2}\left(w_2(t,x)-w_2(t,y)\right)
		\right]\left( \left(w_1(t,x)-w_1(t,y)\right)-\left(w_2(t,x)-w_2(t,y)\right) \right)dxdydt\notag\\
		&+\lambda_3 \iint_{Q_s} \left(K u_1 - K u_2 \right)\left(K w_1 - K w_2 \right) dxdt = \lambda_3 \iint_{Q_s} \left( f_{10} - f_{20} \right)\left(K w_1 - K w_2 \right)dxdt.\notag
	\end{align}
	By combining the two equalities above, we have
	\begin{align}
		&\frac{1}{2}\int_{\Omega}\left(u_1(s,x)-u_2(s,x) \right)^2 dx -\frac{1}{2}\int_{\Omega}\left(f_{10}-f_{20} \right)^2 dx +\frac{1}{2} \int_{\Omega}\left(w_1(s,x)-w_2(s,x) \right)^2 dx\notag\\
		\leq & -\iint_{Q_s}\left( \langle\nu^1,\vec{q}_{\gamma}\rangle - \langle\nu^2,\vec{q}_{\gamma}\rangle\right)\cdot\left(\nabla u_1 - \nabla u_2\right) dxdt+\lambda_2 \iint_{Q_s}\left(u_1-u_2\right)\left(w_1-w_2\right)dxdt \notag\\
		&+ \lambda_3 \iint_{Q_s}\left(Ku_2-Ku_1\right)\left(Kw_1-Kw_2\right) dxdt+\lambda_3 \iint_{Q_s}\left(f_{10}-f_{20}\right)\left(Kw_1-Kw_2\right) dxdt. \label{unique pf}
	\end{align}
	From \eqref{gradientforu}--\eqref{suppfornu}, along with the convexity of the function $\varphi_{\gamma}^{**}$, we deduce
	\begin{align}
		&\iint_{Q_s}\left( \langle\nu^1,\vec{q}_{\gamma}\rangle - \langle\nu^2,\vec{q}_{\gamma}\rangle\right)\cdot\left(\nabla u_1 - \nabla u_2\right) dxdt \notag\\
		= & \iint_{Q_s} \left(\langle\nu^1,\vec{q}_{\gamma}\cdot\mathrm{id}\rangle-\langle\nu^1,\vec{q}_{\gamma}\rangle\cdot\langle\nu^2,\mathrm{id}\rangle-\langle\nu^2,\vec{q}_{\gamma}\rangle\cdot\langle\nu^1,\mathrm{id}\rangle+\langle\nu^2,\vec{q}_{\gamma}\cdot\mathrm{id}\rangle\right) dxdt \notag\\
		= & \iint_{Q_s} \left( \int_{\mathbb{R}^N}\vec{q}_{\gamma}(\hat{\theta})\cdot\lambda d\nu_{t,x}^1(\hat{\theta})\int_{\mathbb{R}^N}d\nu^2_{t,x}(\theta) - \int_{\mathbb{R}^N}\vec{q}_{\gamma}(\hat{\theta})d\nu_{t,x}^1(\hat{\theta})\cdot\int_{\mathbb{R}^N}\xi d\nu_{t,x}^2(\theta)\right.\notag\\
		&\left. - \int_{\mathbb{R}^N}\vec{q}_{\gamma}(\theta)d\nu_{t,x}^2(\theta)\cdot\int_{\mathbb{R}^N}\hat{\theta} d\nu_{t,x}^1(\hat{\theta}) + \int_{\mathbb{R}^N}\vec{q}_{\gamma}(\theta)\cdot\theta d\nu_{t,x}^2(\theta)\int_{\mathbb{R}^N}d\nu^1_{t,x}(\hat{\theta})\right)dxdt \notag\\
		= & \iint_{Q_s} \int_{\mathbb{R}^N}\int_{\mathbb{R}^N} \left( \vec{q}_{\gamma}(\hat{\theta})-\vec{q}_{\gamma}(\theta) \right)\cdot\left(\hat{\theta}-\theta\right) d\nu_{t,x}^1(\hat{\theta})d\nu^2_{t,x}(\theta) dxdt \notag\\
		= & \iint_{Q_s} \int_{\mathbb{R}^N}\int_{\mathbb{R}^N} \left( \vec{\varrho}_{\gamma}(\hat{\theta})-\vec{\varrho}_{\gamma}(\theta) \right)\cdot\left(\hat{\theta}-\theta\right) d\nu_{t,x}^1(\hat{\theta})d\nu^2_{t,x}(\theta) dxdt \notag\\
		= &\iint_{Q_s} \int_{\mathbb{R}^N}\int_{\mathbb{R}^N} \left( \nabla\varphi_{\gamma}^{**}(\hat{\theta})-\nabla\varphi_{\gamma}^{**}(\theta) \right)\cdot\left(\hat{\theta}-\theta\right) d\nu_{t,x}^1(\hat{\theta})d\nu^2_{t,x}(\theta) dxdt\geq 0.\label{nu monotone}
	\end{align}
	Substituting \eqref{nu monotone} into \eqref{unique pf} and applying Cauchy's inequality, we have 
	\begin{align}
		&\int_{\Omega}\left(u_1(s,x)-u_2(s,x) \right)^2 dx + \int_{\Omega}\left(w_1(s,x)-w_2(s,x) \right)^2 dx\notag\\
		\leq& \int_{\Omega}\left(f_{10}-f_{20}\right)^2 dx+2\lambda_2 \iint_{Q_s} \left(u_1 - u_2\right)\left( w_1 - w_2 \right) dxdt + 2\lambda_3\iint_{Q_s} \left(K u_2- Ku_1 \right)\left(K w_1- K w_2\right) dxdt\notag\\
		&+2\lambda_3\iint_{Q_s}\left(f_{10}-f_{20}\right)\left(Kw_1-Kw_2\right) dxdt \notag\\
		\leq & \left(\lambda_2+\lambda_3\left\|K \right\|^2\left|\Omega\right|\right) \iint_{Q_s}\left(u_1-u_2\right)^2 dxdt+\left(\lambda_2+2\lambda_3\left\|K\right\|^2\left|\Omega\right|\right)\iint_{Q_s}\left(w_1-w_2\right)^2 dxdt \notag\\
		&+ \left(1+\lambda_3 T\right)\int_{\Omega}\left(f_{10}-f_{20}\right)^2 dx\notag\\
		\leq & C \int_{0}^{s}\left( \int_{\Omega}\left(u_1(t,x)-u_2(t,x)\right)^2 dx + \int_{\Omega}\left(w_1(t,x)-w_2(t,x)\right)^2 dx \right) dt+C\int_{\Omega}\left(f_{10}-f_{20}\right)^2 dx,\notag
	\end{align}
	where $C>0$ depends only on $T$, $\Omega$, $K$, $\lambda_2$, and $\lambda_3$. Thus, using Gronwall's inequality, we obtain
	\begin{equation}\label{continuous dependence}
		\left\|u_1(s,x)-u_2(s,x) \right\|_{L^2(\Omega)}^2 + \left\|w_1(s,x)-w_2(s,x) \right\|_{L^2(\Omega)}^2 \leq C\left(1+Cs e^{Cs}\right)\left\|f_{10}-f_{20}\right\|_{L^2(\Omega)}^2 ,~s\in[0,T],
	\end{equation}
	which shows the continuous dependence of the Young measure solution on the initial data in the $L^2$ norm. When $f_{10}=f_{20}$ is used in \eqref{continuous dependence}, we have $u_1(t,x)=u_2(t,x)$, $w_1(t,x)=w_2(t,x)$, a.e. $(t,x)\in Q_T$, proving the uniqueness of the Young measure solution to the system \eqref{New system} for $1<p\leq 2$. This completes the proof.
\end{pf}

We next consider the relaxed problem associated with \eqref{New system} to better understand the relationship between Young measure solutions, as a type of generalized solution, and classical weak solutions.
\begin{cor}\label{cor1}
	The pair of functions $(u,w)$ obtained in Theorem \ref{thm2} is the unique weak solution of the following relaxed problem
	\begin{equation}\label{relaxed problem p>1}
		\left\{\begin{array}{ll}  
			\displaystyle\frac{\partial u}{\partial t}=\mathrm{div}\left(\vec{\varrho}_{\gamma}(\nabla u)\right)-\lambda_1 u + \lambda_2 w, & x\in \Omega,~0<t<T,\\
			\displaystyle\frac{\partial w}{\partial t} = \Delta_p^s w - \lambda_3 K^{\prime}\left(Ku-f\right), & x\in \Omega,~0<t<T,\\
			\displaystyle u(0,x)=f,~w(0,x)=0, & x\in\Omega,\\
			\displaystyle \vec{\varrho}_{\gamma}(\nabla u)\cdot\vec{n}=0, & x\in \partial\Omega,~0<t<T.\\
		\end{array} \right.
	\end{equation}
\end{cor}
\begin{pf}
	As proved in Theorem \ref{thm2}, $w$ satisfies the integral equality \eqref{equation-w}, so it suffices to focus on the first equation of system \eqref{relaxed problem p>1} in the following proof. From \eqref{baru_m-Linfty0TW_01p-Linfty-bound}, we have
	\begin{equation}
		\sup_{0\leq t\leq T} \left\|\vec{\varrho}_{\gamma}(\nabla \bar{u}_m(t,\cdot))\right\|_{L^{\frac{\gamma}{\gamma-1}}(\Omega)}^{\frac{\gamma}{\gamma-1}}\leq \sup_{0\leq t\leq T}\int_{\Omega} \left( 2^{\frac{1}{\gamma-1}}\tilde{C}_{\gamma}^{\frac{\gamma}{\gamma-1}}\left|\nabla\bar{u}_m(t,x)\right|^{\gamma}+2^{\frac{1}{\gamma-1}} \right) dx \leq 2^{\frac{1}{\gamma-1}}\tilde{C}_{\gamma}^{\frac{\gamma}{\gamma-1}} M_2^{\gamma} + 2^{\frac{1}{\gamma-1}}\left|\Omega\right|.\notag
	\end{equation}
	Thus, there exists a function $\Lambda(t,x)\in L^{\infty}(0,T;L^{\frac{\gamma}{\gamma-1}}(\Omega))$ such that
	\begin{equation}
		\vec{\varrho}_{\gamma}(\nabla\bar{u}_m)\rightharpoonup \Lambda \text{ weakly }{\ast}\text{ in }L^{\infty}(0,T;L^{\frac{\gamma}{\gamma-1}}(\Omega)).  \notag
	\end{equation}
	Based on the previously introduced approximation solutions, equilibrium equation \eqref{equation u_mj}, and \eqref{add1}, we deduce for any $\zeta \in L^2(Q_T) \cap L^1(0,T;W^{1,\gamma}(\Omega))$ that
	\begin{equation}\label{equation varrho bar u_m}
		\iint_{Q_T}\left(\vec{\varrho}_{\gamma}(\nabla \bar{u}_m)\cdot\nabla\zeta + \frac{\partial u_m}{\partial t}\zeta + \left(\lambda_1 \tilde{u}_m-\lambda_2 \bar{w}_m\right)\zeta \right)dxdt=0.
	\end{equation}
	Letting $m\rightarrow\infty$, we have
	\begin{equation}\label{Lambda equation}
		\iint_{Q_T}\left( \Lambda\cdot\nabla\zeta + \frac{\partial u}{\partial t}\zeta + \left(\lambda_1 u-\lambda_2 w\right)\zeta \right)dxdt=0.
	\end{equation}
	Taking $\zeta=u$, we obtain
	\begin{equation}
		\iint_{Q_T}\left( \Lambda\cdot\nabla u + \frac{\partial u}{\partial t}u + \left(\lambda_1 u-\lambda_2 w\right)u \right)dxdt=0.\notag
	\end{equation}
	Furthermore, choosing $\zeta=\bar{u}_m$ in \eqref{equation varrho bar u_m} leads to
	\begin{equation}\label{equation varrho bar u_m bar u_m}
		\iint_{Q_T}\left(\vec{\varrho}_{\gamma}(\nabla \bar{u}_m)\cdot\nabla \bar{u}_m + \frac{\partial u_m}{\partial t}\bar{u}_m + \left(\lambda_1 \tilde{u}_m-\lambda_2 \bar{w}_m\right)\bar{u}_m \right)dxdt=0.
	\end{equation}
	Combining \eqref{Lambda equation} and \eqref{equation varrho bar u_m bar u_m}, we deduce
	\begin{align}
		\iint_{Q_T}\Lambda\cdot\nabla u dxdt=&-\iint_{Q_T}\frac{\partial u}{\partial t}u dxdt-\iint_{Q_T}\left(\lambda_1 u-\lambda_2 w\right)u dxdt\notag\\
		= &  -\lim_{m\rightarrow\infty}\iint_{Q_T}\frac{\partial u_m}{\partial t}\bar{u}_m dxdt-\lim_{m\rightarrow\infty}\iint_{Q_T}\left(\lambda_1 \tilde{u}_m-\lambda_2 \bar{w}_m\right)\bar{u}_m dxdt \notag\\
		= & \lim_{m\rightarrow\infty}\iint_{Q_T}\vec{\varrho}_{\gamma}(\nabla \bar{u}_m)\cdot\nabla\bar{u}_m dxdt.\notag
	\end{align}
	Given that $\varphi_{\gamma}^{**}$ is convex, for any $\psi\in L^{\infty}(0,T;W^{1,\gamma}(\Omega))$, we have
	\begin{equation}
		\iint_{Q_T}\left( \vec{\varrho}_{\gamma}(\nabla \bar{u}_m)-\vec{\varrho}_{\gamma}(\nabla \psi)\right)\cdot\left(\nabla\bar{u}_m-\nabla\psi\right) dxdt\geq 0.\notag
	\end{equation}
	Taking the limit as $m\rightarrow\infty$ in the inequality above, and applying \eqref{convergence-of-baru_m-tildeu_m}, we obtain
	\begin{equation}
		\iint_{Q_T}\left(\Lambda-\vec{\varrho}_{\gamma}(\nabla\psi)\right)\cdot\left(\nabla u-\nabla\psi\right) dxdt\geq 0.\notag
	\end{equation}
	Then, using the standard monotone method, we derive that $\Lambda = \vec{\varrho}_{\gamma}(\nabla u)$, and substituting this into \eqref{Lambda equation} yields
	\begin{equation}
		\iint_{Q_T}\left(\vec{\varrho}_{\gamma}(\nabla u)\cdot\nabla\zeta + \frac{\partial u}{\partial t}\zeta + \left(\lambda_1 u-\lambda_2 w\right)\zeta \right)dxdt=0,\notag
	\end{equation}
	for every $\zeta\in L^2(Q_T)\cap L^{1}(0,T;W^{1,\gamma}(\Omega))$. The Neumann boundary condition is implicitly incorporated in the weak variational formulation $E_m^{**}$ and can be recovered via integration by parts. Therefore, the pair of functions $(u,w)$ derived from Theorem \ref{thm2} is a weak solution to problem \eqref{relaxed problem p>1}, and the uniqueness of the weak solution has been established in \eqref{continuous dependence}.

	Furthermore, it follows from \eqref{add1} that $\vec{\varrho}_{\gamma}(\nabla \bar{u}_m(t,x))=\langle\nu^m_{t,x},\vec{\varrho}_{\gamma}\rangle=\langle\nu^m_{t,x},\vec{q}_{\gamma}\rangle$ for a.e. $(t,x)\in Q_T$, and Lemma \ref{wcp} shows that $\langle\nu^m,\vec{\varrho}_{\gamma}\rangle$ converges weakly to $\langle\nu,\vec{\varrho}_{\gamma}\rangle$ in $L^{1}(Q_T)$, where $\nu$ is the $W^{1,\gamma}(Q_T)$-gradient Young measure obtained in Theorem \ref{thm2}. Then, we derive from \eqref{Lq/q-1convergence-for-nu_m,alpha} that
	\begin{equation}
		\vec{\varrho}_{\gamma}(\nabla u(t,x))=\langle\nu_{t,x},\vec{\varrho}_{\gamma}\rangle=\langle\nu_{t,x},\vec{q}_{\gamma}\rangle,~\text{a.e. }(t,x)\in Q_T.\label{add2}
	\end{equation}
	By combining \eqref{add2} with \eqref{equation-for-u-nu}, we can deduce that $(u,w)$ with $\nu$ is also a Young measure solution to the relaxed problem \eqref{relaxed problem p>1} in the sense of Definition \ref{defn solution system}, and that the functions $u$ and $w$ are unique.
	%Combining \eqref{add2} with \eqref{equation-for-u-nu}, we can deduce that $(u,w)$ with $\nu$ is also the unique Young measure solution of the relaxed problem \eqref{relaxed problem p>1}.

\end{pf}
\subsection{The critical case $p = 1$}\label{sec3.2}
In this section, we investigate the existence and uniqueness of Young measure solutions to problem \eqref{New system} for $p = 1$ by introducing a regularized problem involving the standard fractional $p$-Laplacian diffusion term ($1 < p < 2$), following the approximation method in \cite{GAO2022Fractional,MAZON2016}.
\begin{theorem}\label{thm3}
	If $f\in L^{\infty}(\Omega)\cap W^{1,\gamma}(\Omega)$, the system \eqref{New system} admits a unique Young measure solution for $p=1$.
\end{theorem}
\begin{pf}
	For any $0 < s < 1$ and $1 < \rho < \frac{N}{N + s - 1}$, take $s_{\rho} := s + N(1 - \frac{1}{\rho})$. Then we have $N + s_{\rho} \rho = (N + s)\rho$, and $s < s_{\rho} < 1$. To analyze the well-posedness of problem \eqref{New system}, we begin with the following regularized problem, obtained by approximating the diffusion term $\Delta_1^s w$
	\begin{equation}\label{regularized problem}
		\left\{\begin{array}{ll}  
			\displaystyle\frac{\partial u}{\partial t}=\mathrm{div}\left(\vec{q}_{\gamma}(\nabla u)\right)-\lambda_1 u + \lambda_2 w, & x\in \Omega,~0<t<T,\\
			\displaystyle\frac{\partial w}{\partial t} = \Delta_{\rho}^{s_{\rho}} w - \lambda_3 K^{\ast}\left(Ku-f\right), & x\in \Omega,~0<t<T,\\
			\displaystyle u(0,x)=f,~w(0,x)=0, & x\in\Omega,\\
			\displaystyle \frac{\partial u}{\partial \vec{n}}(t,x)=0, & x\in \partial\Omega,~0<t<T.\\
		\end{array} \right.
	\end{equation}
	By Theorem \ref{thm2}, there exists a unique Young measure solution $(u_{\rho},w_{\rho})$ to problem \eqref{regularized problem}, associated with a $W^{1,\gamma}(Q_T)$-gradient Young measure $\nu^{\rho}=(\nu^{\rho}_{t,x})_{(t,x)\in Q_T}$, such that
	\begin{align}\notag
		\iint_{Q_T} & \frac{\partial w_{\rho}}{\partial t}\psi dx dt+\frac{1}{2}\iiint_{Q_T^{1,2}}\frac{1}{\left|x-y\right|^{(N+s) \rho}}\left|w_{\rho}(t,y)-w_{\rho}(t,x)\right|^{\rho-2}\left(w_{\rho}(t,y)-w_{\rho}(t,x)\right)\left(\psi(t,y)\right.\notag\\
		&\left.-\psi(t,x)\right)dx dy dt =\lambda_3 \iint_{Q_T}\left(f-K u_{\rho}\right)K\psi dx dt\label{equation-w_p}
	\end{align}
	for every $\psi\in L^2(Q_T)\cap L^1(0,T;W^{s_{\rho},\rho}(\Omega))$, and
	\begin{equation}
		\iint_{Q_T} \left(\frac{\partial u_{\rho}}{\partial t}\zeta+\langle\nu^{\rho},\vec{q}_{\gamma}\rangle\cdot\nabla\zeta+\left(\lambda_1 u_{\rho}-\lambda_2 w_{\rho}\right)\zeta\right) dx dt = 0\label{equation-u_p}
	\end{equation}
	for every $\zeta\in L^2(Q_T) L^{1}(0,T;W^{1,\gamma}(\Omega))$. Furthermore, we have
	\begin{align}
		&\nabla u_{\rho}(t,x)=\langle\nu^{\rho}_{t,x},\mathrm{id}\rangle,~\text{a.e }(t,x)\in Q_T,\notag\\
		&\langle \nu^{\rho}_{t,x},\vec{q}_{\gamma}\cdot\mathrm{id}\rangle=\langle \nu^{\rho}_{t,x},\vec{q}_{\gamma}\rangle\cdot\langle \nu^{\rho}_{t,x},\mathrm{id}\rangle,~\text{a.e }(t,x)\in Q_T,\notag\\
		&\mathrm{supp} \nu^{\rho}_{t,x}\subseteq \{\theta\in\mathbb{R}^N:\varphi_{\gamma}(\theta)=\varphi_{\gamma}^{**}(\theta)\},~\text{a.e }(t,x)\in Q_T,\notag
	\end{align}
	and \eqref{add2} yields
	\begin{equation}
		\vec{\varrho}_{\gamma}(\nabla u_{\rho}(t,x))=\langle\nu^{\rho}_{t,x},\vec{\varrho}_{\gamma}\rangle=\langle\nu^{\rho}_{t,x},\vec{q}_{\gamma}\rangle,~\text{a.e. }(t,x)\in Q_T.\label{add3}
	\end{equation}
	In addition, Lemma \ref{estimates-for-w-m} implies that there exists a positive constant $M$, independent of $s$ and $\rho$, satisfying
	\begin{align}
		&	\left\|u_{\rho}\right\|_{L^{\infty}(0,T;W^{1,\gamma}(\Omega))}+\left\|u_{\rho}\right\|_{L^{\infty}(Q_T)}+\left\|\frac{\partial u_{\rho}}{\partial t}\right\|_{L^2(Q_T)}\leq M,\notag\\
		&\left\|w_{\rho}\right\|_{L^{\infty}(0,T; W^{s_{\rho},\rho}(\Omega))}+\left\|w_{\rho}\right\|_{L^{\infty}(Q_T)}+\left\|\frac{\partial w_{\rho}}{\partial t} \right\|_{L^{2}(Q_T)}\leq M,\notag\\
		&\left\|\nu^{\rho}\right\|_{L^{1}(Q_T;(\mathcal{E}_0^{\gamma}(\mathbb{R}^N))^{\prime})}+\left\|\langle\nu^{\rho},\vec{q}_{\gamma}\rangle\right\|_{L^{\infty}(0,T;L^{\frac{\gamma}{\gamma-1}}(\Omega))}\leq M.\notag
	\end{align}
	Thus, there exists a function $u\in L^{\infty}(Q_T)\cap L^{\infty}(0,T;W^{1,\gamma}(\Omega))\cap C([0,T];L^{\gamma}(\Omega))$, with $\frac{\partial u}{\partial t}\in L^2(Q_T)$, and a monotonically decreasing sequence $\{\rho_n\}_{n=1}^{\infty}$ that converges to $1$, such that as $n\rightarrow\infty$,
	\begin{equation}
		u_{\rho_n}\rightharpoonup u \text{ weakly }{\ast}\text{ in }L^{\infty}(0,T;W^{1,\gamma}(\Omega)),~u_{\rho_n}\rightharpoonup u\text{ weakly in }W^{1,\gamma}(Q_T),~\frac{\partial u_{\rho_n}}{\partial t}\rightharpoonup \frac{\partial u}{\partial t}\text{ weakly in }L^{2}(Q_T).\notag
	\end{equation}
	By the Aubin-Lions Lemma, there exists a function $w\in L^{\infty}(Q_T)\cap C([0,T],L^q(\Omega))$ for any $1\leq q<\frac{N}{N-s}$, with $\frac{\partial w}{\partial t}\in L^2(Q_T)$, and
	a subsequence of $\{\rho_n\}_{n=1}^{\infty}$ (not relabeled), such that as $n\rightarrow\infty$,
	\begin{equation}
		w_{\rho_n}\rightarrow w \text{ strongly in } C([0,T];L^1(\Omega)),~w_{\rho_n}\rightharpoonup w \text{ weakly }{\ast}\text{ in }L^{\infty}(Q_T),~\frac{\partial w_{\rho_n}}{\partial t}\rightharpoonup \frac{\partial w}{\partial t}\text{ weakly in }L^{2}(Q_T).\notag
	\end{equation}
	From Lemma \ref{wcp} and Remark \ref{rmk1}, we can find a $W^{1,\gamma}(Q_T)$-gradient Young measure $\nu=(\nu_{t,x})_{(t,x)\in Q_T}$ that satisfies
	\begin{align}
		&\langle\nu^{\rho_n},\mathrm{id}\rangle\rightharpoonup\langle\nu,\mathrm{id}\rangle\text{ weakly }{\ast}\text{ in }L^{\infty}(0,T;L^{\gamma}(\Omega)),~\langle\nu^{\rho_n},\vec{q}_{\gamma}\rangle\rightharpoonup\langle\nu,\vec{q}_{\gamma}\rangle\text{ weakly }{\ast}\text{ in }L^{\infty}(0,T;L^{\frac{\gamma}{\gamma-1}}(\Omega))\notag\\
		&\langle\nu^{\rho_n},\vec{q}_{\gamma}\cdot\mathrm{id} \rangle\rightharpoonup\langle\nu,\vec{q}_{\gamma}\cdot\mathrm{id}\rangle\text{ in the biting sense}.\notag
	\end{align}
	It follows from H\"{o}lder's inequality that
	\begin{align}
		\sup_{0\leq t\leq T}\left[w_{\rho}(t,\cdot)\right]_{W^{s,1}(\Omega)}=&\sup_{0\leq t\leq T}\int_{\Omega}\int_{\Omega}\frac{1}{\left|x-y\right|^{N+s}}\left|w_{\rho}(t,x)-w_{\rho}(t,y)\right|dxdy\notag\\
		\leq & \sup_{0\leq t\leq T}\left(\int_{\Omega}\int_{\Omega} \frac{1}{\left|x-y\right|^{N+s_{\rho} \rho}}\left|w_{\rho}(t,x)-w_{\rho}(t,y)\right|^{\rho} dxdy\right)^{\frac{1}{\rho}}\left|\Omega\times\Omega\right|^{\frac{\rho-1}{\rho}} \notag\\
		\leq & M \left(\left|\Omega\right|^2+1\right). \notag
	\end{align}
	Therefore, by the strong convergence of $w_{\rho_n}(t,\cdot)$ to $ w(t,\cdot)$ in $L^{1}(\Omega)$ and Fatou's Lemma, we conclude that $w\in L^{\infty}(0,T;W^{s,1}(\Omega))$. 
	
	Next, we prove that $(u,w)$ with $\nu$ is a Young measure solution to problem \eqref{New system} for $p=1$. For any $1<\sigma<\frac{\rho}{\rho-1}$, using H\"{o}lder's inequality, we have
	\begin{align}
		&\left\|\left|\frac{w_{\rho}(t,y)-w_{\rho}(t,x)}{\left|x-y\right|^{N+s}}\right|^{\rho-2}\frac{w_{\rho}(t,y)-w_{\rho}(t,x)}{\left|x-y\right|^{N+s}} \right\|^{\sigma}_{L^{\sigma}(Q_T^{1,2})}\notag\\
		\leq & \left\|\left|\frac{w_{\rho}(t,y)-w_{\rho}(t,x)}{\left|x-y\right|^{N+s}}\right|^{\rho-2}\frac{w_{\rho}(t,y)-w_{\rho}(t,x)}{\left|x-y\right|^{N+s}} \right\|^{\sigma}_{L^{\frac{\rho}{\rho-1}}(Q_T^{1,2})}\left|Q_T^{1,2}\right|^{1-\frac{\rho-1}{\rho}\sigma} \notag\\
		= & \left(\iiint_{Q_T^{1,2}} \frac{1}{\left|x-y\right|^{N+s_{\rho} \rho}}\left|w_{\rho}(t,x)-w_{\rho}(t,y)\right|^{\rho} dxdydt\right)^{\frac{\rho-1}{\rho}\sigma}\left|Q_T^{1,2}\right|^{1-\frac{\rho-1}{\rho}\sigma} \notag\\
		\leq &  \left(T\left\|w_{\rho}\right\|^{\rho}_{L^{\infty}(0,T;W^{s_{\rho},\rho}(\Omega))}\right)^{\frac{\rho-1}{\rho}\sigma}\left|Q_T^{1,2}\right|^{1-\frac{\rho-1}{\rho}\sigma} \notag\\
		\leq & C,\notag
	\end{align}
	where $C>0$ is a constant depending only on $T$, $\Omega$, and $M$. By a diagonal argument, we can extract a subsequence of $\{\rho_n\}_{n=1}^{\infty}$ (still denoted by itself) such that
	\begin{equation}
		\left|\frac{w_{\rho_n}(t,y)-w_{\rho_n}(t,x)}{\left|x-y\right|^{N+s}}\right|^{\rho_n-2}\frac{w_{\rho_n}(t,y)-w_{\rho_n}(t,x)}{\left|x-y\right|^{N+s}}\rightharpoonup \varsigma(t,x,y) \text{ weakly in }L^{\sigma}(Q_T^{1,2}),\notag
	\end{equation}
	where $\varsigma(t,x,y)=-\varsigma(t,y,x)$, and $\varsigma\in L^{\sigma}(Q_T^{1,2})$ for all $1<\sigma<+\infty$. Moreover,
	\begin{align}
		\left\|\varsigma\right\|_{L^{\sigma}(Q_T^{1,2})}^{\sigma}\leq\liminf_{n\rightarrow \infty}\left\|\left|\frac{w_{\rho_n}(t,y)-w_{\rho_n}(t,x)}{\left|x-y\right|^{N+s}}\right|^{\rho_n-2}\frac{w_{\rho_n}(t,y)-w_{\rho_n}(t,x)}{\left|x-y\right|^{N+s}} \right\|^{\sigma}_{L^{\sigma}(Q_T^{1,2})}\leq C.\notag
	\end{align}
	Then we have $\varsigma\in L^{\infty}(Q_T^{1,2})$, and  $\left\|\varsigma\right\|_{L^{\sigma}(Q_T^{1,2})}\leq C^{\frac{1}{\sigma}}$. Letting $\sigma\rightarrow +\infty$, we derive $\left\|\varsigma\right\|_{L^{\infty}(Q_T^{1,2})}\leq 1$. Now we pass to the limit as $\rho\rightarrow 1^{+}$ in \eqref{equation-w_p}. To this end, we choose the test function  $\psi(t,x) = \eta(t)\phi(x)$, where $\eta\in C^{\infty}_0(0,T)$ and $\phi\in C^{\infty}(\overline{\Omega})$. Since
	\begin{align}
		&\frac{1}{2}\iiint_{Q_T^{1,2}}\frac{1}{\left|x-y\right|^{N+s}}\varsigma(t,x,y)\left(\psi(t,y)-\psi(t,x)\right)dxdydt\notag\\
		=& \frac{1}{2}\lim_{n\rightarrow \infty}\iiint_{Q_T^{1,2}}\frac{1}{\left|x-y\right|^{N+s}}\left|\frac{w_{\rho_n}(t,y)-w_{\rho_n}(t,x)}{\left|x-y\right|^{N+s}}\right|^{\rho_n-2}\frac{w_{\rho_n}(t,y)-w_{\rho_n}(t,x)}{\left|x-y\right|^{N+s}}\left(\psi(t,y)-\psi(t,x)\right)dxdydt \notag\\
		=& -\lim_{n\rightarrow \infty}\iint_{Q_T}\frac{\partial w_{\rho_n}}{\partial t}\psi dxdt+\lambda_3\lim_{n\rightarrow \infty} \iint_{Q_T}\left(f-Ku_{\rho_n}\right)K\psi dxdt\notag\\
		=&-\iint_{Q_T}\frac{\partial w}{\partial t}\psi dxdt+\lambda_3 \iint_{Q_T}\left(f-Ku\right)K\psi dxdt.\notag
	\end{align}
	By the density of $C^{\infty}(\overline{\Omega})$ in $W^{s,1}(\Omega)$, it follows
	that the identity
	\begin{align}
		\iint_{Q_T}\frac{\partial w}{\partial t}\psi dxdt+\frac{1}{2}\iiint_{Q_T^{1,2}}\frac{1}{\left|x-y\right|^{N+s}}\varsigma(t,x,y)\left(\psi(t,y)-\psi(t,x)\right)dxdydt=\lambda_3 \iint_{Q_T}\left(f-Ku\right)K\psi dxdt\notag
	\end{align}
	holds for every $\psi\in L^2(Q_T)\cap L^1(0,T;W^{s,1}(\Omega))$. Taking $\psi=w_{\rho_n}$ as the test function in \eqref{equation-w_p} for $\rho_n$, and $\psi=w$ in \eqref{equation-w-1}, respectively, we have
	\begin{align}\notag
		&\frac{1}{2}\iiint_{Q_T^{1,2}}\frac{1}{\left|x-y\right|^{(N+s) \rho_n}}\left|w_{\rho_n}(t,y)-w_{\rho_n}(t,x)\right|^{\rho_n}dx dy dt\notag\\
		=&-\iint_{Q_T} \frac{\partial w_{\rho_n}}{\partial t}w_{\rho_n} dx dt+\lambda_3 \iint_{Q_T}\left(f-K u_{\rho_n}\right)Kw_{\rho_n} dx dt\notag\\
		=&-\iint_{Q_T} \left(\frac{\partial w_{\rho_n}}{\partial t}w_{\rho_n} - \frac{\partial w}{\partial t}w \right)+\lambda_3\iint_{Q_T}f\left(Kw_{\rho_n}-Kw\right)dxdt+\lambda_3\iint_{Q_T}\left(KuKw-Ku_{\rho_n}Kw_{\rho_n}\right)dxdt\notag\\
		&+\frac{1}{2}\iiint_{Q_T^{1,2}}\frac{1}{\left|x-y\right|^{N+s}}\varsigma(t,x,y)\left(w(t,y)-w(t,x)\right)dxdydt. \notag
	\end{align}
	Sending $n\rightarrow\infty$, we arrive at
	\begin{align}
		&\limsup_{n\rightarrow\infty}\frac{1}{2}\iiint_{Q_T^{1,2}}\frac{1}{\left|x-y\right|^{(N+s) \rho_n}}\left|w_{\rho_n}(t,y)-w_{\rho_n}(t,x)\right|^{\rho_n}dx dy dt\notag\\
		\leq & \frac{1}{2}\iiint_{Q_T^{1,2}}\frac{1}{\left|x-y\right|^{N+s}}\varsigma(t,x,y)\left(w(t,y)-w(t,x)\right)dxdydt. \notag
	\end{align}
	By Fatou's Lemma and H\"{o}lder's inequality, we deduce that
	\begin{align}
		&\frac{1}{2}\iiint_{Q_T^{1,2}}\frac{1}{\left|x-y\right|^{N+s}}\left|w(t,x)-w(t,y)\right|dxdydt\notag\\
		\leq & \liminf_{n\rightarrow \infty}\frac{1}{2}\iiint_{Q_T^{1,2}}\frac{1}{\left|x-y\right|^{N+s}}\left|w_{\rho_n}(t,x)-w_{\rho_n}(t,y)\right|dxdydt\notag\\
		\leq & \liminf_{n\rightarrow \infty}\frac{1}{2}\left(\iiint_{Q_T^{1,2}}\frac{1}{\left|x-y\right|^{(N+s)\rho_n}}\left|w_{\rho_n}(t,x)-w_{\rho_n}(t,y)\right|^{\rho_n}dxdydt\right)^{\frac{1}{\rho_n}}\left|Q_T^{1,2}\right|^{\frac{\rho_n-1}{\rho_n}}\notag\\
		\leq & \frac{1}{2}\iiint_{Q_T^{1,2}}\frac{1}{\left|x-y\right|^{N+s}}\varsigma(t,x,y)\left(w(t,y)-w(t,x)\right)dxdydt,\notag
	\end{align}
	which implies that 
	\begin{equation}
		\varsigma(t,x,y)\in\operatorname{sign}\left(w(t,y)-w(t,x)\right),~\text{a.e. }(t,x,y)\in Q_T^{1,2}.\notag
	\end{equation}
	
	On the other hand, taking the limit as $n\rightarrow\infty$ in \eqref{equation-u_p}, it is straightforward to deduce that identity \eqref{equation-for-u-nu} holds for all $\zeta\in L^2(Q_T)\cap L^1(0,T;W^{1,\gamma}(\Omega))$. By the weak convergence of $\nabla u_{\rho_n}=\langle\nu^{\rho_n},\mathrm{id}\rangle$, we obtain \eqref{gradientforu}, while \eqref{suppfornu} follows from the weak $\ast$ compactness of $\nu^{\rho_n}$ in $L^{\infty}(0,T;\mathcal{M}(\mathbb{R}^N))$. To complete the proof, it remains to verify the independence property
	\begin{equation}
		\langle\nu_{t,x},\vec{q}_{\gamma}\cdot\mathrm{id}\rangle=\langle\nu_{t,x},\vec{q}_{\gamma}\rangle\cdot\langle\nu_{t,x},\mathrm{id}\rangle,~\text{a.e. }(t,x)\in Q_T.\notag
	\end{equation}
	Following the approach used in the proof of Theorem \ref{thm2}, we substitute $\zeta = u_{\rho_n}\tilde{\zeta}$, with $\tilde{\zeta}\in C_0^{\infty}(Q_T)$, into \eqref{equation-u_p} for $\rho_n$, and into \eqref{equation-for-u-nu}, respectively, to obtain
	\begin{align}
		&\left|\iint_{Q_T}\left(\langle\nu^{\rho_n},\vec{q}_{\gamma}\rangle\cdot\langle\nu^{\rho_n},\mathrm{id}\rangle-\langle\nu,\vec{q}_{\gamma}\rangle\cdot\langle\nu,\mathrm{id}\rangle\right)\tilde{\zeta} dxdt\right|\notag\\
		\leq&\left| \iint_{Q_T}\left( \langle \nu^{\rho_n},\vec{q}_{\gamma} \rangle - \langle \nu,\vec{q}_{\gamma} \rangle \right)\cdot\langle \nu^{\rho_n},\mathrm{id} \rangle \tilde{\zeta} dxdt \right|+\left| \iint_{Q_T}\left( \langle \nu^{\rho_n},\mathrm{id} \rangle - \langle \nu,\mathrm{id} \rangle \right)\cdot\langle \nu,\vec{q}_{\gamma} \rangle \tilde{\zeta} dxdt \right|\notag\\
		\leq & \left| \iint_{Q_T}\left( \langle \nu^{\rho_n},\vec{q}_{\gamma} \rangle - \langle \nu,\vec{q}_{\gamma} \rangle \right)\cdot\nabla\left(u_{\rho_n} \tilde{\zeta}\right) dxdt \right|+\left| \iint_{Q_T}u_{\rho_n}\left( \langle \nu^{\rho_n},\vec{q}_{\gamma} \rangle - \langle \nu,\vec{q}_{\gamma} \rangle \right)\cdot\nabla\tilde{\zeta} dxdt \right| \notag\\
		&+\left| \iint_{Q_T}\left( \langle \nu^{\rho_n},\mathrm{id} \rangle - \langle \nu,\mathrm{id} \rangle \right)\cdot\langle \nu,\vec{q}_{\gamma} \rangle \tilde{\zeta} dxdt \right|\notag\\
		\leq & \left|\iint_{Q_T}u_{\rho_n}\left(\frac{\partial u}{\partial t}-\frac{\partial u_{\rho_n}}{\partial t}\right)\tilde{\zeta}dxdt\right|+\lambda_1\left|\iint_{Q_T}u_{\rho_n}\left(u-u_{\rho_n}\right)\tilde{\zeta} dxdt\right| +\lambda_2 \left|\iint_{Q_T}u_{\rho_n}\left(w_{\rho_n}-w\right)\tilde{\zeta} dxdt\right| \notag\\
		& + \left| \iint_{Q_T}u_{\rho_n}\left( \langle \nu^{\rho_n},\vec{q}_{\gamma} \rangle - \langle \nu,\vec{q}_{\gamma} \rangle \right)\cdot\nabla\tilde{\zeta} dxdt \right|+\left| \iint_{Q_T}\left( \langle \nu^{\rho_n},\mathrm{id} \rangle - \langle \nu,\mathrm{id} \rangle \right)\cdot\langle \nu,\vec{q}_{\gamma} \rangle \tilde{\zeta} dxdt \right| \notag\\
		\leq & \left|\iint_{Q_T}u\left(\frac{\partial u}{\partial t}-\frac{\partial u_{\rho_n}}{\partial t}\right)\tilde{\zeta}dxdt\right| + \left|\iint_{Q_T}\left(u_{\rho_n}-u\right)\left(\frac{\partial u}{\partial t}-\frac{\partial u_{\rho_n}}{\partial t}\right)\tilde{\zeta}dxdt\right| +\lambda_1\left|\iint_{Q_T}u_{\rho_n}\left(u-u_{\rho_n}\right)\tilde{\zeta} dxdt\right| \notag\\
		& + \lambda_2 \left|\iint_{Q_T}u_{\rho_n}\left(w_{\rho_n}-w\right)\tilde{\zeta} dxdt\right|+\left| \iint_{Q_T}u\left( \langle \nu^{\rho_n},\vec{q}_{\gamma} \rangle - \langle \nu,\vec{q}_{\gamma} \rangle \right)\cdot\nabla\tilde{\zeta} dxdt \right| \notag\\
		& + \left| \iint_{Q_T}\left(u_{\rho_n}-u\right)\left( \langle \nu^{\rho_n},\vec{q}_{\gamma} \rangle - \langle \nu,\vec{q}_{\gamma} \rangle \right)\cdot\nabla\tilde{\zeta} dxdt \right|+ \left| \iint_{Q_T}\left( \langle \nu^{\rho_n},\mathrm{id} \rangle - \langle \nu,\mathrm{id} \rangle \right)\cdot\langle \nu,\vec{q}_{\gamma} \rangle \tilde{\zeta} dxdt \right|.\notag
	\end{align}
	Letting $n\rightarrow\infty$, we have 
	\begin{equation}
		\lim_{n\rightarrow \infty}\left|\iint_{Q_T}\left(\langle\nu^{\rho_n},\vec{q}_{\gamma}\rangle\cdot\langle\nu^{\rho_n},\mathrm{id}\rangle-\langle\nu,\vec{q}_{\gamma}\rangle\cdot\langle\nu,\mathrm{id}\rangle\right)\tilde{\zeta} dxdt\right|=0.\notag
	\end{equation}
	Thus, $\langle\nu^{\rho_n},\vec{q}_{\gamma}\rangle\cdot\langle\nu^{\rho_n},\mathrm{id}\rangle$ converges weakly to $\langle\nu,\vec{q}_{\gamma}\rangle\cdot\langle\nu,\mathrm{id}\rangle$ in $L^1(Q_T)$. Combining the biting convergence of $\langle\nu^{\rho_n},\vec{q}_{\gamma}\cdot\mathrm{id}\rangle$ with the fact that  $\vec{q}_{\gamma}\cdot\mathrm{id}\in \mathcal{E}_0^{\gamma}(\mathbb{R}^N)$, we obtain \eqref{independentfornu}. Therefore, $(u,w)$ with $\nu$ is a Young measure solution to problem \eqref{New system} for $p=1$. The uniqueness of the solution follows from the monotonicity of $\operatorname{sign}(r)$ and the independence property \eqref{independentfornu}, as established in Theorem \ref{thm2}, and hence the details are omitted.
	
	Moreover, the weak $\ast$ $L^{\infty}(0,T;L^{\frac{\gamma}{\gamma-1}}(\Omega))$ convergence of $\vec{\varrho}_{\gamma}(\nabla u_{\rho_n})$ to $\vec{\varrho}_{\gamma}(\nabla u)$ can be derived by the standard monotone method, while the convergence of $\langle\nu^{\rho_n},\vec{\varrho}_{\gamma}\rangle$ and $\langle\nu^{\rho_n},\vec{q}_{\gamma}\rangle$ relies on Lemma \ref{wcp}. Thus, from \eqref{add3}, we obtain $\vec{\varrho}_{\gamma}(\nabla u(t,x))=\langle\nu_{t,x},\vec{\varrho}_{\gamma}\rangle=\langle\nu_{t,x},\vec{q}_{\gamma}\rangle$ for a.e. $(t,x)\in Q_T$. Similar to the proof of Corollary \ref{cor1}, we can further deduce that the pair $(u,w)$ is the unique classical weak solution to relaxed problem \eqref{relaxed problem p>1} associated with problem \eqref{New system} for $p=1$.

\end{pf}

\begin{rmk}
	For $\delta>0$, the well-posedness of system \eqref{New system} is established in Theorem \ref{theorem 1}. Indeed, at least for $\delta < \frac{1}{8}$, $\varphi_{\gamma}$ is of the convex-concave-convex type, and problem \eqref{New system} has no classical weak solution but only the Young measure solution described in Definition \ref{defn solution system}. However, as mentioned in \cite{GUIDOTTI2012backward,GUIDOTTI2013restoration}, there exists a positive constant $\delta_0=\delta_0(\gamma)$ such that for $\delta \geq\delta_0$, $\varphi_{\gamma}$ becomes convex, that is, its convexification satisfies $\varphi_{\gamma}^{**} = \varphi_{\gamma}$, where
	\begin{equation}
		\delta_0=\frac{1}{\gamma-1}\kappa_{\gamma}^{\frac{2-\gamma}{2}}\frac{\kappa_{\gamma}-1}{(\kappa_{\gamma}+1)^2},~\kappa_\gamma=\frac{1}{\gamma}\left(3+\sqrt{\gamma^2-2\gamma+9}\right).\notag
	\end{equation}
	In this case, the heat flux $\vec{q}_{\gamma} = \vec{\varrho}_{\gamma}$, and the first equation in \eqref{New system} ceases to exhibit backward diffusion, with only forward diffusion remaining. Combining the proof of Corollary \ref{cor1} with Theorem \ref{thm3}, it is straightforward to conclude that under the condition $f\in L^{\infty}(\Omega)\cap W^{1,\gamma}(\Omega)$, problem \eqref{New system} has a unique weak solution for any $\delta\geq \delta_0$ and $1\leq p\leq 2$.

	%In this case, the backward diffusion region in the first equation of \eqref{New system} vanishes, and only forward diffusion remains.

	%and problem \eqref{New system} has no classical weak solution, admitting only the Young measure solution described in Definition \ref{defn solution system}. 

	%Next, we will further investigate the existence and uniqueness of weak solutions to problem \eqref{New system} for $\delta\geq\delta_0$. In particular, for $\gamma=2$, we have $\delta_0=\frac{1}{8}$.
\end{rmk}

\section{Numerical experiments}\label{sec4}
In this section, we provide several numerical examples to illustrate the effectiveness of the proposed model for image restoration. The difference explicit scheme is implemented to solve problem \eqref{New system} in $\mathbb{R}^2$. We assume the observed image of size $I \times J$ pixels, $\Delta t$ the time step size, and $h$ the spatial grid size. Then the spatio-temporal discretization is given by
\begin{align}
	& t_n = n\Delta t,~n=0,1,2,\cdots,\notag\\
	& x_i = ih,~y_j = jh,~i=0,1,\cdots,I,~j=0,1,\cdots,J.\notag
\end{align}
Denote by $w_{i,j}^n$ the approximation of $w(n\Delta t, ih, jh)$. The discrete fractional $p$-Laplacian operator $\Delta_p^s$ is defined by
\begin{equation}
	\Delta_p^s w_{i,j}^n = \sum_{\stackrel{1\leq \tilde{i}\leq I,1\leq\tilde{j}\leq J}{\tilde{i}\neq i,\tilde{j}\neq j}}\frac{\operatorname{sign}_0\left(w_{\tilde{i},\tilde{j}}^n-w_{i,j}^n\right)}{\left(h\sqrt{(\tilde{i}-i)^2+(\tilde{j}-j)^2}\right)^{2+sp}}\left|w_{\tilde{i},\tilde{j}}^n-w_{i,j}^n\right|^{p-1},\notag
\end{equation}
where the univalued function $\operatorname{sign}_0(r)=\frac{r}{|r|}$ for $r\neq 0$ and $\operatorname{sign}_0(0)=0$. We also write $u_{i,j}^n$ for the approximation of $u(n\Delta t, ih, jh)$ and give some additional notations and assumptions for the numerical scheme that follows.
\begin{align}
	&\nabla_{\pm}^{x}u_{i,j}^n:=\pm\frac{u_{i\pm 1,j}^n-u_{i,j}^n}{h},~\nabla_{\pm}^{y}u_{i,j}^n:=\pm\frac{u_{i,j\pm 1}^n-u_{i,j}^n}{h},\notag\\
	&m(a,b)=\frac{1}{2}\left(\operatorname{sign}_0(a)+\operatorname{sign}_0(b)\right)\min\left(\left|a\right|,\left|b\right|\right),\notag\\
	&d_{i,j}^{x,n} = \left(\nabla_{+}^{x}u_{i,j}^n\right)^2+\epsilon\left(\nabla_{+}^{y}u_{i,j}^n\right)^2+\left(1-\epsilon\right)\left(m\left(\nabla_{+}^{y} u_{i,j}^n,\nabla_{-}^{y} u_{i,j}^n\right)\right)^2,\notag\\
	&d_{i,j}^{y,n} = \left(\nabla_{+}^{y}u_{i,j}^n\right)^2+\epsilon\left(\nabla_{+}^{x}u_{i,j}^n\right)^2+\left(1-\epsilon\right)\left(m\left(\nabla_{+}^{x} u_{i,j}^n,\nabla_{-}^{x} u_{i,j}^n\right)\right)^2,\notag\\
	&q_{i,j}^{x,n} = \frac{\nabla_{+}^x u^n_{i,j}}{1+\left(\nabla_{+}^{x}u_{i,j}^n\right)^2+\left(\nabla_{+}^{y}u_{i,j}^n\right)^2}+\delta\left(e+d_{i,j}^{x,n}\right)^{\frac{\gamma}{2}-1}\nabla_{+}^x u^n_{i,j},\notag\\
	&q_{i,j}^{y,n} = \frac{\nabla_{+}^y u^n_{i,j}}{1+\left(\nabla_{+}^{x}u_{i,j}^n\right)^2+\left(\nabla_{+}^{y}u_{i,j}^n\right)^2}+\delta\left(e+d_{i,j}^{y,n}\right)^{\frac{\gamma}{2}-1}\nabla_{+}^y u^n_{i,j},\notag
\end{align}
where $\epsilon \in \left\{0, 1\right\}$, and $e \ll 1$ is a small constant introduced to prevent division by zero. Furthermore, let $\mathcal{K}=\left(k_{i,j}\right)$ denote the discrete convolution kernel, with adjoint $\mathcal{K}^{\prime}$, derived from the continuous kernel $k$ via $k_{i,j} = k(x_i, y_j)$. The numerical approximation for problem \eqref{New system} is then expressed as follows:
\begin{align}
	&w^{n+1}=w^{n}+\Delta t\Delta_p^s w^n-\lambda_3\Delta  t \mathcal{K}^{\prime}\ast\left(\mathcal{K}\ast u^n-u^0\right),\notag\\
	&u^{n+1}=u^{n}+\Delta t \left(\nabla_{-}^{x}q^{x,n}+\nabla_{-}^{y}q^{y,n}\right)-\lambda_1 \Delta t u^n+\lambda_2 \Delta t w^{n+1},\notag\\
	&u_{i,j}^0=f_{i,j}:=f(x_i,y_j),~w_{i,j}^0=0,~i=0,1,\cdots,I,~j=0,1,\cdots,J,\notag\\
	&u^n_{i,0} = u^n_{i,1},~u^n_{0,j} = u^n_{1,j},
	~u^n_{i,J} = u^n_{i,J-1},~u^n_{I,j} = u^n_{I-1,j}.\notag
\end{align}
Through the above lines, we can obtain $u^1$ and  $w^1$ from $u^0$ and $w^0$. The restoration quality is evaluated using the peak signal noise ratio (PSNR)
\begin{equation}
	\text{PSNR}\left(u,f\right)=10\log_{10}\left(\frac{\sum_{i,j} 255^2}{\sum_{i,j}(u_{i,j}-f_{i,j})^2}\right),\notag
\end{equation}
and the structural similarity index measure (SSIM)
\begin{equation}
	\text{SSIM}\left(u,f\right)=\frac{(2\mu_u \mu_f + c_1)(2\sigma_{uf} + c_2)}{\mu_u^2 + \mu_f^2 + c_1 + \sigma_u^2 + \sigma_f^2 + c_2},
	\notag
\end{equation}
where $f$ denotes the observed image and $u$ the compared (restored) image. The terms $\mu_u$, $\mu_f$, $\sigma_u$, $\sigma_f$, and $\sigma_{uf}$ denote the local means, standard deviations, and cross-covariance of images $u$ and $f$, respectively. The constants $c_1$ and $c_2$ are employed to stabilize the division when the denominator is close to zero. Higher PSNR and SSIM values indicate better image quality. The stopping criteria of our algorithm is designed through reaching the maximal PSNR.

\begin{figure}%[]
	\centering
	\subfigure[Texture 1]{
		\includegraphics[width=0.23\linewidth]{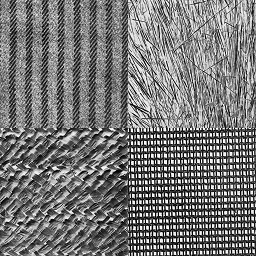}
	}
	\centering
	\subfigure[Texture 2]{
		\includegraphics[width=0.23\linewidth]{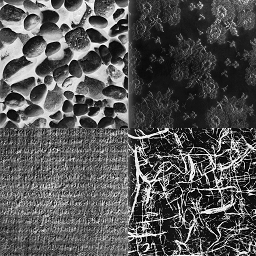}
	}
	\centering
	\subfigure[House]{
		\includegraphics[width=0.23\linewidth]{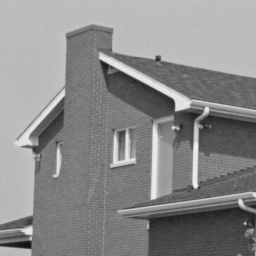}
	}
	\quad
	\centering
	\subfigure[Cameraman]{
		\includegraphics[width=0.23\linewidth]{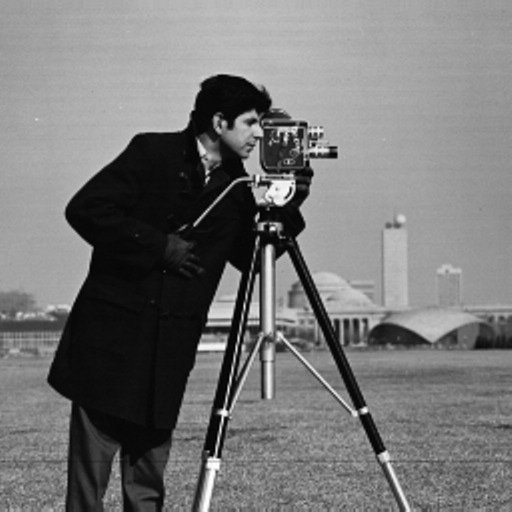}
	}
	\centering
	\subfigure[Satellite 1]{
		\includegraphics[width=0.23\linewidth]{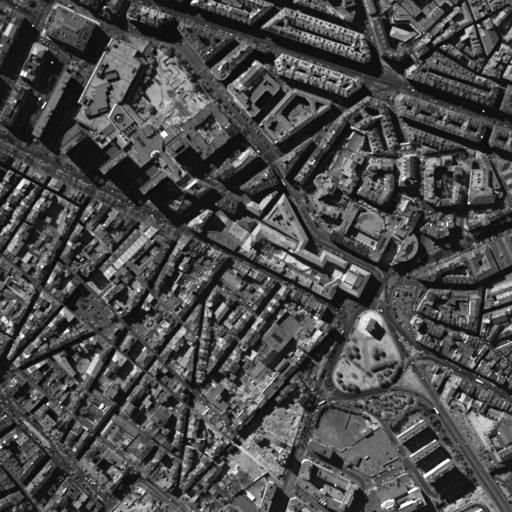}
	}
	\centering
	\subfigure[Satellite 2]{
		\includegraphics[width=0.23\linewidth]{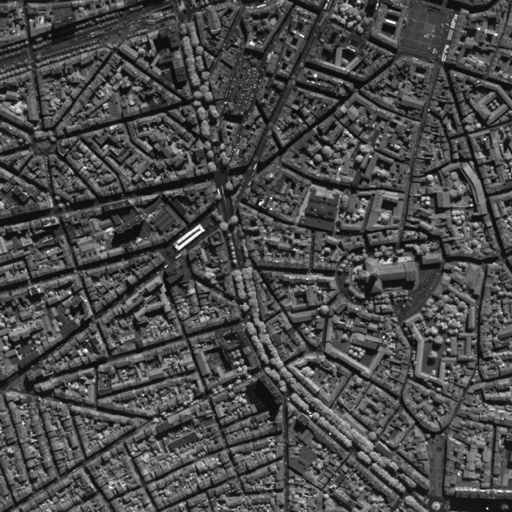}
	}
	\caption{Test images.}\label{fig1}
\end{figure}

In experiments, we consider the following three common types of blur kernels: the disk kernel with a radius of $3$, the $5 \times 5$ average kernel, and the motion blur kernel with length $20$ and angle $\frac{\pi}{3}$. To evaluate the effectiveness of the proposed model, we conduct tests on several noisy blurred images, simultaneously degraded by Gaussian white noise with a standard deviation $\sigma=3$ and one of the blur kernels described above. For illustration, the results for the
texture 1 ($256 \times 256$), texture 2 ($256 \times 256$), house ($256 \times 256$), cameraman ($512 \times 512$), satellite 1 ($512 \times 512$) and satellite 2 ($512 \times 512$) are presented, see the original images in Figure \ref{fig1}. The proposed model is compared with the NLTV method \cite{GILBOA2009}, the FBRPM model \cite{GUIDOTTI2013restoration}, and the nonlinear fractional diffusion (NFD) model \cite{GUO2019deblur}, all implemented with parameter settings recommended in the respective references. For the proposed algorithm, we set $\Delta t=0.01$, $\delta = 0.28$, $\gamma\in [1.28,1.38]$, and $s\in [0.74,0.87]$ in all experiments. The parameters $\lambda_1$, $\lambda_2$, and $\lambda_3$ are selected based on the specific image and the type of blur kernel. To quantitatively assess the restoration performance, the PSNR and SSIM values of the restored images are presented in Tables \ref{tbl1} and \ref{tbl2}.

\begin{figure}%[]
	\centering
	\subfigure[Original]{
		\includegraphics[width=0.23\linewidth]{figs/clean/texture}
	}
	\centering
	\subfigure[Blurred]{
		\includegraphics[width=0.23\linewidth]{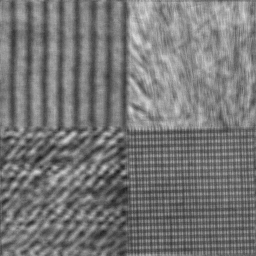}
	}
	\centering
	\subfigure[NLTV]{
		\includegraphics[width=0.23\linewidth]{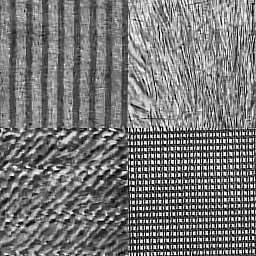}\label{texture1c}
	}
	\quad
	\centering
	\subfigure[NFD]{
		\includegraphics[width=0.23\linewidth]{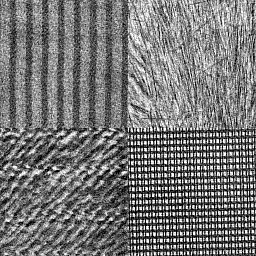}\label{texture1d}
	}
	\centering
	\subfigure[FBRPM]{
		\includegraphics[width=0.23\linewidth]{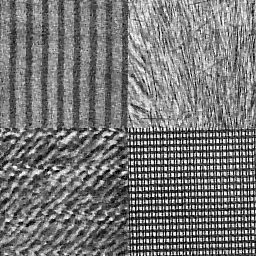}\label{texture1e}
	}
	\centering
	\subfigure[Ours]{
		\includegraphics[width=0.23\linewidth]{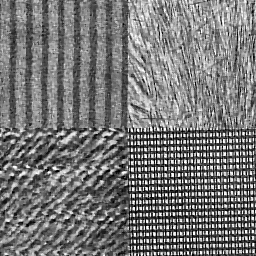}\label{texture1f}
	}
	\caption{Recovered results of NLTV, NFD, FBRPM and our model for the texture 1 image with average blur and corrupted by Gaussian noise of standard deviation $\sigma = 3$.}\label{texture1recover}
\end{figure}

\begin{figure}%[]
	\centering
	\subfigure[Original]{
		\includegraphics[width=0.23\linewidth]{figs/clean/texture2}
	}
	\centering
	\subfigure[Blurred]{
		\includegraphics[width=0.23\linewidth]{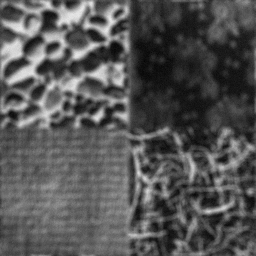}
	}
	\centering
	\subfigure[NLTV]{
		\includegraphics[width=0.23\linewidth]{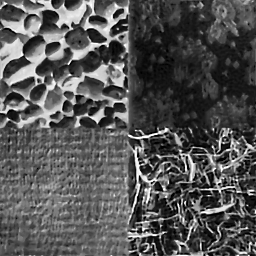}\label{texture2c}
	}
	\quad
	\centering
	\subfigure[NFD]{
		\includegraphics[width=0.23\linewidth]{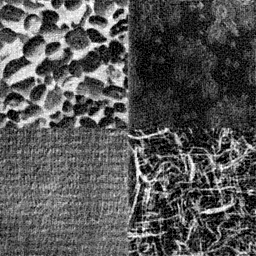}\label{texture2d}
	}
	\centering
	\subfigure[FBRPM]{
		\includegraphics[width=0.23\linewidth]{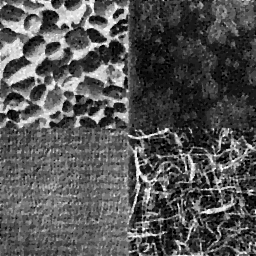}\label{texture2e}
	}
	\centering
	\subfigure[Ours]{
		\includegraphics[width=0.23\linewidth]{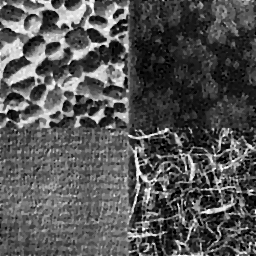}\label{texture2f}
	}
	\caption{Recovered results of NLTV, NFD, FBRPM and our model for the texture 2 image with disk blur and corrupted by Gaussian noise of standard deviation $\sigma = 3$.}\label{texture2recover}
\end{figure}

\begin{figure}%[]
	\centering
	\subfigure[Original]{
		\includegraphics[width=0.23\linewidth]{figs/clean/house}
	}
	\centering
	\subfigure[Blurred]{
		\includegraphics[width=0.23\linewidth]{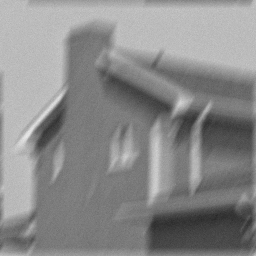}
	}
	\centering
	\subfigure[NLTV]{
		\includegraphics[width=0.23\linewidth]{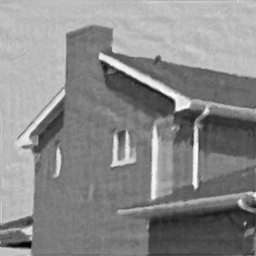}\label{housec}
	}
	\quad
	\centering
	\subfigure[NFD]{
		\includegraphics[width=0.23\linewidth]{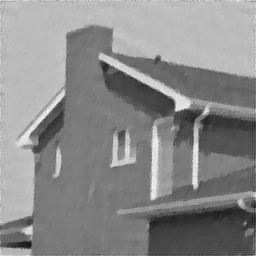}\label{housed}
	}
	\centering
	\subfigure[FBRPM]{
		\includegraphics[width=0.23\linewidth]{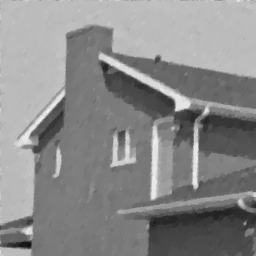}\label{housee}
	}
	\centering
	\subfigure[Ours]{
		\includegraphics[width=0.23\linewidth]{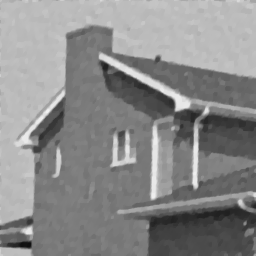}\label{housef}
	}
	\caption{Recovered results of NLTV, NFD, FBRPM and our model for the house image with motion blur and corrupted by Gaussian noise of standard deviation $\sigma = 3$.}\label{houserecover}
\end{figure}

\begin{figure}%[]
	\centering
	\subfigure[Original]{
		\includegraphics[width=0.23\linewidth]{figs/clean/Cameraman}
	}
	\centering
	\subfigure[Blurred]{
		\includegraphics[width=0.23\linewidth]{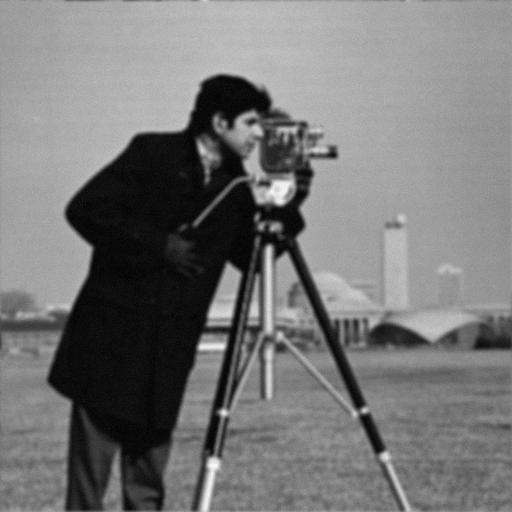}
	}
	\centering
	\subfigure[NLTV]{
		\includegraphics[width=0.23\linewidth]{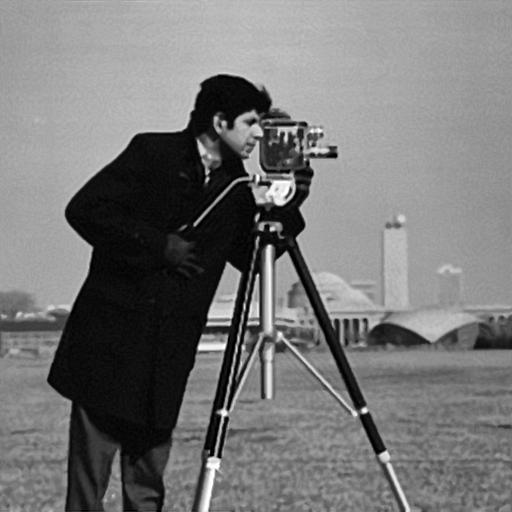}\label{cameramanc}
	}
	\quad
	\centering
	\subfigure[NFD]{
		\includegraphics[width=0.23\linewidth]{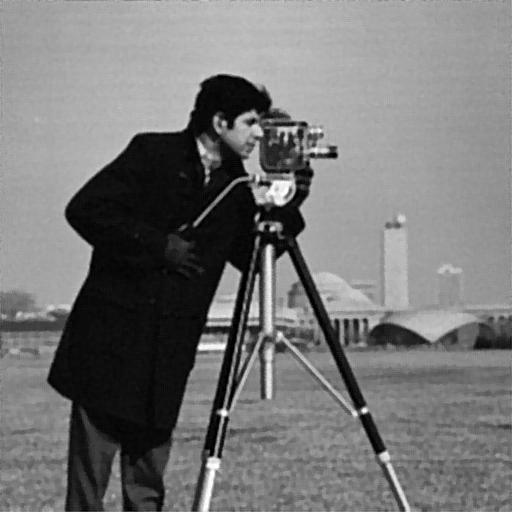}\label{cameramand}
	}
	\centering
	\subfigure[FBRPM]{
		\includegraphics[width=0.23\linewidth]{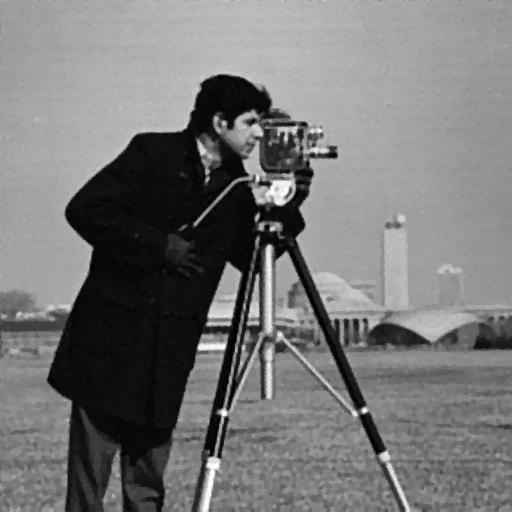}\label{cameramane}
	}
	\centering
	\subfigure[Ours]{
		\includegraphics[width=0.23\linewidth]{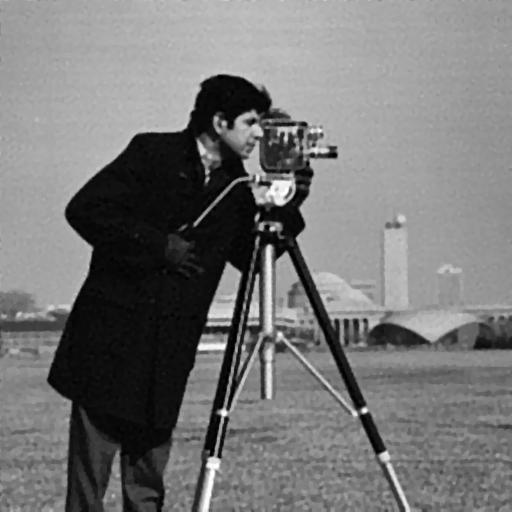}\label{cameramanf}
	}
	\caption{Recovered results of NLTV, NFD, FBRPM and our model for the cameraman image with disk blur and corrupted by Gaussian noise of standard deviation $\sigma = 3$.}\label{cameramanrecover}
\end{figure}

\begin{figure}%[]
	\centering
	\subfigure[Original]{
		\includegraphics[width=0.23\linewidth]{figs/clean/satellite1}
	}
	\centering
	\subfigure[Blurred]{
		\includegraphics[width=0.23\linewidth]{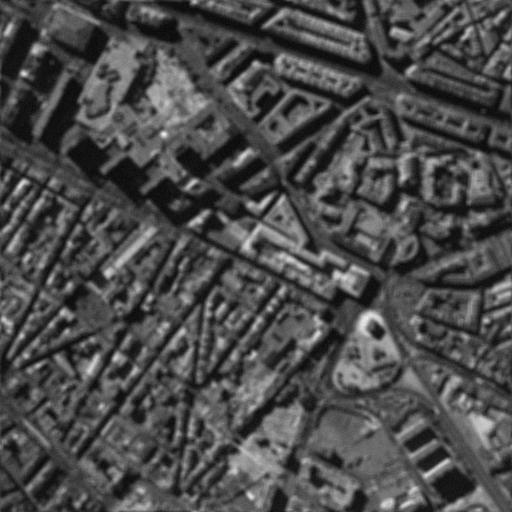}
	}
	\centering
	\subfigure[NLTV]{
		\includegraphics[width=0.23\linewidth]{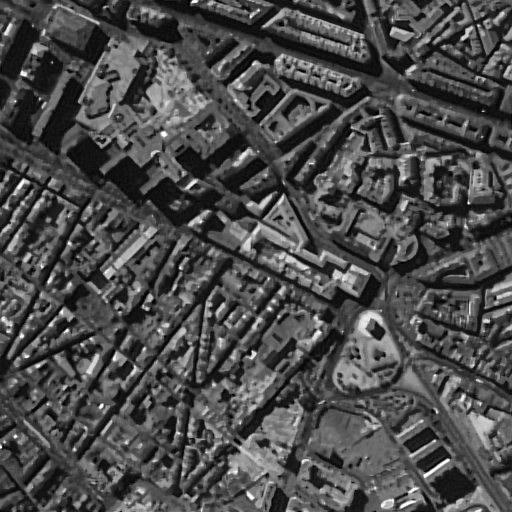}\label{satellite1c}
	}
	\quad
	\centering
	\subfigure[NFD]{
		\includegraphics[width=0.23\linewidth]{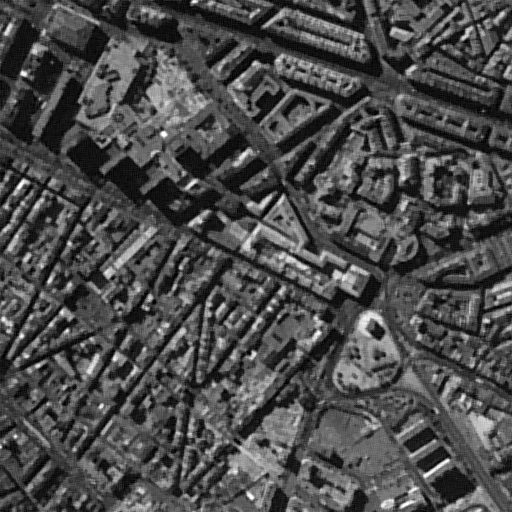}\label{satellite1d}
	}
	\centering
	\subfigure[FBRPM]{
		\includegraphics[width=0.23\linewidth]{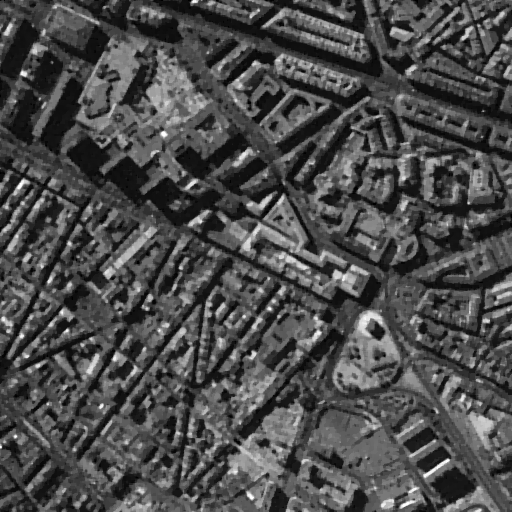}\label{satellite1e}
	}
	\centering
	\subfigure[Ours]{
		\includegraphics[width=0.23\linewidth]{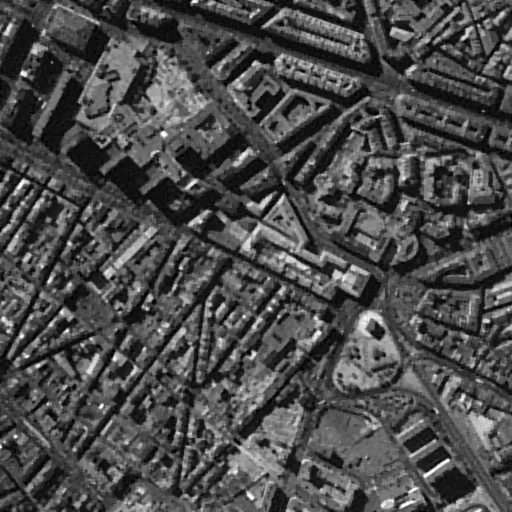}\label{satellite1f}
	}
	\caption{Recovered results of NLTV, NFD, FBRPM and our model for the satellite 1 image with average blur and corrupted by Gaussian noise of standard deviation $\sigma = 3$.}\label{satellite1recover}
\end{figure}

\begin{figure}%[]
	\centering
	\subfigure[Original]{
		\includegraphics[width=0.23\linewidth]{figs/clean/satellite2}
	}
	\centering
	\subfigure[Blurred]{
		\includegraphics[width=0.23\linewidth]{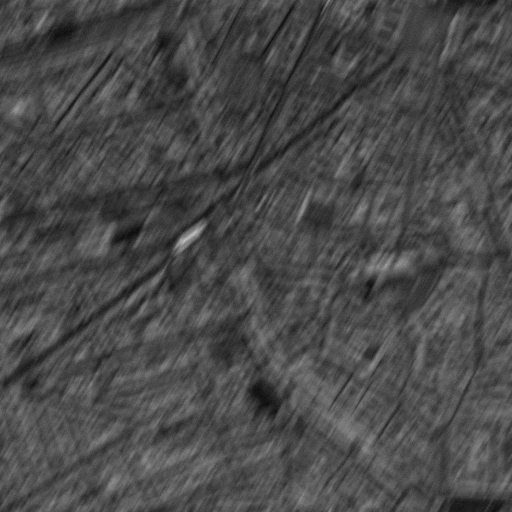}
	}
	\centering
	\subfigure[NLTV]{
		\includegraphics[width=0.23\linewidth]{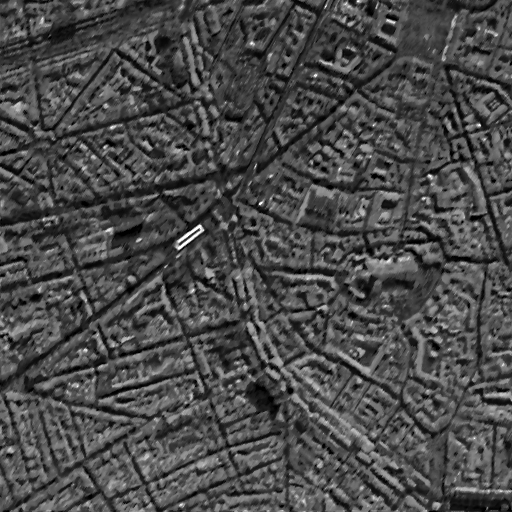}\label{satellite2c}
	}
	\quad
	\centering
	\subfigure[NFD]{
		\includegraphics[width=0.23\linewidth]{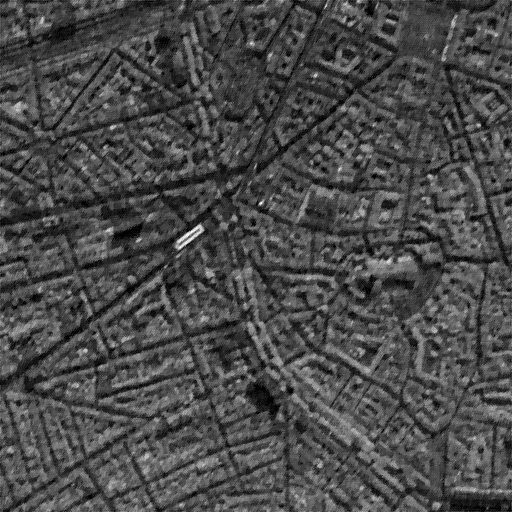}\label{satellite2d}
	}
	\centering
	\subfigure[GLRPM]{
		\includegraphics[width=0.23\linewidth]{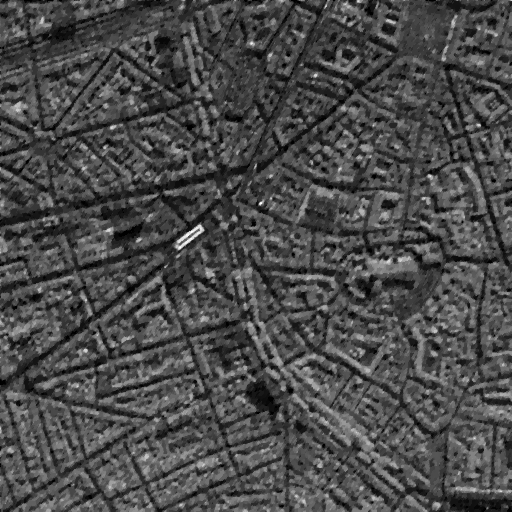}\label{satellite2e}
	}
	\centering
	\subfigure[Ours]{
		\includegraphics[width=0.23\linewidth]{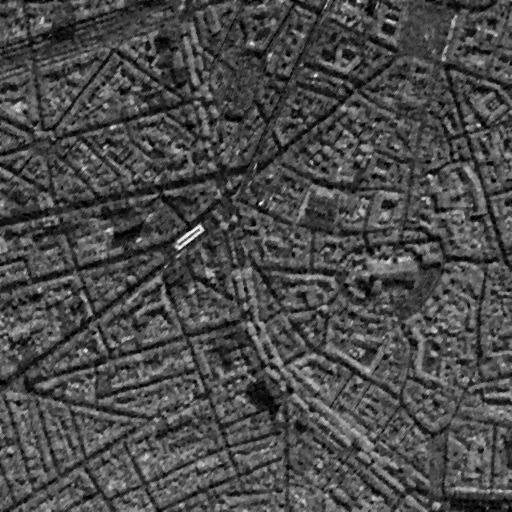}\label{satellite2f}
	}
	\caption{Recovered results of NLTV, NFD, GLRPM and our model for the satellite 2 image with motion blur and corrupted by Gaussian noise of standard deviation $\sigma = 3$.}\label{satellite2recover}
\end{figure}

We now report the numerical experiments on image deblurring and denoising for the test images shown in Figure \ref{fig1}. The corresponding results are illustrated in Figures \ref{texture1recover}-\ref{satellite2recover}. First, Figures \ref{texture1recover} and   \ref{texture2recover} compare the restored results of the texture 1 and texture 2 images, which are blurred by the average kernel and the disk kernel, respectively. We observe that Figures \ref{texture1c} and \ref{texture1d} better preserve texture features; however, the former exhibits a slightly higher noise level, whereas the latter shows reduced edge sharpness. Compared with Figure \ref{texture1e}, the proposed model retains more texture details, as shown in Figure \ref{texture1f}. Although Figure \ref{texture2c} achieves higher PSNR and SSIM values, it visually loses more edge and detail information  compared to Figures \ref{texture2e} and \ref{texture2f}. Regarding texture preservation, the restoration result of the proposed model outperforms that in Figure \ref{texture2d} and effectively suppresses noise.

The restoration results for the house and cameraman images are shown in Figures \ref{houserecover} and \ref{cameramanrecover}. The house image is blurred by the motion kernel, while the cameraman image is blurred by the disk kernel. As shown in Figures \ref{housec} and \ref{housed}, the restoration results produced by NLTV and NFD exhibit noticeable artifacts and fail to preserve certain texture features. Although the proposed model is slightly less effective than Figure \ref{housee} in homogeneous regions, it demonstrates superior texture preservation, see Figure \ref{housef}. Compared with Figure \ref{cameramanc}, Figures \ref{cameramane} and \ref{cameramanf} show slightly lower contrast in fine details, but they reduce noise more effectively in smooth areas and preserve clearer edge structures than Figure \ref{cameramand}. Finally, Figures \ref{satellite1recover} and   \ref{satellite2recover} compare the restored results of the satellite 1 and satellite 2 images, which are blurred by the average kernel and the motion kernel, respectively. In Figures \ref{satellite1c}, \ref{satellite1d}, \ref{satellite2c}, and \ref{satellite2d}, NLTV results in blurred edges, and the images restored by NFD contain higher noise levels. Overall, the proposed method demonstrates improved visual quality relative to FBRPM, while maintaining a favorable balance between image restoration and edge preservation, as shown in Figures \ref{satellite1f} and \ref{satellite2f}.

\begin{table}[h]
	\centering
	\caption{Comparison of PSNR and SSIM for different models in the experiments on texture 1, texture 2 and house image. Bold values indicate the best result.}\label{tbl1}
	\begin{tabular}{cccccccc}
		\toprule
		& \multicolumn{2}{c}{Figure \ref{texture1recover}} & \multicolumn{2}{c}{Figure \ref{texture2recover}} & \multicolumn{2}{c}{Figure \ref{houserecover}} \\
		\cmidrule(lr){2-3} \cmidrule(lr){4-5} \cmidrule(lr){6-7}
		& PSNR & SSIM & PSNR & SSIM & PSNR & SSIM \\
		\midrule
		Blurred & 12.73 & 0.2040 & 12.52 & 0.3511 & 21.61 & 0.5865 \\
		NLTV & 17.34 & 0.7172 & \textbf{18.63} & \textbf{0.6640} & 29.42 & 0.7700 \\
		NFD & 17.46 & \textbf{0.7291} & 18.23 & 0.6142 & 28.61 & 0.7848 \\
		FBRPM & 17.53 & 0.7188 & 18.53 & 0.6429 & 30.17 & \textbf{0.8217} \\
		Ours & \textbf{17.55} & 0.7121 & 18.55 & 0.6447 & \textbf{30.29} & 0.8132 \\
		\bottomrule
	\end{tabular}
\end{table}
\begin{table}[h]
	\centering
	\caption{Comparison of PSNR and SSIM for different models in the experiments on cameraman, satellite 1 and satellite 2 image. Bold values indicate the best result.}\label{tbl2}
	\begin{tabular}{cccccccc}
		\toprule
		& \multicolumn{2}{c}{Figure \ref{cameramanrecover}} & \multicolumn{2}{c}{Figure \ref{satellite1recover}} & \multicolumn{2}{c}{Figure \ref{satellite2recover}} \\
		\cmidrule(lr){2-3} \cmidrule(lr){4-5} \cmidrule(lr){6-7}
		& PSNR & SSIM & PSNR & SSIM & PSNR & SSIM \\
		\midrule
		Blurred & 27.63 & 0.9053 & 22.49 & 0.8547 & 17.55 & 0.2369 \\
		NLTV & 32.89 & 0.9423 & 26.44 & 0.9575 & 22.62 & 0.8485 \\
		NFD & 32.06 & 0.9530 & 25.71 & 0.9525 & 22.00 & 0.8230 \\
		FBRPM & 32.70 & \textbf{0.9565} & 26.33 & 0.9575 & 23.08 & 0.8656 \\
		Ours & \textbf{32.92} & 0.9510 & \textbf{26.55} & \textbf{0.9588} & \textbf{23.46} & \textbf{0.8771} \\
		\bottomrule
	\end{tabular}
\end{table}

\section{Conclusion}\label{}
A class of nonlinear reaction-diffusion systems with local and nonlocal coupling has been proposed for image restoration, inheriting the advantages of forward-backward diffusion in preserving high-frequency contour features and fractional diffusion in capturing texture details. The existence of Young measure solutions to the proposed model for $1<p\leq 2$ has been discussed through Rothe's method, the relaxation theorem, and Moser's iteration, while the uniqueness relies on the independence property satisfied by the solution. For the case $p=1$, the well-posedness of the model has been established by means of the regularization method. The explicit finite difference scheme has been adopted for numerical implementation. Through numerical experiments, we demonstrate that our model performs satisfactorily across the entire image, effectively balancing noise and blur removal with edge and texture preservation.

\section*{Acknowledgements}
This work is partially supported by the National Natural Science Foundation of China (12301536, U21B2075, 12171123, 12371419, 12271130), the Fundamental Research Funds for the Central Universities (HIT.NSRIF202302, 2022FRFK060020, 2022FRFK060031, 2022FRFK060014), the Natural Science Foundation of Heilongjiang Province (ZD2022A001).

\section*{\small Declaration of competing interest}

{\small
	The authors declare that they have no conflict of interest.}

%% Use figure environment to create figures
%% Refer following link for more details.
%% https://en.wikibooks.org/wiki/LaTeX/Floats,_Figures_and_Captions
%\begin{figure}[t]%% placement specifier
%% Use \includegraphics command to insert graphic files. Place graphics files in 
%% working directory.
%\centering%% For centre alignment of image.
%\includegraphics{example-image-a}
%% Use \caption command for figure caption and label.
%\caption{Figure Caption}\label{fig1}
%% https://en.wikibooks.org/wiki/LaTeX/Importing_Graphics#Importing_external_graphics
%\end{figure}

%% The Appendices part is started with the command \appendix;
%% appendix sections are then done as normal sections
%\appendix
%\section{Example Appendix Section}
%\label{app1}

%Appendix text.

%% For citations use: 
%%       \citet{<label>} ==> Lamport [21]
%%       \citep{<label>} ==> [21]
%%
%Example citation, See \citet{WELK2005}.

%% If you have bib database file and want bibtex to generate the
%% bibitems, please use
%%
%%  \bibliographystyle{elsarticle-num-names} 
%%  \bibliography{<your bibdatabase>}

%% else use the following coding to input the bibitems directly in the
%% TeX file.

%% Refer following link for more details about bibliography and citations.
%% https://en.wikibooks.org/wiki/LaTeX/Bibliography_Management

%\begin{thebibliography}{00}

%% For authoryear reference style
%% \bibitem[Author(year)]{label}
%% Text of bibliographic item

%\bibitem[Lamport(1994)]{lamport94}
%  Leslie Lamport,
%  \textit{\LaTeX: a document preparation system},
%  Addison Wesley, Massachusetts,
%  2nd edition,
%  1994.

%\end{thebibliography}
\bibliographystyle{elsarticle-names}
\bibliography{cas-refs}

\end{document}